\documentclass[12pt]{article} 

\usepackage{latexsym, amsmath, amscd, amsthm, graphicx, 
epsfig, amssymb}
\usepackage{xypic}
\usepackage{psfrag}
\DeclareGraphicsExtensions{.ps}
   
\newtheorem{thm}{Theorem}[section]
\newtheorem{defn}[thm]{Definition}
\newtheorem{prop}[thm]{Proposition}
\newtheorem{cor}[thm]{Corollary}
\newtheorem{rema}[thm]{Remark}
\newtheorem{lemma}[thm]{Lemma}
\newtheorem{ass}[thm]{Assumption}

\newcommand{\halmos}{\rule{1ex}{1.4ex}}

\newcommand{\beq}{\begin{equation}}
\newcommand{\eeq}{\end{equation}}
\newcommand{\bnu}{\begin{enumerate}}
\newcommand{\enu}{\end{enumerate}}
 \newcommand{\bea}{\begin{eqnarray}}
\newcommand{\eea}{\end{eqnarray}}
 \newcommand{\nn}{\nonumber \\}
        \newcommand{\nno}{\nonumber}
        \newcommand{\lbar}{\bigg\vert}

\renewcommand{\hom}{\mbox{\rm Hom}}

\newcommand{\edo}{\mbox{\rm End}\;}
 \newcommand{\pf}{{\it Proof.}\hspace{2ex}}
 \newcommand{\epf}{\hspace*{\fill}\mbox{$\halmos$}}

\newcommand{\wt}{\mbox{\rm wt}\ }

\newcommand{\one}{\mathbf{1}}

\newcommand{\ds}{\displaystyle}
\newcommand{\C}{\mathbb{C}}
 
\newcommand{\HH}{\mathbb{H}}
\newcommand{\N}{\mathbb{N}}
\newcommand{\Q}{\mathbb{Q}}
\newcommand{\R}{\mathbb{R}}
\newcommand{\Z}{\mathbb{Z}}
\newcommand{\Y}{\mathcal{Y}}

\newcommand{\A}{\mathcal{A}}

\newcommand{\I}{\mathcal{I}}
\newcommand{\V}{\mathcal{V}}

\title{ {\bf Full field algebras, operads and tensor categories} }
\author{Liang Kong}
\date{}

\begin{document}

\bibliographystyle{alpha}
\maketitle

\begin{abstract} 
We study the operadic and categorical 
formulations of (conformal) full field algebras.
In particular, we show that a grading-restricted 
$\R\times \R$-graded
full field algebra is equivalent to an algebra over a
partial operad constructed from spheres with punctures and
local coordinates. This result is generalized to 
conformal full field algebras over $V^L\otimes V^R$, 
where $V^L$ and $V^R$ are two vertex operator algebras
satisfying certain finiteness and reductivity conditions. 
We also study the geometry interpretation of 
conformal full field algebras over $V^L\otimes V^R$
equipped with a nondegenerate invariant bilinear form. 
By assuming slightly stronger conditions on $V^L$ and $V^R$, 
we show that a conformal 
full field algebra over $V^L\otimes V^R$ equipped with 
a nondegenerate invariant bilinear form exactly corresponds to 
a commutative Frobenius algebra with a trivial twist in the 
category of $V^L\otimes V^R$-modules. 
The so-called diagonal constructions \cite{HK2}
of conformal full field algebras are given in 
tensor-categorical language. 
\end{abstract}

\renewcommand{\theequation}{\thesection.\arabic{equation}}
\renewcommand{\thethm}{\thesection.\arabic{thm}}
\setcounter{equation}{0}
\setcounter{thm}{0}
\setcounter{section}{-1}

\section{Introduction}

In \cite{HK2}, Huang and the author 
introduced the notion of conformal 
full field algebra and some variants of this notion.
We also studied their basic properties and gave constructions.
We explained briefly without
giving details in the introduction of \cite{HK2} that
our goal is to construct conformal field theories 
\cite{BPZ}\cite{MS}. 
It is one of the purpose of this work to
explain the connection between conformal full field algebras and 
genus-zero conformal field theories. 
\cite{HK2} and this work are actually a part of 
Huang's program (\cite{H1}-\cite{H12}) 
of constructing rigorously conformal field theories 
in the sense of Kontsevich and Segal \cite{S1}\cite{S2}.

Around 1986, I. Frenkel initiated a program 
of using vertex operator
algebras to construct, in a suitable sense,  
geometric conformal field theories, 
the precise mathematical 
definition of which was actually given independently by 
Kontsevich and Segal \cite{S1}\cite{S2} in 1987. 
According to Kontsevich and Segal, 
a conformal field theory is a projective tensor functor
from a category consisting of 
finite many ordered copies of $S^1$ as objects and the
equivalent classes of Riemann surfaces with 
parametrized boundaries as morphisms to the category of 
Hilbert spaces. 
This beautiful and compact definition 
of conformal field theory encloses enormously rich structures 
of conformal field theory.
In order to construct such theories, 
it is more fruitful to look at some substructures of 
conformal field theories at first.
In the first category, the set of morphisms which have
arbitrary number of copies of $S^1$ 
in their domains and only one copy of $S^1$ in their
codomain has a structure of operad \cite{Ma}, 
which is induced by gluing surfaces
along their parametrized boundaries. 
We denoted this operad as $K_{\mathfrak{H}}$.  
Let $H$, a Hilbert space, be the image object of $S^1$. 
The projective tensor functor in the definition of conformal 
field theory endows $H$ with
a structure of algebra over an operad, 
which is a $\C$-extension of the operad $K_{\mathfrak{H}}$ \cite{H4}.

It is very difficult to study this algebra-over-operad
structure on $H$ directly. It was suggested by Huang 
\cite{H11}\cite{H12}
that one should first look at a dense subset of $H$, which
carries a structure of an algebra over a partial operad 
$\tilde{K}^{c}\otimes \overline{\tilde{K}^{\bar{c}}}$ for 
$c\in\C$. The so-called sphere partial operad (\cite{H4})
is denoted by $K$. It contains $K_{\mathfrak{H}}$ as a suboperad. 
$\tilde{K}^c$ is the 
$\frac{c}{2}$-power of the determine line bundle over the
partial operad $K$. 
The restriction of the line bundle $\tilde{K}^c$
on $K_{\mathfrak{H}}$, denoted as 
$\tilde{K}_{\mathfrak{H}}^c$, is also a suboperad. 
$\tilde{K}^{c}\otimes \overline{\tilde{K}^{\bar{c}}}$ is simply
the tensor product of the bundle $\tilde{K}^c$
and the complex conjugate of the bundle $\tilde{K}^{\bar{c}}$. 
It is shown by Huang in \cite{H1}\cite{H4} that the category of
algebras over partial operad $\tilde{K}^c$ 
(or $\tilde{K}^c$-algebras), satisfying some natural conditions 
and the condition that all the correlation 
functions are rational functions,
is isomorphic to the category of vertex operator algebras with 
central charge $c$. 
Importantly, Huang also showed in \cite{H11}\cite{H12} that one can 
obtain algebras over $\tilde{K}_{\mathfrak{H}}^c$, which is
a true operad instead of partial operad, 
by completing the space of vertex operator algebra properly.

The results of Huang suggests that vertex operator algebras
(or its proper completion) describe genus-zero conformal 
field theories whose correlation functions are
rational functions. 
It turns out that vertex operator algebras, in general, are not
enough for the genus-one properties of conformal field theories, 
which require that all the genus-one correlation functions
are modular invariant 
\cite{Z}\cite{DLM}\cite{Mi1}\cite{Mi2}\cite{H8}. 
It is well-known in physics that one needs to combine both the 
chiral theory \cite{H3}\cite{H5}\cite{H6} 
and the antichiral theory to obtain a full conformal field theory.  
As a result
the correlation functions in a conformal field theory, in general,
are neither holomorphic nor antiholomorphic. 
The conformal full field algebras introduced 
in \cite{HK2} are exactly
$\tilde{K}^{c^L}\otimes \overline{\tilde{K}^{\overline{c^R}}}$-algebras
with neither holomorphic nor antiholomorphic 
correlation functions in general, and are capable of
producing modular invariant genus-one correlation functions
when $c^L=c^R$.  
We show in Section 1 that the grading-restricted 
$\R\times \R$-graded
full field algebras are exactly smooth algebras over 
$\widehat{K}$, which is partial suboperad of $K$,
and conformal full field algebra over $V^L\otimes V^R$ are
exactly smooth algebras over 
$\tilde{K}^{c^L}\otimes \overline{\tilde{K}^{\overline{c^R}}}$.
The modular invariance property of conformal full field algebras
will be studied in \cite{HK3}.

In order to cover the entire genus-zero conformal field theories, 
one still needs to consider
the Riemann surfaces with more than one copy of $S^1$ 
in the codomain. In chiral theory, it was studied by Hubbard
in the framework of so-called vertex operator coalgebras 
\cite{Hub1}\cite{Hub2}.  
In a theory including both chiral and antichiral parts, 
it amounts to have a conformal full field algebra 
equipped with a nondegenerate invariant bilinear form \cite{HK2}. 
Conversely,  a conformal full field algebra with 
a nondegenerate invariant bilinear form 
gives a complete set of data needed for 
a genus-zero conformal field theory (nonunitary). In particular, we 
show in Section 2 that such conformal full field algebras
are just algebraic representations of the sewing operations 
among spheres with arbitary number of negatively oriented and
positively oriented punctures as long as the resulting surfaces
after sewing are still of genus-zero. 
The invariant property of bilinear form used 
in this work is slightly different from that in \cite{HK2}. 
Both definitions have
clear geometric meanings. But we need the new definition
for a reason which is explained later.

Although the notion of 
conformal full field algebras has the advantage of 
being a pure algebraic formulation 
of a part of genus-zero conformal field theories, 
its axioms are still very hard to check directly. 
The theory of the tensor products of the modules of
vertex operator algebras is developed by Huang and Lepowsky
\cite{HL1}-\cite{HL4}\cite{H2}\cite{H7}. 
The following Theorem proved by Huang in \cite{H7} is very 
crucial for our constructions of conformal full field algebras. 
\begin{thm} \label{ioa}
Let $V$ be a vertex operator algebra satisfying 
the following conditions: 
\bnu
\item Every $\C$-graded generalized $V$-module is a direct sum of 
$\C$-graded irreducible $V$-modules,
\item There are only finitely many inequivalent $\C$-graded
irreducible $V$-modules, 
\item Every $\R$-graded irreducible $V$-module satisfies the 
$C_1$-cofiniteness condition. 
\enu
Then the direct sum of all (in-equivalent) 
irreducible $V$-modules has
a natural structure of intertwining operator algebra and 
the category of $V$-modules, denoted as 
$\mathcal{C}_V$ has a natural structure of 
vertex tensor category. 
In particular, $\mathcal{C}_V$ 
has natural structure of braided tensor category. 
\end{thm}

\begin{ass} \label{assumption}
{\rm 
All the vertex operator algebras appeared in this work, $V$, $V^L$
and $V^R$ are all assumed to satisfy the conditions in 
Theorem \ref{ioa} without further announcement.  
Sometimes we even assume stronger
conditions on them, as we will do explicitly in Section 4 and 
Section 5. 
}
\end{ass}

In \cite{HK2}, we have studied in detail the properties of 
conformal full field algebras over $V^L\otimes V^R$. 
In particular, an equivalent definition of conformal
full field algebra over $V^L\otimes V^R$ is also given.
We recalled this result in Theorem \ref{ffa-thm}.   
The axiom in this definition are much easier to verify than
those in the original definition of conformal full field algebra.  
Also by Theorem \ref{ioa},  
the categories of $V^L$-modules, $V^R$-modules and 
$V^L\otimes V^R$-modules, denoted as 
$\mathcal{C}_{V^L}$, $\mathcal{C}_{V^L}$ and $\mathcal{C}_{V^L\otimes V^R}$
respectively, all have the structures 
of braided tensor category \cite{H7}\cite{HK2}.
In this work, the braiding structure in 
$\mathcal{C}_{V^L\otimes V^R}$ is chosen to
be different from the one obtained from Theorem \ref{ioa}.   
Using Theorem \ref{ffa-thm} and 
this new braiding structure on $\mathcal{C}_{V^L\otimes V^R}$, 
we can show, without much effort, that
a conformal full field algebra over $V^L\otimes V^R$
is equivalent to a commutative associative 
algebra in $\mathcal{C}_{V^L\otimes V^R}$ with a trivial twist.

We would also like to give a categorical formulation of 
conformal full field algebra over $V^L\otimes V^R$ equipped 
with a
nondegenerate invariant bilinear form, where $V^L$ and $V^R$
are assumed to satisfy 
the conditions in Theorem \ref{MTC}, which are
slightly stronger than those in Theorem \ref{ioa}. 
It turns out that the way we define the invariant property of 
bilinear form of conformal full field algebra in 
\cite{HK2} is not easy to work with categorically. 
This is the reason why a slightly modified notion of invariant
bilinear form is introduced in Section 2. 
This modification leads us to consider 
on the graded dual space of a module 
over a vertex operator algebra a module structure 
which is different from (but equivalent to) 
the usual contragredient module structure \cite{FHL}.
Because of this, we redefine the duality maps \cite{H14}\cite{H10} 
in this new convention and prove the rigidity (in appendix). 
Then we show in detail, in Section 4,  
that a conformal full field algebra 
over $V^L\otimes V^R$ equipped 
with a nondegenerate invariant bilinear form
in the new sense exactly amounts to a commutative Frobenius
algebra with a trivial twist in $\mathcal{C}_{V^L\otimes V^R}$.

Once the categorical formulation is known. We can give a
categorical construction of conformal full field algebras
equipped with a nondegenerate invariant bilinear form. 
This construction was previously given by Huang and the author 
in \cite{HK2} by using intertwining operator algebras,
and was also known to 
physicists as diagonal construction (see for example 
\cite{FFFS} and references therein).

Recently, Fuchs, Runkel, Schweigert and Fjelstad
have proposed a very general 
construction of all correlation functions of 
boundary conformal conformal field theories using 
3-dimensional topological field theories in a series of 
papers \cite{FRS1}-\cite{FRS4}\cite{FjFRS}\cite{RFFS}\cite{SFR}. 
In particular, a construction of commutative associative 
algebras in $\mathcal{C}_{V\otimes V}$ is explicitly given
in \cite{RFFS}. Our approach is somewhat complementary 
to their approach (see \cite{RFFS} for comments on
the relation of two approachs). We hope that two approachs
can be combined to obtain a rather complete picture of
conformal field theory in the near future.

The layout of this paper is as follow. In Section 1, 
we study the operadic formulation of  
grading-restricted $\R\times \R$-graded full field algebra and
its variants, following the work of Huang \cite{H4}.  
In Section 2, we give a the geometric description of
a conformal full field algebra over $V^L\otimes V^R$ 
equipped with a nondegenerate invariant bilinear form.
In Section 3, we give a categorical formulation of 
conformal full field algebra over $V^L\otimes V^R$. 
In Section 4, we give a categorical formulation 
of conformal full field algebra over $V^L\otimes V^R$ with
a nondegenerate invariant bilinear form for $V^L$ and $V^R$ 
satisfying the conditions in Theorem \ref{MTC}. 
In Section 5, we give a
categorical construction of conformal full field algebras 
over $V^L\otimes V^R$ and prove that such obtained 
conformal full field algebras over $V^L\otimes V^R$
are naturally equipped with a nondegenerate 
invariant bilinear form. 

For the convenience of readers, the materials in 
Section 3, 4, 5 are completely independent of 
those in Section 1, 2. For those who is only 
interested in categorical formulation of conformal full field 
algebras over $V^L\otimes V^R$, 
it is harmless to start from Section 3 directly.

Convention of notations: 
$\N, \R, \R_+, \C, \HH, \hat{\C}, \hat{\HH}$ 
denote the set of
natural numbers, real numbers, positive real numbers and
complex numbers, and upper half plane, 
one point compactification of $\C$ and $\HH\cup \R$, respectively. 
We also use $I_F$ to denote the identity map on a vector space $F$.

Note. After this paper appeared online in math arxiv,
the author noticed that the categorical construction of
commutative associative algebra in $\mathcal{C}_{V\otimes V}$ 
given in this work is nothing but a basis independent 
version of that in \cite{FrFRS} which appeared earlier. 
There is another natural point of view of this
construction in terms of adjoint functors. It will be discussed elsewhere.

\paragraph{Acknowledgment}
This work grows from a chapter in author's thesis. 
I want to thank my advisor Yi-Zhi Huang 
for his constant support and many important suggestions. 
I also want to thank him for spending incredible amount of 
time in helping me to improve the writting
of my thesis and this work. 
I thank J. Lepowsky and C. Schweigert for 
many inspiring conversations related to this work.

\setlength{\unitlength}{1cm}

\renewcommand{\theequation}{\thesection.\arabic{equation}}
\renewcommand{\thethm}{\thesection.\arabic{thm}}
\setcounter{equation}{0}
\setcounter{thm}{0}

\section{Operadic formulations of full field alegbras}

In this section, we study the operadic formulation of 
full field algebra. In section 1.1, 
we recall the notion of sphere partial operad $K$ and its
partial suboperad $\widehat{K}$. We introduce the notion of 
smooth function on $K$. In section 1.2, 
we recall the notions of determinant line bundle over $K$ 
and the $\C$-extensions of $K$, such as 
$\tilde{K}^c$ and 
$\tilde{K}^{c^L} \otimes \overline{\tilde{K}^{\overline{c^R}} }$
for $c, c^L, c^R\in \C$. 
Section 1.1 and 1.2 are mainly taken from 
\cite{H4}. The readers who is interested in knowing more on 
this subject should 
consult with \cite{H4} for details. 
In section 1.3, we recall the notion of algebra over
partial operad \cite{H4}, and explain what it means for an 
algebra over $\widehat{K}$ and 
$\tilde{K}^{c^L} \otimes \overline{\tilde{K}^{\overline{c^R}} }$
to be smooth. Then we give two isomorphism theorems. 
The first one says that the category of grading-restricted 
$\R\times \R$-graded full field algebras is isomorphic to
the category of smooth $\widehat{K}$-algebras. 
The second one says that
the category of conformal full field algebras over $V^L\otimes V^R$
is isomorphic to the category of smooth
$\tilde{K}^{c^L} \otimes \overline{\tilde{K}^{\overline{c^R}}}$-algebras
over $V^L\otimes V^R$. We give a selfcontent proof of the 
first isomorphism theorem. The proof of the second isomorphism
theorem is technical, and heavily depends on 
the results in \cite{H4}.

\subsection{Sphere partial operad $K$}

A sphere with tubes of type $(n_-,n_+)$ is a sphere $S$ with 
$n_-$ ordered punctures $p_i, i=1,\dots, n_-$ 
and $n_+$ ordered punctures $q_j, j=1, \dots, n_+$, together 
with a negatively oriented local chart $(U_i, \varphi_i)$
around each $p_i$ and 
a positively oriented local chart $(V_j, \psi_j)$ around
each $q_i$, where $U_i$ and $V_j$ are neighborhood of 
$p_i$ and $q_j$ respectively 
and the local coordinate maps
$\varphi_i: U_i\rightarrow \C$ and $\psi_j: V_j\rightarrow \C$
are conformal maps so that $\varphi(p_i)=\psi_j(q_j)=0$.
 
The conformal equivalence of sphere with tubes is defined to be
the conformal maps between two spheres so that the germs of 
local coordinate map $\varphi_i$ and $\psi_j$ are preserved. 
We then obtain
a moduli space of sphere with tubes of type $(n_-, n_+)$. In this
section, We are
only interested in spheres with tubes of type $(1,n)$. 
For this type of spheres with tubes, 
we label the only negative oriented puncture as the $0$-th puncture.
We denote the moduli space of sphere with tubes of type $(1,n)$ as
$K(n)$. 

Using automorphisms of sphere, 
we can select a canonical representative
from each conformal equivalence class in $K(n)$ 
for all $n>0$ 
by fixing the $n$-th puncture at $0\in \C$,  the 
$0$-th puncture at $\infty$, and $\psi_0$, the local 
coordinate map at $\infty$, to be so that
\beq \label{can-rep-cond}
\lim_{w\rightarrow \infty} w\psi_0(w) = -1.
\eeq
As a consequence, the moduli space $K(n), n\in \Z_+$ 
can be identified with
$$
K(n) = M^{n-1} \times H \times (\C^{\times} \times H)^n 
$$
where 
$$
M^{n-1} = \{ (z_1, \dots, z_{n-1}) | z_i\in \C^{\times}, z_i\neq z_j,
\mbox{ for $i\neq j$} \}
$$
and 
\bea
H &=& \{ A=(A_1, A_2, \dots) \in \prod_{i=1}^{\infty} \C \quad
| A_i\in \C, \quad e^{\sum_{j=1}^{\infty} 
A_j x^{j+1}\frac{d}{dx} }x   \nn
&&\hspace{0.5cm} \mbox{ is an absolute convergent series 
in some neighborhood of $0$ } \}.  \nonumber
\eea
Similarly, using automorphisms of sphere, we can choose 
the canonical representatives of conformal 
equivalence classes of sphere
with tube of type (1,0) so that 
the moduli space $K(0)$ can be identified with the set 
$$K(0) = \{ A \in  \prod_{i=1}^{\infty}  \C  |  A_1=0\}.$$
We will denote a general element in $K(n), n\in \N$ as 
\beq  \label{ele-K-n}
P=(z_1,\ldots, z_{n-1}; A^{(0)}, (a_0^{(1)}, A^{(1)}), \ldots, 
(a_0^{(n)}, A^{(n)})),
\eeq
where $z_1,\dots, z_{n-1}$ 
denotes the location of the first $n-1$ punctures
and the rest of the data give the local coordinate maps at 
the negatively oriented ($0$-th) puncture and 
positively oriented punctures respectively as follow:
\bea  
f_{0}(w) &=& -e^{\sum_{j\in \Z_+} A_j^{(0)}x^{j+1} \frac{d}{dx} } x 
\lbar_{x=\frac{1}{w}} \label{f-0} \\
f_{i}(w) &=& a_0^{(i)} e^{\sum_{j\in \Z_+} A_j^{(i)}w^{j+1} \frac{d}{dw} } w, \quad
\forall i=1,\dots, n.  \label{f-i}
\eea

Let $P\in K(m)$ and $Q\in K(n)$. Let $\bar{B}^r$ 
be the closed ball in $\C$ centered at $0$ 
with radius $r$,   $\varphi_i$ 
the germs of local coordinate map at $i$-th puncture $p_i$
of $P$, and $\psi_0$ the germs of local coordinate 
map at $0$-th puncture $q_0$ of $Q$.  
Then we say that the $i$-th tube of $P$ can
be sewn with $0$-th tube of $Q$ if there is a $r\in \R_+$ such 
that $p_i$ and $q_0$ are the only punctures in 
$\varphi_i^{-1}(\bar{B}^r)$
and $\psi_0^{-1}(\bar{B}^{1/r})$ respectively. 
A new sphere with tubes in $K(m+n-1)$, denoted as 
$P _{^i} \infty_{^0} Q$,  can be obtained by
cutting out $\varphi_i^{-1}(\bar{B}^r)$
and $\psi_0^{-1}(\bar{B}^{1/r})$ from $P$ and $Q$ respectively, and
then identifying the boundary circle via the map 
$$
\psi^{-1} \circ J_{\hat{\HH}} \circ \varphi_i
$$
where $J_{\hat{\HH}}: w\rightarrow \frac{-1}{w}$.

\begin{rema}
{\rm
$\HH$ denotes the upper hand plane. We use it to remind 
us the fact that $w\rightarrow \frac{-1}{w}$ is an 
automorphism of upper hand plane. 
}
\end{rema}

Therefore, we have sewing operations:
$$
_{^i} \infty_{^0} : K(m) \times K(n) \rightarrow K(m+n-1)
$$
partially defined on the entire $K$ between two spheres
with tubes along two oppositely oriented tubes.

\begin{rema}
{\rm
Our definition of sewing operation is 
defined differently from that defined in \cite{H4}, where 
$J_{\hat{\C}}: w\mapsto \frac{1}{w}$ 
is used in the place of $J_{\hat{\HH}}$. 
Also notice that our 
convention of local coordinate map (\ref{f-0}) is 
also different from that in \cite{H4} by a sign. 
Combining effect of these two difference is a trivial one. 
Namely, our definition of sewing operation is 
equivalent to that in \cite{H4} if one identify 
each sphere with tube with local coordinate map at $\infty$
being $f_0(w)$ used in this paper, with the same sphere 
with tube but with local coordinate map at $\infty$ be $-f_0(w)$
used in \cite{H4}. This identification actually gives 
an isomorphism of partial operad.
}
\end{rema}

\begin{rema}
{\rm
One reason for introducing this
new convention is to make sure that 
the invariant bilinear form on conformal full field algebra 
introduced in Section 2 has a clear geometric meaning. 
The notion of invariant bilinear form on a conformal
full field algebra used in this work is slightly different from 
that in \cite{HK2} for a categorical reason 
(see also the Remark \ref{rema-inv-biform}).
} 
\end{rema}

\begin{rema} \label{J-H-C}
{\rm 
Another reason for our choice being geometrically natural is that 
$J_{\hat{\HH}}$ is an automorphism of $\hat{\HH}$ 
while $J_{\hat{\C}}$ is not. 
This will be important in boundary conformal field theories, 
where we study algebras over a partial operad consisting 
of disks with strips \cite{HK1}, which is often modeled as 
upper half planes. Hence it is more natural to use 
$J_{\hat{\HH}}$ instead of $J_{\hat{\C}}$ there.  
Moreover, in the study of boundary conformal field theories, 
it is necessary to embed disks with strips 
into spheres with tubes by certain 
doubling maps (\cite{HK1}\cite{Ko}). 
In order for the embedding becoming a morphism of partial operad, 
it is also natural to define the sewing operations in $K$ using 
$J_{\hat{\HH}}$ instead of $J_{\hat{\C}}$. 
}
\end{rema}

Let $\mathbf{0}$ be the sequence $(0,0,\dots )$. The element 
$$
I_K:=(\mathbf{0}, (1, \mathbf{0})) \in K(1)
$$ 
is called identity element. Its tubes are sewable with
any oppositely oriented tubes in other sphere. 
Moreover, it satisfies the following identity property:
$$
Q _{^i}\infty_{^0} I = I _{^1}\infty_{^0} Q = Q
$$ 
for any $Q\in K(n), n\in \N$. 
The subset 
\beq  \label{res-group}
\{ (\mathbf{0}, (a, \mathbf{0}) \in H\times (\C^{\times} \times H)
\}
\eeq
of $K(1)$ together with the sewing operation $_{^1}\infty_{^0}$
is a group isomorphic to $\C^{\times}$. 
There is also an obvious action of permutation group $S_n$ on $K(n)$. 
The following result is proved in \cite{H4}. 

\begin{prop}
The collection of sets
$$
K = \{ K(n) \}_{n\in \N}
$$
together with $I_K$, sewing operations, the actions of $S_n$ on $K(n)$ 
and the group (\ref{res-group}), 
is a $\C^{\times}$-rescalable partial operad. 
\end{prop}

Let $A(a;i) =\{ A_j | A_i=a, A_j=0, j\in \Z_+, j\neq i \}$.
For simplicity, we will also use $A(a;i)$ to 
denote the element in $K(0)$ such that the local coordinate
map at $\infty$ is given by 
$$
-\exp \left( a \left( \frac{1}{w} \right)^{i+1} 
\frac{d}{d\frac{1}{w}}
\right) \frac{1}{w} = -\exp \left( -a w^{-i+1} \frac{d}{dw} \right)
\frac{1}{w}. 
$$

Let $\widehat{K}(n)$ be the subset of $K(n)$ consisting of
elements of the form
\beq  \label{K-hat-ele}
(z_1, \dots, z_{n-1}; A(a;1), (a_0^{(1)}, \mathbf{0}), 
\dots, (a_0^{(n)}, \mathbf{0}) ).
\eeq
Then $\widehat{K} = \{ \widehat{K}(n) \}_{n\in \N}$
is a partial suboperad of $K$ \cite{H4}.

We use ``overline'' to denote complex conjugation. 
A function $f$ on $K(n), n\in \N$ is called {\it smooth}  
if there is a $N\in \N$ such that $f$ can be written as 
\bea \label{HA-smooth}
&&\sum_{k=1}^N 
g_k(A^{(0)}, 
a_0^{(1)}, A^{(1)},  \dots, a_0^{(n)}, A^{(n)}; 
\overline{A^{(0)}}, \overline{a_0^{(1)}}, \overline{A^{(1)}},  
\dots, \overline{a_0^{(n)}},\overline{A^{(n)}} )  \nn
&&\hspace{1cm} 
h_k(z_1, \dots, z_{n-1}; \bar{z}_1, \dots, \bar{z}_{n-1})
\eea
where $g_k$ are polynomial functions of 
$A_j^{(i)}$ and $\overline{A_j^{(i)}}$, for $i=0, \dots, n, j\in \Z_+$ 
and linear combination of
$\prod_{i=1}^n (a_0^{(i)})^{r_i}\overline{a_0^{(i)}}^{s_i}$ for 
some $r_i, s_i \in \R$,  
and $h_k$ are smooth functions of $z_i, i=1, \dots, n-1$. 
It is clear that there is also a naturally induced notion of 
smooth function on $\widehat{K}$.

A tangent space $TK(n)$ of $K(n), n\in \N$ 
can be defined naturally. Let $\epsilon \in \C$.
We define two tangent vectors in $T_IK(1)$ as follow:
for any smooth function $f$ on $K(1)$, 
\bea
\mathcal{L}_I(z)f &:=& \frac{\partial}{\partial \epsilon} 
\lbar_{\epsilon=0} 
f(P(z) _{^1}\infty_{^0} A(\epsilon; 2)),   \nn
\bar{\mathcal{L}}_I(\bar{z}) f &:=& 
\frac{\partial}{\partial \bar{\epsilon}} 
\lbar_{\epsilon=0}  f(P(z) _{^1}\infty_{^0} A(\epsilon; 2)). 
\eea
The following proposition is a generalization 
of Proposition 3.2.5 in \cite{H6}. 
\begin{prop} 
\bea
\mathcal{L}_I(z) &=& 
z^{-2} \frac{\partial}{\partial a^{(1)}} \lbar_I 
+ \sum_{i=0}^1 \sum_{j=1}^{\infty} z^{-(2i-1)j-2} 
\frac{\partial}{\partial A_{j}^{(i)}}
\lbar_I , \nn
\bar{\mathcal{L}}_I(\bar{z})
&=& \bar{z}^{-2}\frac{\partial}{\partial \overline{a^{(1)}} } \lbar_I 
+ \sum_{i=0}^1 \sum_{j=1}^{\infty} \bar{z}^{-(2i-1)j-2} 
\frac{\partial}{\partial \overline{A_{j}^{(i)}} }
\lbar_I. \label{L-I-equ}
\eea
\end{prop}
\pf
Let $P=P(z)_{^1}\infty_{^0} A(\epsilon,2)$ with its coordinates in 
moduli space given by 
$$
(z; A^{(0)}, (a_0^{(1)}, A^{(1)}), (a_0^{(2)}, A^{(2)}) ).
$$
It is shown in \cite{H4} that 
$z$, $A_j^{(0)}$, $a_0^{(1)}$, $A_j^{(1)}$, $a_0^{(2)}$ and $A_j^{(2)}$, 
$j\in \Z_+$ are holomorphic functions of $\epsilon$. 
Hence their complex conjugation 
$\bar{z}$, $\overline{A_j^{(0)}}$, $\overline{a_0^{(1)}}$, 
$\overline{A_j^{(1)}}$, $\overline{a_0^{(2)}}$ 
and $\overline{A_j^{(2)}}$ are holomorphic with respect to 
$\bar{\epsilon}$. Therefore (\ref{L-I-equ}) 
follows from Proposition 3.2.5 in \cite{H4}. 
\epf

\subsection{Determinant line bundle over $K$}

The determinant line bundle over $K$ and 
the $\C$-extensions of $K$ are studied in \cite{H4}. 
For each $n\in \N$, the determinant line bundle 
$\text{Det}(n)$ over $K(n)$ is a trivial bundle over $K(n)$.
We denote the fiber at $Q\in K(n)$ as $\text{Det}_Q$. 
There is a canonical section of $\text{Det}(n)$, 
denoted by $\psi_n$, for each $n\in \N$. 
For any element $Q\in K(n)$, 
let $\mu_n(Q)$ be the element of the fiber over $Q$ given by
$$
\psi_n(Q)=(Q, \mu_n(Q)).
$$
Then there is a $\lambda_{Q}$ for each element $\tilde{Q}$ of 
$\text{Det}_Q$
such that $\tilde{Q}=(Q, \lambda_Q)$ and $\lambda_{Q} = \alpha \mu_n(Q)$ 
for some $\alpha \in \C$. 

Consider the following two general elements in $K(m)$ and $K(n)$.
\bea  
P&=&(z_1,\ldots, z_m; 
A^{(0)}, (a_0^{(1)}, A^{(1)}), \ldots, (a_0^{(m)}, A^{(m)}) ) \nn
Q&=&(\xi_1,\ldots, \xi_n; 
B^{(0)}, (b_0^{(1)}, B^{(1)}), \ldots, (b_0^{(n)}, B^{(n)}) ) 
\eea
If $P _{^i}\infty_{^0} Q$ exists,  
then there is a canonical isomorphism:
$$
l_{P, Q}^i: \text{Det}_{P} \otimes \text{Det}_{Q} 
\rightarrow \text{Det}_{P _{^i}\infty_{^0} Q},
$$
given as:
\bea  \label{sew-det-fiber}
&& l_{P, Q}^i (a_1\mu_m(P) \otimes a_2 \mu_n(Q) ) \nn
&& \hspace{2cm}= a_1a_2 e^{2\Gamma(A^{(i)}, B^{(0)}, a_0^{(i)})} 
\mu_{m+n-1}(P _{^i}\infty_{^0} Q),
\eea
where $\Gamma$ is a $\C$-valued analytic function 
of complex variables $A^{(i)}_j, B^{(0)}_k, a_0^{(i)}, j,k\in \N$. 
If we expand $\Gamma$ as formal series, we have
$$
\Gamma( A^{(i)}, B^{(0)}, \alpha) 
\in \Q [\alpha, \alpha^{-1}] [[ A^{(i)}, B^{(0)} ]].
$$
For a detailed discussion of $\Gamma$, see chapter $4$ in \cite{H4}. 
In particular, we have
\begin{equation} \label{gamma-real-analytic}
\overline{\Gamma(A^{(i)}, B^{(0)}, a_0^{(i)})} = 
\Gamma(\overline{A^{(i)}}, \overline{B^{(0)}}, \overline{a_0^{(i)}}).
\end{equation}
For $(P,\lambda_P)\in \text{Det}(m)$ and 
$(Q,\lambda_Q)\in \text{Det}(n)$ 
such that $P _{^i}\infty_{^0} Q$ exists, then we define a partially
defined map
$$
_{^i}\widetilde{\infty}_{^0}^2: 
\text{Det}(m) \times \text{Det}(n) \rightarrow \text{Det}(m+n-1)
$$
by 
$$
(P,\lambda_P) _{^i}\widetilde{\infty}_{^0}^2 (Q,\lambda_Q) = 
(P _{^i}\infty_{^0} Q, l_{P, Q}^i(\lambda_P \otimes \lambda_Q) ).
$$
Using this partial operation, one obtain a 
$\C^{\times}$-rescalable partial operad structure on 
$\text{Det}=\{ \text{Det}(n) \}_{n\in \N}$.

For $c\in \C$, the so-called 
{\em vertex partial operad of central charge $c$}, 
denoted as $\tilde{K}^c$, is the $\frac{c}{2}$-th power of 
the determinant line bundle $\text{Det}$ over $K$. 
$\tilde{K}^c$ also has a structure of partial operad. 
The construction of sewing operations on $\tilde{K}^c$
is same as that on $\text{Det}$ except 
replacing (\ref{sew-det-fiber}) by 
\bea \label{sew-k-c-fiber}
&&(l_{P, Q}^i)^c (a_1\mu_m(P) \otimes a_2 \mu_n(Q) ) \nn
&& \hspace{2cm}= a_1a_2 e^{\Gamma(A^{(i)}, B^{(0)}, a_0^{(i)})c} 
\mu_{m+n-1}(P _{^i}\infty_{^0} Q).
\eea
We denote the corresponding sewing operation 
as $_{^i}\widetilde{\infty}_{^0}^c$. 
It is proved in \cite{H4} that 
$\tilde{K}^c$ is also $\C^{\times}$-rescalable 
associative partial operad and a $\C$-extension of $K$ (see
\cite{H4}).

The complex conjugation
$\overline{\tilde{K}^{\bar{c}}}$ of holomorphic line bundle 
$\tilde{K}^{\bar{c}}$ also has a natural structure of partial operad.
The section of $\psi$ on
$\tilde{K}^{\bar{c}}$ canonically gives a section $\bar{\psi}$ 
on $\overline{\tilde{K}^{\bar{c}}}$.
Using the global section $\bar{\psi}$, the canonical isomorphism
$$
\bar{l}_{P,Q}^{\, i} : \overline{\text{Det}^{\bar{c}} }_{P} \otimes 
\overline{\text{Det}^{\bar{c}} }_{Q} \rightarrow 
\overline{\text{Det}^{\bar{c}} }_{P _{^i}\infty_{^0} Q}, 
$$
for any pair of $P,Q\in K$ so that $P _{^i}\infty_{^0} Q$ exists,
can be written as
\begin{equation} \label{conj-sew-fiber}
\bar{l}_{P,Q}^{\, i} ( \lambda_1 \otimes \lambda_2) = 
\lambda_1\lambda_2 e^{\Gamma(\overline{A^{(i)}}, 
\overline{B^{(0)}}, \overline{a_0^{(i)}} )c } 
\end{equation} 
where $(a_0^{(i)},A^{(i)})$ gives the local coordinate map at 
$i$-th positively oriented puncture
in $P$ and $B^{(0)}$ gives the local coordinate map at 
$\infty$ in $Q$. Observe that (\ref{conj-sew-fiber}) implies that 
$\overline{\tilde{K}^{\bar{c}}}=\overline{\tilde{K}}^{\,\, c}$
as partial operads. 

In this work, we are also interested in the 
tensor product bundle $\tilde{K}^{c^L} 
\otimes \overline{\tilde{K}^{\overline{c^R}}}$ for $c^L, c^R\in \C$. 
It is clear that it is a $\C^{\times}$-rescalable partial operad
as well. The natural section induced from $\tilde{K}^{c^L}$ and 
$\overline{\tilde{K}^{\overline{c^R}}}$
is simply $\psi \otimes \bar{\psi}$. We will denote 
the sewing operation on $\tilde{K}^{c^L} 
\otimes \overline{\tilde{K}^{\overline{c^R}}}$ simply as 
$\widetilde{\infty}$ without making its dependents on
$c^L, c^R$ explicit. 

A function on $\tilde{K}^c$ or $\tilde{K}^{c^L} 
\otimes \overline{\tilde{K}^{\overline{c^R}}}$ is called {\it smooth} 
if it is smooth on the base space $K$ and linear on fiber.

\subsection{Two isomorphism theorems}
In this subsection we discuss the operadic formulations of 
grading-restricted $\R\times \R$-graded full
field algebras and conformal full field algebras over $V$.
We apologize for not recalling the definitions of 
these two notions here. They can be found in \cite{HK2}.

We first recall the definition of algebra over partial operad. 
Let $G$ be a group and $U$ a complete reducible $G$-module and 
$W$ a $G$-submodule of $U$.
We will use $\overline{U}$ to denote the algebraic completion of 
$U$. It should be distinguishable with complex conjugation from
the context.  
Let $\mathcal{H}_{UW}^G(n)$ be the set of 
multilinear maps $U^{\otimes n} \rightarrow \overline{U}$ such that
image of $W^{\otimes n}$ is in $\overline{W}$. 
$\mathcal{H}_{UW}^G:=\{ \mathcal{H}_{UW}^G(n)\}_{n\in \N}$ 
has a natural structure of partial pseudo-operad
\cite{H4}, which satisfies all the axioms of partial operad
except the associativity. 

\begin{defn}  \label{def-P-alg}
{\rm 
Let $\mathcal{P}$ be a partial operad with rescaling group $G$. A 
$\mathcal{P}$-pseudo-algebra is a triple $(U,W,\nu)$, 
where $U$ is a completely reducible $G$-module 
$$
U = \coprod_{M\in A} U_{(M)}
$$
in which $A$ is the set of equivalent 
class of irreducible $G$-modules,
$W$ is a submodule of $U$, and $\nu$ is a morphism from
$\mathcal{P}$ to the partial pseudo-operad 
$\mathcal{H}_{U,W}^G$ (\cite{H3}),
satisfying the following conditions
\bnu
\item $\dim U_{(M)}<\infty$ for all $M\in A$. 

\item The submodule of $U$ generated by the 
homogeneous components of 
the elements of $\nu_0(\mathcal{P}(0))$ is $W$ \footnote{Here, 
$\nu_n$ for $n\in \N$ is simply the restriction of $\nu$ on 
$\mathcal{P}(n)$.}.
\item The map from $G$ to $\mathcal{H}_{U,W}^G(1)$ induced from 
$\nu_1$ is the given representation of $G$ on $U$. 
\enu
If $\mathcal{P}$ is rescalable \cite{H4}, 
we call a $\mathcal{P}$-pseudo-algebra a 
$\mathcal{P}$-algebra.  
}
\end{defn}

In this work, we are interested in studying $\widehat{K}$-algebras  and 
$\tilde{K}^{c^L} \otimes \overline{\tilde{K}^{\overline{c^R}}}$-algebras,
both of which are $\C^{\times}$-rescalable partial operads. 
\begin{defn}
{\rm A $\widehat{K}$-algebra (or 
$\tilde{K}^{c^L} \otimes \overline{\tilde{K}^{\overline{c^R}}}$-algebra) 
$(F, W, \nu)$ is called {\em smooth} if it satisfies 
the following conditions:
\bnu
\item $F_{(m,n)}=0$ if the real part of $m$ or $n$ is sufficiently small. 

\item For any $n\in \N$, $w'\in F', w_1, \dots, w_{n+1}\in F$,
an element $Q$ in $K(n)$ (or  
$\tilde{K}^{c^L} \otimes \overline{\tilde{K}^{\overline{c^R}}}(n)$),
the function
$$
Q \mapsto 
\langle w', \nu(Q)(w_1\otimes \dots \otimes w_{n+1})\rangle 
$$
is smooth on $K(n)$ (or on 
$\tilde{K}^{c^L} \otimes \overline{\tilde{K}^{\overline{c^R}}}(n)$). 
\enu 
}
\end{defn}

We choose the branch cut of logarithm as follow
\beq
\log z = \log |z| + \text{Arg} z, \quad  0\leq \text{Arg} z < 2\pi. 
\eeq
We define the power functions, 
$z^m$ and $\bar{z}^n$ for $m,n \in \R$, 
to be $e^{m\log z}$ and $e^{n\overline{\log z}}$ respectively. 
This convention is used throughout this work.

The following lemma must be well-known. But we are not aware of 
any source of reference. So we give a proof here. 
\begin{lemma}  \label{double-grad}
If $\rho: \C^{\times} \rightarrow \C^{\times}$ 
is an irreducible representation of group $\C^{\times}$, then 
$\rho(z)=z^{-m} \bar{z}^{-n}$  for some  $m,n\in \C$ and $m-n\in \Z$. 
\end{lemma}
\pf
We view $\C^{\times}$ as a topological group. 
All the irreducible representations of $\C^{\times}$ are one-dimensional. 
They are continuous maps $\C^{\times} \rightarrow \C^{\times}$ 
as group homomorphisms.

$\C$ is the universal covering space of $\C^{\times}$ with the projection map 
$\pi$ given by $z\mapsto e^{z}$. By the lifting property of covering space, there is
a unique continuous map $\tilde{\rho}: \C \rightarrow \C$ such that 
$\rho \pi = \pi \tilde{\rho}$ if we choose $\tilde{\rho}: 0\mapsto 0$. 
Let $s_1, s_2\in \C$. 
Since $\rho$ is a group homomorphism, hence we have 
\beq
e^{\tilde{\rho }(s_1+s_2)} = \rho( e^{s_1+s_2} )  
= \rho(e^{s_1}) \rho(e^{s_2}) = e^{\tilde{\rho } (s_1)} e^{\tilde{\rho } (s_2)} 
= e^{\tilde{\rho } (s_1) + \tilde{\rho } (s_2)}.
\eeq
It implies that 
$\tilde{\rho }(s_1+s_2)= \tilde{\rho } (s_1) + \tilde{\rho } (s_2) + k2\pi i$
for some $k\in \Z$. Notice that $k$ must be unique because $\tilde{\rho}$ is continuous. 
Recall that we have already chosen $\tilde{\rho}$ to be so that $\tilde{\rho}(0)=0$. 
Therefore $k=0$. $\tilde{\rho}$ 
is actually a linear map from $\C \rightarrow \C$. 
Namely, $\tilde{\rho}$ can be written as 
\beq \label{rho}
\tilde{\rho}:  \left( \begin{array}{c} x \\ y \end{array} \right) \mapsto 
\left( \begin{array}{cc} a & b \\ c & d \end{array} \right)
\left( \begin{array}{c} x \\ y \end{array} \right) 
\eeq
where $x, y\in \R$ are the real part and the imaginary 
part of a complex number in $\C$. 

Moreover, the group homomorphism preserves the identity, 
i.e. $\rho: 1\mapsto 1$. It 
implies that $\tilde{\rho}(2\pi i) = l 2\pi i$ for some $l\in \Z$. 
Applying this result to (\ref{rho}), we obtain that
$b=0$ and $d=l\in \Z$. Conversely, it is easy to see
that every $\tilde{\rho}$ of form (\ref{rho}) with $b=0$ and $d\in \Z$
gives arise to a group homomorphism $\rho: \C^{\times} \rightarrow \C^{\times}$ of 
the following form: 
\beq
z=e^{x+iy} \mapsto e^{ax + i(cx+dy)} = e^{(a+ic)x} e^{d i y} 
= |z|^{a+ic} \left( \frac{z}{\bar{z}}\right)^{d/2}= 
z^{-m} \bar{z}^{-n}.
\eeq
where $m=-\frac{a+ic}{2}-\frac{d}{2}$ and 
$n=-\frac{a+ic}{2}+\frac{d}{2}$
and $m-n=-d\in \Z$.  
\epf

Now we study the basic properties of a smooth 
$\widehat{K}$-algebras. We fix a smooth $\widehat{K}$-algebras
$(F, W, \nu)$. 

The set $I:=\{ (m,n)\in \C\times \C | m-n\in \Z\}$ together 
with the usual addition operation gives an abelian group. 
By Lemma \ref{double-grad} 
any $\widehat{K}$-algebras must be $I$-graded. Namely,  
$$
F=\coprod_{(m,n)\in I} F_{(m,n)}.
$$
Since $\nu_1( (\mathbf{0},(a,\mathbf{0})))$ for $a\in \C^{\times}$ 
gives the representation 
of $\C^{\times}$ on $F$ by the definition of algebra over operad, 
we must have 
$\nu_1((0,a))u= a^{-m} \bar{a}^{-n} u$ for $u\in F_{(m,n)}$. 
Let $\mathbf{d}^L$ and $\mathbf{d}^R$ be the grading operators such that 
$\mathbf{d}^L u = m u$ and $\mathbf{d}^R u =n u$ for $u\in F_{(m,n)}$. 
We call $m$ the left weight of $u\in F_{(m,n)}$
and $n$ the right weight of $u$, and 
denote them as $\wt^L u$ and $\wt^R u$ respectively.
We also define $\wt u:=\wt^L u + \wt^R u$ which is 
called total weight. These two grading operators 
can also be obtained from $\nu_1((\mathbf{0}, (a,\mathbf{0})))$
as follow: 
\bea  \label{d-L-R-nu}
\mathbf{d}^L &=& \frac{\partial}{\partial a} \lbar_{a=1} 
\nu_1((\mathbf{0}, (a,\mathbf{0}))), \nn
\mathbf{d}^R &=& \frac{\partial}{\partial \bar{a}} 
\lbar_{a=1} \nu_1((\mathbf{0}, (a,\mathbf{0}))).
\eea
Conversely, using $\mathbf{d}^L$ and $\mathbf{d}^R$, 
we can also express the action of 
$\nu_1((\mathbf{0},(a,\mathbf{0})))$ on $F$ as 
\beq  \label{nu-1-d}
\nu_1((\mathbf{0},(a,\mathbf{0}))) = 
a^{-\mathbf{d}^L} \bar{a}^{-\mathbf{d}^R}. 
\eeq

By the definition of smooth function on $K$, the correlation functions are
linear combination of 
$\prod_i (a_0^{(i)})^{r_i}\overline{a_0^{(i)}}^{s_i}$ where $r_i, s_i \in \R$.
Therefore, by (\ref{nu-1-d}), $F$ must be $\R\times \R$-graded.

\begin{rema}
{\rm 
When we define the smoothness, we only allow 
the real powers of $a_0^{(i)}$ and $\bar{a}_0^{(i)}$ to 
appear exactly for the sake of restricting $F$ to 
a $\R\times \R$-grading instead of a $\C\otimes \C$-grading.  
This is a physically natural condition because 
the operator $\mathbf{d}^L+\mathbf{d}^R$ 
is the Hamiltonian in physics. 
It can only has positive real eigenvalues 
in a unitary theory. If the total weights are real, then by Lemma 
\ref{double-grad}, both the left weights and the 
right weights must be real as well. 
}
\end{rema}

The image of $(\mathbf{0})\in \widehat{K}(0)$ 
under the morphism $\nu$ gives arise to
a very special element $\one \in \overline{F}$. 
Namely, 
\beq  \label{one-K}
\one :=\nu_0((\mathbf{0})).
\eeq
We call $\one$ the vacuum state. We have 
$(\mathbf{0},(a,\mathbf{0})) _{^1}\infty_{^0} (\mathbf{0}) 
= (\mathbf{0})$ 
for all $a\in \C^{\times}$. This implies $\one \in F_{(0,0)}$. 
Since $\widehat{K}(0) = \{ (\mathbf{0}) \}$, $W$ is just $\C \one$.

Let $P(z)=(z; \mathbf{0},(1,\mathbf{0}),(1,\mathbf{0}))
\in \widehat{K}(2)$. 
We denote the linear map 
$$
\nu_2(P(z)): F \otimes F \rightarrow 
\overline{F}
$$
as $\mathbb{Y}(\cdot; z, \bar{z}) \cdot$. 
Here we use both $z$ and its complex conjugation 
$\bar{z}$ in $\mathbb{Y}$ to emphasis that $\nu$ 
is not holomorphic in general. 
Because
$$
P(z)
_{^1}\infty_{^0} (\mathbf{0})= (\mathbf{0}, (1, \mathbf{0})),
$$
we have
\beq  \label{vacuum}
\mathbb{Y}(\one; z, \bar{z}) = I_{F}. 
\eeq
The equation (\ref{vacuum}) is called {\it vacuum property}. 
Moreover, since
$$
\lim_{z\rightarrow 0} P(z) _{^2}\infty_{^0} (\mathbf{0})=
(\mathbf{0},(1,\mathbf{0})) = I_K
$$ 
and $\nu$ maps identity to identity, we have 
\beq   \label{creat}
\lim_{z\rightarrow 0} \mathbb{Y}(\one; z, \bar{z}) = I_{F},
\eeq
which is called {\it creation property}.

\begin{prop} 
For $a, z\in \C^{\times}$ and $u\in F$, we have 
\beq   \label{d-l-d-r-conj}
a^{\mathbf{d}^L}\bar{a}^{\mathbf{d}^R} \mathbb{Y}(u; z, \bar{z}) 
a^{-\mathbf{d}^L}\bar{a}^{-\mathbf{d}^R} =
\mathbb{Y}(a^{\mathbf{d}^L}\bar{a}^{\mathbf{d}^R}u; az, \bar{a}\bar{z}).
\eeq
\end{prop}
\pf
First, for $u,v\in F$, we have 
\bea
&& \nu_2 ( (z; \mathbf{0}, (a_1, \mathbf{0}), (a_2, \mathbf{0}) )
(u\otimes v)  \nn
&&\hspace{1cm}= \nu_2((P(z) _{^1}\infty_{^0} 
(\mathbf{0}, (a_1, \mathbf{0}))) _{^2}\infty_{^0} 
(\mathbf{0}, (a_2, \mathbf{0}))) (u\otimes v)
\nn
&&\hspace{1cm}= ( \nu_2(P(z)) \, _{^1}\ast_{^0} 
\nu_1((\mathbf{0}, (a_1,\mathbf{0}))) )  
_{^2}\ast_{^0} \nu_1((\mathbf{0}, (a_2, \mathbf{0}))) (u\otimes v)\nn
&&\hspace{1cm}= 
\mathbb{Y}(a_1^{-\mathbf{d}^L}\bar{a}_1^{-\mathbf{d}^R} u ; z, \bar{z}) \,
a_2^{-\mathbf{d}^L}\bar{a}_2^{-\mathbf{d}^R}v. \label{d-conj-equ-1}
\eea
On the other hand, we also have 
\beq  \label{preimage-nu}
(\mathbf{0}, (a^{-1},\mathbf{0})) _{^1}\infty_{^0} P(z) = 
(az; \mathbf{0}, (a^{-1}, \mathbf{0}), (a^{-1}, \mathbf{0})).
\eeq
Using (\ref{d-conj-equ-1}), the image of 
(\ref{preimage-nu}) under the morphism $\nu$ gives 
$$
a^{\mathbf{d}^L}\bar{a}^{\mathbf{d}^R} \mathbb{Y}(u; z, \bar{z}) 
= \mathbb{Y}(a^{\mathbf{d}^L}\bar{a}^{\mathbf{d}^R} u ; az, \bar{a}\bar{z})
\, a^{\mathbf{d}^L}\bar{a}^{\mathbf{d}^R},
$$
which implies (\ref{d-l-d-r-conj}).
\epf

\begin{cor} \label{d-l-r-proof} 
For $z\in \C^{\times}$ and $u\in F$, we have
\bea
\left[ \mathbf{d}^L, \mathbb{Y}(u; z,\bar{z})\right] 
&=& z\frac{\partial}{\partial z} \mathbb{Y}(u; z, \bar{z})
+\mathbb{Y}(\mathbf{d}^L u; z,\bar{z})     \label{d-l} \\
\left[ \mathbf{d}^R, \mathbb{Y}(u; z, \bar{z})\right]
&=& \bar{z}\frac{\partial}{\partial \bar{z}} \mathbb{Y}(u; z,\bar{z})
+ \mathbb{Y}(\mathbf{d}^R u; z, \bar{z}), \label{d-r}
\eea
\end{cor}
\pf
Replace $a$ and $\bar{a}$ in (\ref{d-l-d-r-conj}) by $e^s$ and $e^{\bar{s}}$ 
respectively, and then take derivative of 
$\frac{\partial}{\partial s}|_{s=0}$
and $\frac{\partial}{\partial \bar{s}}|_{s=0}$ on both sides of 
(\ref{d-l-d-r-conj}), we immediate obtain (\ref{d-l}) and (\ref{d-r}). 
\epf

We further define two operators 
$D^L, D^R\in \hom(F, \overline{F})$ to be 
\bea \label{D-L-R-const}
D^L &:=& -\frac{\partial}{\partial a} \lbar_{a=0} 
\nu_1( (A(a; 1), (1, \mathbf{0})), \nn
D^R &:=& -\frac{\partial}{\partial \bar{a}} \lbar_{a=0} 
\nu_1( (A(a; 1), (1,\mathbf{0})).
\eea

\begin{prop}
$[\mathbf{d}^L, D^L] = D^L, [\mathbf{d}^R, D^R] = D^R$, 
$[\mathbf{d}^L, D^R] = [\mathbf{d}^R, D^L] =0$ and $[D^L, D^R]=0$.
As a consequence $D^R, D^L\in \edo F$
and have weights $(1,0)$ and $(0,1)$ respectively.  
\end{prop}

\pf
First we prove $[\mathbf{d}^L, D^L] = D^L$. 
Consider the sewing identity:
$$
(A(a;1), (a_1, \mathbf{0})) _{^1}\infty_{^0} (A(b;1), (b_1,\mathbf{0})) 
= (A(a+b/a_1;1), (a_1b_1, \mathbf{0}) ).
$$
It implies the following  identity: 
\bea  
&&\langle w', (\nu_1((A(a;1), (1,\mathbf{0}))) 
_{^1}\ast_{^0} \nu_1((\mathbf{0}, (a_1, \mathbf{0}))) \nn
&&\hspace{2cm} _{^1}\ast_{^0}
(\nu_1((A(b;1),(1,\mathbf{0})))  _{^1}\ast_{^0} 
\nu_1( (\mathbf{0}, (b_1, \mathbf{0})) ) (w) \rangle  \nn
&&\hspace{1cm} = \langle w', \nu_1((A(a+b/a_1;1), (1,\mathbf{0})))  
_{^1}\ast_{^0} \nu_1((\mathbf{0}, (a_1b_1, \mathbf{0})))(w)
\rangle  \label{d-D-comm-equ-1} 
\eea
for all $w\in F, w'\in F'$. Apply 
$$
\left( -\frac{\partial}{\partial a_1}\right) \lbar_{a_1=1}
\left( -\frac{\partial}{\partial b}\right) \lbar_{b=0,a=0,b_1=1}
-
\left( -\frac{\partial}{\partial a}\right) \lbar_{a=0}
\left( -\frac{\partial}{\partial b_1}\right) \lbar_{b=0,a_1=1,b_1=1}
$$
to both sides of equation (\ref{d-D-comm-equ-1}). The left hand side of 
(\ref{d-D-comm-equ-1}) gives
$$
\sum_{(m,n)\in I} \langle w', (\mathbf{d}^L P_{(m,n)} D^L - 
D^L \mathbf{d}^L) w\rangle,
$$
while the right hand of (\ref{d-D-comm-equ-1}) gives
\bea
&&\hspace{-0.5cm}
\left( -\frac{\partial}{\partial a_1}\right) \lbar_{a_1=1}
\left( -\frac{\partial}{\partial b}\right) \lbar_{b=0} 
\langle w', \nu_1((A(b/a_1;1) , (1,\mathbf{0}))) 
a_1^{-\mathbf{d}^L}\bar{a}_1^{-\mathbf{d}^R} w\rangle  \nn
&&\hspace{1.5cm} - \left( -\frac{\partial}{\partial a}\right) \lbar_{a=0}
\left( -\frac{\partial}{\partial b_1}\right) \lbar_{b_1=1}
\langle w', \nu_1((A(a;1), (1,\mathbf{0}))) b_1^{-\mathbf{d}^L}
\bar{b}_1^{-\mathbf{d}^R} w\rangle  \nn
&&\hspace{1cm} = \left( -\frac{\partial}{\partial a_1}\right) 
\langle w', D^L a_1^{-\mathbf{d}^L-1}
\bar{a}_1^{-\mathbf{d}^R} w\rangle \lbar_{a_1=1} 
- \langle w', D^L \mathbf{d}^L w\rangle
\nn
&&\hspace{1cm} = \langle w', D^L w\rangle.
\eea
Compare above two equations. We obtain
$[\mathbf{d}^L, D^L] = D^L$.

Second we show $[\mathbf{d}^R, D^L]=0$.  Apply 
$$
\left( -\frac{\partial}{\partial \bar{a}_1}\right) \lbar_{a_1=1}
\left( -\frac{\partial}{\partial b }\right) \lbar_{b=0,a=0,b_1=1}
-
\left( -\frac{\partial}{\partial a }\right) \lbar_{a=0}
\left( -\frac{\partial}{\partial \bar{b}_1}\right) 
\lbar_{b=0,a_1=1,b_1=1}
$$
to both sides of (\ref{d-D-comm-equ-1}). The left hand side of 
(\ref{d-D-comm-equ-1}) gives
$$
\sum_{(m,n)\in I} \langle w', (\mathbf{d}^R P_{(m,n)} D^L - 
D^L \mathbf{d}^R) w\rangle,
$$
while the right hand of (\ref{d-D-comm-equ-1}) gives
$$
\left( -\frac{\partial}{\partial a_1}\right) \lbar_{a_1=1} 
\langle w', D^L \frac{1}{a_1}a_1^{-\mathbf{d}^L}
\bar{a}_1^{-\mathbf{d}^R}  
w\rangle - 
\langle w', D^L \mathbf{d}^R w\rangle  = 0. 
$$
Compare above two equations. 
It is clear that $[\mathbf{d}^R, D^L]=0$. 

Combining above two results, we conclude that
$D^L\in \text{End} V$ and has weight $(1,0)$.
Similarly, we can show $[\mathbf{d}^L, D^R]=0$ and
$[\mathbf{d}^R, D^R]=D^R$. Therefore, $D^R$ is in 
$\text{End} V$ as well and has weight $(0,1)$.

Next we show that $[D^L, D^R]=0$. The identity 
\bea
(A(a+b;1), (1,\mathbf{0})) &=& 
(A(a;1), (1,\mathbf{0})) _{^1}\infty_{^0} (A(b;1), (1,\mathbf{0})) \nn
&=& (A(b;1), (1,\mathbf{0})) 
_{^1}\infty_{^0} (A(a;1), (1,\mathbf{0}))  
\eea
in $\widehat{K}$ implies that for $w', w\in F$,
\bea  \label{D-D-comm-equ-1}
&&\langle w', \nu_1((A(a;1), (1,\mathbf{0}))) 
_{^1}\ast_{^0} \nu_1((A(b;1), (1,\mathbf{0}))) w\rangle  \nn
&&\hspace{1.5cm} =
\langle w', \nu_1((A(b;1),(1,\mathbf{0}))) _{^1}\ast_{^0} 
\nu_1((A(a;1), (1,\mathbf{0})))w \rangle.
\eea
Apply 
$$
\left( -\frac{\partial}{\partial a}\right) \lbar_{a=0} 
\left( -\frac{\partial}{\partial \bar{b}}\right) \lbar_{b=0}
$$  
to both sides of (\ref{D-D-comm-equ-1}). We obtain $[D^L, D^R]=0$.
\epf

\begin{prop} \label{D-l-r-b-d-proof}
The following $D^L$- and $D^R$-bracket-derivative formula 
\begin{eqnarray}
&&\left[ D^L, \mathbb{Y}(u; z,\bar{z})\right] = 
\mathbb{Y}(D^{L}u; z,\bar{z}) =
\frac{\partial}{\partial z}\mathbb{Y}(u; z, \bar{z}),\label{D-L-b-d} \\
&&\left[ D^R, \mathbb{Y}(u; z, \bar{z})\right] = 
\mathbb{Y}(D^{R}u; z, \bar{z}) =
\frac{\partial}{\partial \bar{z}}\mathbb{Y}(u; z, \bar{z}),
\label{D-R-b-d}
\end{eqnarray}
holds for all $u\in F$ and $z\in \C^{\times}$. 
\end{prop}
\pf
We have the following sewing identity:
\beq
P(a) _{^1}\infty_{^0} (A(-z;1), (1,\mathbf{0})) = P(a+z).
\eeq
Notice that 
\beq
\mathbb{Y}(D^Lu, z,\bar{z}) = \frac{\partial}{\partial z} \lbar_{z=0} 
\nu_2(P(a) _{^1}\infty_{^0}  (A(-z;1), (1,\mathbf{0}))), 
\eeq
and 
\beq
\frac{\partial}{\partial z}\lbar_{z=0} 
\mathbb{Y}(u, z+a, \bar{z}+\bar{a}) 
=\frac{\partial}{\partial z}\lbar_{z=0} \nu_2(P(z+a)).  
\eeq
Hence it is clear that 
\beq
\mathbb{Y}(D^L u, z,\bar{z}) = \frac{\partial}{\partial z} 
\mathbb{Y}(u, z, \bar{z}).
\eeq

Similarly, we can show that 
the $D^L$-bracket property follows from the following sewing identity:
$$
(A(-z;1), (1,\mathbf{0})) _{^1}\infty_{^0} 
(P(a) _{^2}\infty_{^0} (A(z;1), (1,\mathbf{0}))) = P(z+a).
$$
We omit the detail. 

The proof of $D^R$-bracket and $D^R$-derivative properties is
similar. 
\epf

\begin{prop} \label{K-hat-ffa-prop} 
A smooth $\widehat{K}$-algebra $(F, \nu)$ has a natural 
structure of grading-restricted 
$\R\times \R$-graded full field algebra. 
\end{prop}
\pf
We have already proved the $\R\times \R$-grading. The grading-restriction condition 
on $F$ are automatic by the definition of smooth $\widehat{K}$-algebra.

The correlation function maps:
\beq \label{m-n-K}
m_n(u_1, \dots, u_n; z_1, \bar{z}, \dots, z_n, \bar{z}_n)
\eeq
can be defined as 
\beq \label{con-ffa-K-1}
\nu_{n+1}( (z_1, \dots, z_n; \mathbf{0}, (1,\mathbf{0}), \dots,
(1,\mathbf{0}))) (u_1\otimes \dots \otimes u_n \otimes \one)
\eeq
if $z_i\neq 0$, and as
\beq  \label{con-ffa-K-2}
\nu_n ( (z_1, \dots, \widehat{z_i}, \dots, z_n; \mathbf{0}, (1,\mathbf{0}), \dots,
(1,\mathbf{0}))) (u_1\otimes \dots \widehat{u_i} \dots \otimes u_n \otimes u_i)
\eeq
if $z_i=0$. By the definition of $\one$ given in (\ref{one-K}) and the fact that
\bea  \label{lim-z-0-K}
&&\lim_{z_i\rightarrow 0} (z_1, \dots, z_n; \mathbf{0}, (1,\mathbf{0}), \dots,
(1,\mathbf{0})) _{^i}\infty_{^0} (\mathbf{0}) \nn
&&\hspace{2cm}= 
(z_1, \dots, \widehat{z_i}, \dots, z_n; \mathbf{0}, (1,\mathbf{0}), \dots,
(1,\mathbf{0})), 
\eea
we see that (\ref{m-n-K}) are smooth with respect to $z_1, \dots, z_n$. 
Moreover, we also have 
$m_1(u; 0,0)= \nu_1( (1, \mathbf{0}) )(u) = I_F(u) = u$. Another 
identity property of full field algebra: 
\bea
&&m_{n+1}(u_1, \dots, u_n, \one; z_1, \bar{z}_1, \dots, z_{n+1},\bar{z}_{n+1})\nn
&&\hspace{3cm} = m_n(u_1, \dots, u_n; z_1, \bar{z}_1, \dots, z_n, \bar{z}_n).
\eea
also holds because the following identity: 
$$
Q _{^i} \infty _{^0} (\mathbf{0}) = (z_1, \dots, \widehat{z_i}, \dots, z_{n-1}; 
A(a;1), (a_0^{(1)}, \mathbf{0}), \dots, \widehat{(a_0^{(i)}, \mathbf{0})}, 
\dots, (a_0^{(n)}, \mathbf{0}) ), 
$$
where $Q$ is of form (\ref{K-hat-ele}), holds in $\widehat{K}$.

For the convergence property of 
$\R\times \R$-graded full field algebra, we use the weaker version
of convergence property discussed in the remark 1.2 in \cite{HK2}. 
Then it is clear that this weaker version of 
convergence property is automatically true
by the definition of algebra over partial operad. 
The permutation axiom of full field algebra follows automatically
from that of partial operad. 

The single-valuedness property follows from Lemma \ref{double-grad}. 
The rest of axioms follows from  (\ref{d-l-d-r-conj}), (\ref{d-l}),
(\ref{d-r}), (\ref{D-L-b-d}) and (\ref{D-R-b-d}). 
\epf

\begin{lemma}
Given a grading-restricted $\R\times \R$-graded full field algebra $F$, 
we have
\beq  \label{m-1-D}
m_1(u;z,\bar{z}) = e^{zD^L+\bar{z}D^R}u, 
\eeq
for $u\in F$, and 
\bea \label{D-m-n}
&&e^{aD^L+\bar{a}D^R} m_n(u_1, \dots, u_n; z_1, \bar{z}_1, \dots, z_n, \bar{z}_n) \nn
&&\hspace{1cm} = m_n (u_1, \dots, u_n; z_1+a, \overline{z_1+a}, 
\dots, z_n+a, \overline{z_n+a}).
\eea
for all $(z_1, \dots, z_n) \in M^{n}$ and $a\in \C$. 
\end{lemma}
\pf
Recall the following formula 
$$
\mathbb{Y}(u;z,\bar{z})\one = e^{zD^L+\bar{z}D^R}u
$$
which is proved in \cite{HK2}. Then we have 
$$
m_1(u;z,\bar{z}) = m_2(u,\one, z, \bar{z}, 0,0) = \mathbb{Y}(u; z, \bar{z})\one
=e^{zD^L+\bar{z}D^R}u.
$$
Thus we have proved (\ref{m-1-D}). (\ref{D-m-n}) can be proved as follow: 
\bea   \label{m-tilde-m}
&&e^{aD^L+\bar{a}D^R} m_n(u_1, \dots, u_n; z_1, \bar{z}_1, \dots, z_n, \bar{z}_n) \nn
&&\hspace{1cm} = 
m_1(m_n(u_1, \dots, u_n; z_1, \bar{z}_1, \dots, z_n, \bar{z}_n), a) \nn
&&\hspace{1cm} = 
m_n(u_1, 
\dots, u_n; z_1+a, \overline{z_1+a}, \dots, z_n+a, \overline{z_n+a}). 
\eea
for all $(z_1, \dots, z_n)\in M^n$ and $a\in \C$. 
\epf

\begin{thm}  \label{K-hat-ffa-iso-thm}
The category of grading-restricted $\R\times \R$-graded 
full field algebras is isomorphic to the category of smooth 
$\widehat{K}$-algebras.
\end{thm}
\pf
The proof is similar to that of Theorem 5.4.5 in \cite{H4}. 
Our case is much simpler. 

Given a grading-restricted 
$\R\times \R$-graded full field algebra 
$$(F,m,\one, \mathbf{d}^L, \mathbf{d}^R, D^L, D^R).$$ 
We define a map $\nu$ from $\widehat{K}$ to 
$\mathcal{H}_{F\, \C\one}^{\C^{\times}}$ as follow. 
We define
\beq  \label{nu-0-one}
\nu_0( (\mathbf{0} ) ) := \one.
\eeq
For $n>0$, $u_i\in F, a, z_i\in \C, i=1, \dots, n$ and 
$Q \in \widehat{K}(n)$ as (\ref{K-hat-ele}), we define
\bea \label{con-K-hat-alg}
&&\hspace{-1cm}\nu(Q)(u_1 \otimes \dots \otimes u_n) := 
e^{-aD^L-\bar{a}D^R}
m_n( (a_0^{(1)})^{-\mathbf{d}^L}(\overline{a_0^{(1)}})^{-\mathbf{d}^R} u_1, 
\dots, \nn
&&\hspace{3cm} (a_0^{(n)})^{-\mathbf{d}^L}(\overline{a_0^{(n)}})^{-\mathbf{d}^R} u_n; 
 z_1, \overline{z_1}, \dots, z_{n-1}, \overline{z_{n-1}}, 0, 0).
\eea
In particular, 
\bea
\nu_1( (\mathbf{0}, (a, \mathbf{0})) &=&
m_1((a_0^{(1)})^{-\mathbf{d}^L}(\overline{a_0^{(1)}})^{-\mathbf{d}^R} u_1; 0,0)\nn
&=& (a_0^{(1)})^{-\mathbf{d}^L}(\overline{a_0^{(1)}})^{-\mathbf{d}^R} u_1. 
\eea
This gives the representation of rescaling group $\C^{\times}$ 
according to the $\R\times \R$-grading of $F$. 
The $\R\times \R$-grading, 
together with the single-valuedness property, 
guarantees that the grading is in the set 
$\{ (m,n) | m,n\in \R, m-n\in \Z\}$ 
as required by the axioms of smooth $\widehat{K}$-algebra. 
The grading restriction conditions and 
the smoothness of all correlation functions are
all automatically satisfied.

Let $P \in \widehat{K}(m)$ and $Q\in \widehat{K}(n)$
given as follow:  
\bea
P &=& (z_1, \dots, z_{m-1}; A(a;1), (a_0^{(1)}, \mathbf{0}), \dots, 
(a_0^{(m)}, \mathbf{0}) ) \nn
Q &=& (\xi_1,\ldots, \xi_{n-1}; A(b;1), (b_0^{(1)}, \mathbf{0}), \dots, 
(b_0^{(n)}, \mathbf{0}) ). 
\eea
We assume that $P _{^i}\infty_{^0} Q$ exists. Then we have 
\bea
&&\hspace{-1cm}P _{^i}\infty_{^0} Q =
(z_1, \dots, z_{i-1}, \frac{\xi_1-a}{a_0^{(i)}}+z_i, \dots, 
\frac{\xi_m-a}{a_0^{(i)}} + z_i, z_{i+1}, \dots, z_n; \nn
&&\hspace{1cm} A(a;1), (a_0^{(1)}, \mathbf{0}), \dots, 
(a_0^{(i-1)}, \mathbf{0}),   (a_0^{(i)}b_0^{(1)}, \mathbf{0}), \dots,
(a_0^{(i)}b_0^{(m)}, \mathbf{0}),  \nn 
&&\hspace{5cm} (a_0^{(i-1)}, \mathbf{0}), \dots,
(a_0^{(n)}, \mathbf{0})\,\, ).
\eea
By (\ref{D-m-n}), we have
\bea
&&\hspace{-1cm} \nu_{m}(P)  _{^i}*_{^0} \nu_n(Q) \nn
&&\hspace{0cm} = \sum_{k,l\in \R} 
e^{-aD^L-\bar{a}D^R} 
m_m( (a_0^{(1)})^{-\mathbf{d}^L}(\overline{a_0^{(1)}})^{-\mathbf{d}^R}u_1, \dots, \nn 
&&\hspace{1cm} (a_0^{(i)})^{-\mathbf{d}^L}(\overline{a_0^{(i)}})^{-\mathbf{d}^R}P_{k,l}
m_n((b_0^{(1)})^{-\mathbf{d}^L}(\overline{b_0^{(1)}})^{-\mathbf{d}^R}v_1, \dots, 
(b_0^{(n)})^{-\mathbf{d}^L}(\overline{b_0^{(n)}})^{-\mathbf{d}^R}v_n;  \nn
&&\hspace{3cm} \xi_1-b, \overline{\xi_1-b}, \dots, \xi_m-b, \overline{\xi_m-b}), \nn
&&\hspace{4cm} \dots, 
(a_0^{(m)})^{-\mathbf{d}^L}(\overline{a_0^{(m)}})^{-\mathbf{d}^R}u_m; 
z_1,\bar{z}_1, \dots, z_m,\bar{z}_m)  \nn
&&\hspace{0cm}= \sum_{k,l\in \R}  e^{-aD^L-\bar{a}D^R} 
m_m( (a_0^{(1)})^{-\mathbf{d}^L}(\overline{a_0^{(1)}})^{-\mathbf{d}^R}u_1, \dots, \nn 
&&\hspace{1cm} P_{k,l} m_n((a_0^{(i)}b_0^{(1)})^{-\mathbf{d}^L}
(\overline{a_0^{(i)}b_0^{(1)}})^{-\mathbf{d}^R}v_1, \dots, 
(a_0^{(i)}b_0^{(n)})^{-\mathbf{d}^L}(\overline{a_0^{(i)}b_0^{(n)}})^{-\mathbf{d}^R}v_n;  \nn
&&\hspace{1cm} (\xi_1-b)/a_0^{(i)}, \overline{(\xi_1-b)/a_0^{(i)}}, \dots, 
(\xi_m-b)/a_0^{(i)}, \overline{(\xi_m-b)/a_0^{(i)}}), \nn 
&&\hspace{4cm} \dots, 
(a_0^{(m)})^{-\mathbf{d}^L}(\overline{a_0^{(m)}})^{-\mathbf{d}^R}u_m; 
z_1,\bar{z}_1, \dots, z_m,\bar{z}_m)  \nn
&&\hspace{0cm}= e^{-aD^L-\bar{a}D^R}
m_{m+n-1}((a_0^{(1)})^{-\mathbf{d}^L}(\overline{a_0^{(1)}})^{-\mathbf{d}^R}u_1, \dots, 
(a_0^{(i)}b_0^{(1)})^{-\mathbf{d}^L}
(\overline{a_0^{(i)}b_0^{(1)}})^{-\mathbf{d}^R}v_1, \nn
&&\hspace{2cm} \dots, 
(b_0^{(n)})^{-\mathbf{d}^L}(\overline{b_0^{(n)}})^{-\mathbf{d}^R}v_n, \dots, 
(a_0^{(m)})^{-\mathbf{d}^L}(\overline{a_0^{(m)}})^{-\mathbf{d}^R}u_m; z_1,\bar{z}_1, 
\dots,\nn
&&\hspace{1.5cm} 
\frac{\xi_1-b}{a_0^{(i)}}+z_i, \overline{\frac{\xi_1-b}{a_0^{(i)}}+z_i}, \dots, 
\frac{\xi_m-b}{a_0^{(i)}}+z_i, \overline{\frac{\xi_m-b}{a_0^{(i)}}+z_i}, 
\dots, z_m,\bar{z}_m)  \nn
&&\hspace{0cm}= \nu_{m+n-1} ( P _{^i}\infty_{^0} Q ). \nonumber
\eea

We have checked all the axioms of smooth $\widehat{K}$-algebra. 
Therefore the triple $(F, \C\one, \nu)$ gives
a smooth $\widehat{K}$-algebra. Thus we obtain a functor
from the category of $\R\times \R$-graded full field algebras 
to that of smooth $\widehat{K}$-algebras. 

On the other hand, by Proposition \ref{K-hat-ffa-prop}, 
we also have a functor from the category of smooth 
$\widehat{K}$-algebras to the category of 
of grading-restricted $\R\times \R$-graded full field algebras. 

Now we show that these two functors give isomorphisms 
between two categories. 
Let $(F, \tilde{m}, \tilde{\one}, \tilde{D}^L, \tilde{D}^R)$ 
be full field algebra
obtained according to Proposition (\ref{K-hat-ffa-prop}) 
from a smooth $\widehat{K}$-algebra,
which is further obtained from 
a grading-restricted $\R\times \R$-graded full field algebra 
$(F, m, \one, D^L, D^R)$ according to (\ref{con-K-hat-alg}).  
First $\tilde{\one}=\one$ is obvious by 
our constructions (\ref{nu-0-one}) and (\ref{one-K}).  

By construction (\ref{con-ffa-K-1}) and
(\ref{con-ffa-K-2}), for $z_i\neq 0, i=1, \dots, n$, we have
\bea
&&\tilde{m}_n(u_1, \dots u_n; z_1, \bar{z}_1, \dots, z_n, \bar{z}_n) \nn
&&\hspace{1cm} = \nu_{n+1}( (z_1, \dots, z_n; \mathbf{0}, (1,\mathbf{0}), \dots,
(1,\mathbf{0}))) (u_1\otimes u_n \otimes \one) \nn
&&\hspace{1cm} = m_{n+1}(u_1,\dots, u_n, \one; 
z_1, \bar{z}_1, \dots, z_n, \bar{z}_n, 0, 0) \nn
&&\hspace{1cm} = m_{n}(u_1, \dots,  u_n; 
z_1, \bar{z}_1, \dots, z_n, \bar{z}_n). 
\eea
The cases when $z_i=0$ for some $i=1, \dots, n$ follows from 
smoothness. Therefore $\tilde{m}=m$. 

By \ref{D-L-R-const},  we also have
\bea
\tilde{D}^L &=& -\frac{\partial}{\partial a} \lbar_{a=0} v_1((A(a;1), (1,\mathbf{0}) ) 
= -\frac{\partial}{\partial a} \lbar_{a=0} e^{-aD^L-\bar{a}D^R} = D^L,  \nn
\tilde{D}^R &=& -\frac{\partial}{\partial \bar{a}} \lbar_{a=0} 
v_1((A(a;1), (1,\mathbf{0}) )
= -\frac{\partial}{\partial a} \lbar_{a=0} e^{-aD^L-\bar{a}D^R} = D^R. 
\eea 
The coincidence of two gradings is also obvious. 
Therefore, we have proved that one way of 
composing two functors gives the identity functor on the category of 
grading restricted $\R\times \R$-graded full field algebra. 

Similarly, one can show that the opposite way of composing these two functors also
gives the identity functor on the category of smooth $\widehat{K}$-algebras. 
\epf

Let us now turn our attention to the conformal case. 
Using the result of Huang \cite{H4}, 
it is easy to see that the tensor product $V^L\otimes V^R$ of 
two vertex operator algebras $V^L$ and $V^R$ with central charge
$c^L$ and $c^R$ respectively, has a canonical 
structure of smooth $\tilde{K}^{c^L} 
\otimes \overline{\tilde{K}^{\overline{c^R}}}$-algebras.
If a smooth $\tilde{K}^{c^L} 
\otimes \overline{\tilde{K}^{\overline{c^R}}}$-algebras $(F, W, \nu)$
is equipped with
an embedding $\rho: V^L\otimes V^R \hookrightarrow F$ as 
a morphism of smooth $\tilde{K}^{c^L} 
\otimes \overline{\tilde{K}^{\overline{c^R}}}$-algebra, 
we call it a smooth $\tilde{K}^{c^L} 
\otimes \overline{\tilde{K}^{\overline{c^R}}}$-algebras over 
$V^L\otimes V^R$ and denote it as $(F, \nu, \rho)$.

We consider a conformal full field algebra $F$ 
over $V^L\otimes V^R$, denoted as a triple $(F, m,\rho)$ where
$\rho: V^L\otimes V^R \hookrightarrow F$ is a monomorphism 
of conformal full field algebra. 
By Assumption \ref{assumption}, 
$V^L, V^R$ satisfy the conditions in 
Theorem \ref{ioa}. In this case, the products and iterates 
of intertwining operators of $V^L, V^R$ 
satisfies very nice convergence and 
analytic extension properties \cite{H7}. As a consequence, 
$m_n$ also have certain analytic extension properties \cite{HK2}. 
Namely, for $u_1, \dots, u_n\in F, w'\in F'$, 
$$
\langle w', \mathbb{Y}(u_1;z_1, \zeta_1) \dots 
\mathbb{Y}(u_n; z_n, \zeta_n)\one\rangle
$$
is absolutely convergent when $|z_1|>\dots >|z_n|>0$
and $|\zeta_1|>\dots > |\zeta_n|>0$ and can be analytically
extended to a multivalued analytic function in 
the region given by $z_i\neq z_j$, $z_i\neq 0$,
$\zeta_i\neq \zeta_j$, $\zeta_i\neq 0$. We use
$$
E(m)_n(w', u_1,\dots, u_n; z_1, \zeta_1, \dots, z_n, \zeta_n)
$$
to denote this function.

For the simplicity of notation, we will not distinguish
$L^{L}(n)$ with $L^{L}(n)\otimes 1$ and $L^{R}(n)$ with
$1\otimes L^R(n)$ in this work.

\begin{prop}  \label{cffa-K-alg-prop}
{\rm 
A conformal full field algebra over $V^L\otimes V^R$, 
$(F, m, \rho)$,
has a canonical structure of smooth 
$\tilde{K}^{c^L} \otimes \overline{\tilde{K}^{\overline{c^R}}}$-algebras
over $V^L\otimes V^R$.
}
\end{prop}
\pf
We first define a map $\nu_n:  \tilde{K}^{c^L}\otimes 
\overline{\tilde{K}^{\overline{c^R}}}(n) 
\rightarrow \hom (F^{\otimes n}, \overline{F})$ for $n\in \N$. 
For $(A^{(0)})\in K(0)$, we define
\bea  \label{nu-0-L-one}
\nu_0( \psi\otimes \bar{\psi} ((A^{(0)})) ): 
= e^{-L_-^L(A^{(0)}) -L_-^R(A^{(0)}) }\one. 
\eea

We will define other cases indirectly. 
Recall that we choose the canonical representative
of $Q\in K(n)$ to have $0$-th puncture at $\infty$, 
$n$-th puncture at $0$ and satisfy the condition 
(\ref{can-rep-cond}). 
We relax the condition a little. We call 
all the representatives of an equivalent class $Q$ 
with $0$-th puncture sitting at $\infty$ and 
satisfying the condition (\ref{can-rep-cond}) 
as quasi-canonical representatives. 
The set of quasi-canonical representatives are clearly
parametrized by $z_n$ the location of the $n$-th puncture.  
When $z_n=0$, it is nothing but the canonical representative.
In order to distinguish the notation from those of canonical
representative (\ref{ele-K-n}), we use 
\beq \label{quasi-can-rep}
( (\infty; 1, A^{(0)}); 
(z_1; a_0^{(0)}, A^{(1)}), \dots, (z_n; a_0^{(n)}, A^{(n)}) )
\eeq
to denote a quasi-canonical representative of a sphere
with tubes with punctures at 
$\infty$ and $z_1, \dots, z_n\in \C, z_i\neq z_j$
with local coordinate maps giving by (\ref{f-0}) and 
(\ref{f-i}) respectively, and $z_n\neq 0$ in general. 

For $\tilde{Q}=\lambda\psi\otimes \bar{\psi}(Q) 
\in \tilde{K}^{c^L}\otimes \overline{\tilde{K}^{\overline{c^R}}}(n), 
n\geq 1$ where $Q\in K(n)$ given by a quasi-canonical representative
of form (\ref{quasi-can-rep}), we define $\nu_n$, the 
restriction of $\nu$ on 
$\tilde{K}^{c^L} \otimes \overline{\tilde{K}^{\overline{c^R}}}(n)$, 
as follow:
\bea
&&\hspace{-0.5cm}
\nu_n(\tilde{Q}): = \lambda e^{-L_-^L(A^{(0)}) - L_-^R(\overline{A^{(0)}})}  
m_{n}( e^{-L_+^L(A^{(1)}) - L_+^R( \overline{A^{(1)} } ) } 
(a_0^{(1)})^{-L^L(0)} \overline{a_0^{(1)}}^{-L^R(0)} u_1,  \dots, \nn
&&\hspace{1cm}  
e^{-L_+^L(A^{(n)}) - L_+^R( \overline{A^{(n)} })} 
(a_0^{(1)})^{-L^L(0)} \overline{a_0^{(1)}}^{-L^R(0)}u_n;  
z_1, \bar{z}_1, \dots, z_{n-1}, \bar{z}_{n-1}, z_n,\bar{z}_n),  \nn
\label{nu-cffa-K-alg}
\eea
where 
\beq
L_{\pm}^L(A) = \sum_{j=1}^{\infty} A_j L^L(\pm j), \quad 
L_{\pm}^R(A) = \sum_{j=1}^{\infty} A_j L^R(\pm j)  \nonumber
\eeq
for $A\in \prod_{n\in \N} \C$. 

Of course, we have to show that $\nu$ is well-defined on 
$\tilde{K}^{c^L} \otimes \overline{\tilde{K}^{\overline{c^R}}}$. 
It follows immediately from 
the identity (\ref{D-m-n}) and Huang's proof of 
case 1 (as we will recall later) in the Proof of the sewing axiom in 
Proposition 5.4.1 in \cite{H4}\footnote{Actually, it is only a
trivial case of the case 1 in Huang's proof of Proposition 5.4.1 
in \cite{H4}, namely when 
$P$ and $Q$ in (\ref{case-1-P-Q}) is so that
$a_0=1, A^{(1)}=B^{(0)}=\mathbf{0}$.}.

The advantage of this way of 
defining $\nu$ is that we see immdiately that
$\nu_n$ is $S_n$-equivariant because of the permutation 
axiom of full field algebra. Another advantage of this
definition is that we can always assume that $z_i, i=1, \dots, n$
have distinct absolute values. Indeed, 
if for the quasi-canonical representative we start with 
there are some $z_k$ having same absolute values, then
we can always choose another quasi-canonical representative
such that $|z_i|\neq |z_j|$ for $i\neq j$. This fact will
be useful later.

We will show that such defined $(F, \nu, \rho)$ gives
a smooth 
$\tilde{K}^{c^L} \otimes \overline{\tilde{K}^{\overline{c^R}}}$-algebra
over $V^L\otimes V^R$. 

First of all, by the construction (\ref{nu-cffa-K-alg}), we have
$$
\nu_1((\mathbf{0}, (a, \mathbf{0}))) = a^{-L^L(0)}\bar{a}^{-L^R(0)}.
$$
Since $L^L(0)$ and $L^R(0)$ are two grading operators on $F$, 
$F$ is exactly graded by inequivalent irreducible 
representations of the rescaling group $\C^{\times}$. 
The grading restriction conditions are automatic.

From the construction (\ref{nu-cffa-K-alg}), 
$\nu_1((\mathbf{0}, (1, \mathbf{0})))(u) = m_1(u; 0,0)=u$. 
Thus $\nu$ maps the identity $I_K$ to the identity $I_F$. 
We also see that
$\nu_0(\tilde{K}^{c^L}\otimes \overline{\tilde{K}^{\overline{c^R}}}(0))$
is the subspace of $F$ generated by the actions of 
$\{ L^L(n), L^R(n) \}_{n\in \Z}$ on 
$\one^L\otimes \one^R$. We denote
this space as $\langle \omega^L, \omega^R\rangle$.

Notice that (\ref{nu-cffa-K-alg}) is completely 
compactible with the smooth 
$\tilde{K}^{c^L} \otimes \overline{\tilde{K}^{\overline{c^R}}}$-algebra
structure on $V^L\otimes V^R$. In other words, if 
$(F, \langle \omega^L, \omega^R\rangle, \nu)$ indeed gives a
smooth 
$\tilde{K}^{c^L} \otimes \overline{\tilde{K}^{\overline{c^R}}}$-algebra,
then $\rho$ must be an embedding of $V^L\otimes V^R$ into $F$
as a smooth 
$\tilde{K}^{c^L} \otimes \overline{\tilde{K}^{\overline{c^R}}}$-algebra
monomorphism.

It remains to prove the sewing properties of $\nu$, i.e.
\beq  \label{sew-nu}
\nu_{m+n-1} (\tilde{P} _{^i}\widetilde{\infty}_{^0} \tilde{Q}) =
\nu_m(\tilde{P}) _{^i}*_{^0} \nu_n(\tilde{Q}).
\eeq
for all $\tilde{P} \in \tilde{K}^{c^L} \otimes \overline{\tilde{K}^{\overline{c^R}}}(m)$ and $\tilde{Q}\in \tilde{K}^{c^L} \otimes \overline{\tilde{K}^{\overline{c^R}}}(n)$ as long as 
$\tilde{P} _{^i}\widetilde{\infty}_{^0} \tilde{Q}$ exists. 
The proof is essentially same as that of the sewing axiom in 
Proposition 5.4.1 in \cite{H4}. 
We assume that readers are familar with the proof of 
Proposition 5.4.1 in \cite{H4}. 
We will only point out where the differences are.

In the rest of the proof, 
we always use $\tilde{P}, \tilde{Q}$ to denote 
$\psi\otimes \bar{\psi}(P), \psi\otimes \bar{\psi}(Q)$
respectively 
for $P\in K(m), Q\in K(n), m\in \Z_+, n\in \N$. 

In the proof of Proposition 5.4.1 in \cite{H4}, 
the sewing axiom is proved in cases. 
The first case is when $m=n=1$. Let 
\bea
P &:=& ( A^{(0)}, (a_0, A^{(1)}) ),  \nn
Q &:=& ( B^{(0)}, (b_0, B^{(1)}) ). \label{case-1-P-Q}
\eea
Assume that $P _{^1}\infty_{^0} Q$ exists. In this case, only 
Virasoro algebras are involved. Since two Virasoro algebras
generated by $\{ L^L(n)\}_{n\in \Z}$ and $\{ L^R(n)\}_{n\in \Z}$
are mutually commutative, 
we can study these two Virasoro algebras seperately. 
The left Virasoro algebra is completely same 
as the Virasoro algebra in the proof of Proposition 5.4.1 
in \cite{H4}. 
For the right Virasoro algebra, we need one more piece of fact. 
Let $\mathcal{A}=\{ \mathcal{A}_1, \dots \}$ and
$\mathcal{B}=\{ \mathcal{B}_1, \dots, \}$ be two sequences of 
formal variables. It was proved in \cite{H4} that, 
\beq  \label{sew-equ-1}
e^{-L_-^R(\mathcal{A})} \alpha^{-L^L(0)} e^{L_-^R (\mathcal{B})}
= e^{L_-^R(\Psi_-)} e^{L_+^R(\Psi_+)} e^{\Psi_0 L^R(0)} \alpha^{-L^R(0)} 
e^{\Gamma(\mathcal{A}, \mathcal{B}, \alpha)c^R},
\eeq
where 
$\Psi_{\pm,0}=\Psi_{\pm,0}(\mathcal{A}, \mathcal{B}, \alpha) \in 
\C [\alpha, \alpha^{-1}][[\mathcal{A}, \mathcal{B}]]$.
From Huang's study of $\Psi_{\pm,0}$, 
it is easy to see that they are actually 
in $\Q [\alpha, \alpha^{-1}][[\mathcal{A}, \mathcal{B}]]$. 
Hence we have 
\beq  \label{Psi-pm-real}
\overline{\Psi_{\pm,0}(A^{(1)}, B^{(0)}, a)} = 
\Psi_{\pm,0}(\overline{A^{(1)}}, \overline{B^{(0)}}, \bar{a} ). 
\eeq
Using this fact, we see that (\ref{sew-nu}) holds in 
this case. 

Moreover, let $P_t:= ( A^{(0)}, (t^{-1}a_0, A^{(1)}) )$. 
Then $P_t \, _{^1}\infty_{^0} Q$ also exists for all $1\geq |t|>0$. 
We have $\nu_1(\tilde{P}_t) _{^1}*_{^0} \nu_1(\tilde{Q})$ equals to the 
following double series
\bea  \label{double-series}
&&(\nu_1(\tilde{P}) t_1^{L^L(0)}s_1^{L^R(0)}) _{^1}*_{^0} 
\nu_1(\tilde{Q}) \nn
&&\hspace{0.5cm}
=\sum_{m,n} (e^{-L_-^L(A^{(0)})} e^{-L_+^L(A^{(1)})} a^{-L^L(0)})
\otimes (e^{-L_-^R(\overline{A^{(0)})}} e^{-L_+^R(\overline{A^{(1)}})} 
\bar{a}^{-L^R(0)})  \nn
&&\hspace{1cm} P_{m,n} (e^{-L_-^L(B^{(0)})} e^{-L_+^L(B^{(1)})} b^{-L^L(0)} \otimes 
e^{-L_-^R(B^{(0)})} e^{-L_+^R(B^{(1)})} \bar{b}^{-L^R(0)}) t_1^m s_1^n \nn
\eea
when $t_1=t, s_1=\bar{t}$. 
Hence Huang's results in \cite{H4} implies that 
(\ref{double-series}) is absolutely convergent to a
multivalued analytic function of $t_1, s_1$ 
when $1\geq |t_1|, |s_1|>0$.

The second case is when $i=2$, 
$$
P = P(z), \quad Q = (B^{(0)}, (1, \mathbf{0})) 
$$
or $i=1$, 
$$
P= (\mathbf{0}, (a_0^{(1)}, A^{(1)})), \quad  Q=P(z). 
$$
By the definition of $\nu$ in (\ref{nu-cffa-K-alg}), 
$\nu_2(P(z)) = \mathbb{Y}(\cdot; z, \bar{z})\cdot$. In our case, 
$\mathbb{Y}$ is an intertwining operator for $V^L\otimes V^R$
viewed as vertex operator algebra. It is proved in \cite{HK2} that
$\mathbb{Y}$ split as follow: 
\beq
\mathbb{Y} = \sum_{i=1}^{N} \Y_i^L \otimes \Y_i^R, \quad 
\mbox{for some $N\in \Z_+$},
\eeq
where $\Y_i^L$ and $\Y_i^R$ are intertwining operators for 
$V^L$ and $V^R$ respectively. 
Using this splitting property of $\mathbb{Y}$, we can 
again consider the left and right seperately. Moreover, 
it is harmless to only consider the case $N=1$ because
$\Y_i^L$ and $\Y_i^R$ do not change when $L^L(m), L^R(n)$ 
exchange positions with $\Y_i^L$ and $\Y_i^R$. 
Only the arguments of $\Y_i^L$ and $\Y_i^R$ have changed.
Then Huang's results in \cite{H4} immediately imply that
(\ref{sew-nu}) holds in this case.  
Moreover, 
since $F$ is a module over $V^L\otimes V^R$, the grading of $F$ is
of form $(S^L+\N, S^R +\N)$ where $S^L$ and $S^R$ are two 
subsets of $\R$ with finite cardinalities. 
Then for $u_1, u_2\in F$, 
$(\nu(\tilde{P}) t_1^{L^L(0)}s_1^{L^R(0)}) _{^i}*_{^0} \nu(\tilde{Q})
(u_1\otimes u_2)$
is of following form: 
$$
\sum_{j=1}^M  \sum_{m,n\in \N} c_{m,n}^{(j)}\, \, t_1^{m+r_1^{(j)}}s_1^{n+r_2^{(j)}}
$$
where $r_1^{(j)}\in S^L, r_2^{(j)}\in S^R$ for some $M\in \N$. 
We will called such series as 
generalized power series. Such series share the same 
convergence and analytic properties as ordinary series
except for multivaluedness. 
Then Huang's result also implies that the generalized power
series 
$(\nu(\tilde{P}) t_1^{L^L(0)}s_1^{L^R(0)}) _{^i}*_{^0} \nu(\tilde{Q})$,
is absolutely convergent for $1\geq|t_1|, |s_1|>0$.

The third case is when $i=1, 
P=P(z)$ and $Q=(B^{(0)}, (1, \mathbf{0}))$.  
This case  is proved in \cite{H4}
by using the skew-symmetry of vertex operator algebra
to reduce it to the second case.  
For us, a similar skew-symmetry ((1.41) in \cite{HK2})
holds for conformal full field algebra. Hence we  
can again reduce the third case to the second case.

We do slightly differently in the fourth case.
Let $P$ and $Q$ be as follow
\bea
P &=& ( (\infty; 1,A^{(0)}); (z_1; a_0^{(1)}, A^{(1)}), 
\dots, (z_m; a_0^{(m)}, A^{(m)}) ),  \nn
Q &=& (z; \mathbf{0}, (b_0^{(1)}, B^{(1)}), (b_0^{(2)}, B^{(2)})). 
\eea
Notice that $P$ is quasi-canonical, while $Q$ is canonical. 
By the definition of $\nu$ in (\ref{nu-cffa-K-alg}), 
it is clear that 
$$
\nu_m(\tilde{P}) = \nu_m(\tilde{P}') _{^i}*_{^0} \nu_1 (\tilde{P_i}),
$$
where 
\bea
P' &=& 
((\infty; 1,A^{(0)}); (z_1; a_0^{(1)}, A^{(1)}), \dots, 
(z_i; 1, \mathbf{0}), 
\dots, (z_m; a_0^{(m)}, A^{(m)}) ), \nn
P_i &=& ( (\infty; 1, \mathbf{0}), (0; a_0^{(i)}, A^{(i)})).
\eea
Then $\nu_m(\tilde{P}) \, _{^i}*_{^0} \nu_2(\tilde{Q})$ equals to
the following iterate series
\bea  \label{preswitch-series}
&&\big( (\nu_m(\tilde{P}')t_1^{L^L(0)}s_1^{L^R(0)}) \, _{^i}*_{^0} 
\nu_1(\tilde{P_i}) t_2^{L^L(0)}s_2^{L^R(0)}\big) \, _{^i}*_{^0}
\nu_2(\tilde{Q})\nn
&&= \sum_{p_2, q_2} \left( \sum_{p_1, q_1} 
\big( (\nu_m(\tilde{P}') \, _{^i}*_{^0} 
P_{p_1,q_1} \nu_1(\tilde{P_i}) \big)  \, _{^i}*_{^0}
P_{p_2,q_2} \nu_2(\tilde{Q})  \, \,  
t_1^{p_1}s_1^{q_1} t_2^{p_2}s_2^{q_2} \right)
\eea
when $t_1=s_1=t_2=s_2=1$. We want to show that above iterate
series is absolutely convergent. 

We first consider a different iterate series
\beq \label{switch-series}
\sum_{p_1, q_1} \left( \sum_{p_2, q_2} 
\big( \nu_m(\tilde{P}') \, _{^i}*_{^0} 
P_{p_1,q_1} \nu_1(\tilde{P_i}) \big)  \, _{^i}*_{^0}
P_{p_2,q_2} \nu_2(\tilde{Q})  \, \,  t_1^{p_1}s_1^{q_1} t_2^{p_2}s_2^{q_2} 
\right) 
\eeq
which is obtained by switching order of multiple sum in 
(\ref{preswitch-series}). 
Since $\nu_1(\tilde{P}_i)\in \text{End} F$, we can apply
associativity for each term of the 
iterate sum (\ref{switch-series}). 
We obtain the following series:
\beq  \label{swith-series-equ-1}
\sum_{p_1, q_1} \left( \sum_{p_2, q_2} 
\nu_m(\tilde{P}') \, _{^i}*_{^0} 
P_{p_1,q_1} \big( \nu_1(\tilde{P_i})   \, _{^1}*_{^0}
P_{p_2,q_2} \nu_2(\tilde{Q}) \big)  \, \,  
t_1^{p_1}s_1^{q_1} t_2^{p_2}s_2^{q_2} \right).
\eeq

Solving the sewing equation for the case 
$P_i \, _{^1}\infty_{^0} Q$, we obtain that the 
canonical representative of 
$Q':=P_i \, _{^1}\infty_{^0} Q$ can be written in the 
following form
$$
(f_i^{-1}(z); \mathbf{0}, (c_0^{(1)}(a_0^{(i)}), C^{(1)}(a_0^{(i)})), 
(c_0^{(2)}(a_0^{(i)}), C^{(2)}(a_0^{(i)}))
$$
where $c_0^{(k)}(a_0^{(i)}), k=1,2$ and $C^{(k)}(a_0^{(i)}), k=1,2$ 
are analytic functions of 
$a_0^{(i)}$ valued in $\C^{\times}$ and $\C$ respectively, and 
$f_i$ is the local coordinate map at $0\in P_i$. 
Using the results proved in the case 2, we see that, for 
$u_1, u_2\in F$, the following series
\beq  \label{series-it-asso}
\sum_{p_2, q_2} \nu_1(\tilde{P_i})   \, _{^1}*_{^0}
P_{p_2,q_2} \nu_2(\tilde{Q})  (u_1\otimes u_2)  \, \,  t_2^{p_2}s_2^{q_2}
\eeq
is absolutely convergent when $1\geq |t_2|,|s_2|>0$ 
to a multivalued analytic function 
\beq
\mathbb{Y}(T_iu_1; f_i^{-1}(t_2z), 
\overline{f_i^{-1}}(s_2\bar{z}))T_{i+1}u_2
\eeq
where 
\beq  \label{T-ii+1}
T_k= e^{-L_+^L(C^{(k)}(a_0^{(i)}t_2^{-1})) 
-L_+^R(\overline{C^{(k)}}(\overline{a_0^{(i)}}s_2^{-1}))} 
(c_0^{(k)}(a_0^{(i)}t_2^{-1}))^{-L^L(0)}
(\overline{c_0^{(k)}}(\overline{a_0^{(i)}}s_2^{-1}) )^{-L^R(0)}
\eeq
for $k=i,i+1$. When $t_2=s_2=1$, the series 
(\ref{series-it-asso}) simply converges to 
$\nu_2(Q')$.

As we mentioned before, we can always 
assume $z_1, \dots, z_{m}$ to have distinct
absolute values by choosing a suitable quasi-canonical 
representative of $P'$.
In this case, $\nu_m(\tilde{P}')$ can 
be written as a product of $\mathbb{Y}$ which is 
an intertwining operator of $V^L\otimes V^R$.   
Since $f_i^{-1}$ map neighborhoods of $0$ to 
neighborhoods of $0$, for fixed $z$, 
$|f_i^{-1}(t_2z)|$ is sufficiently small as long as $|t_2|$ is 
sufficiently small. Similarly, 
$|\overline{f_i^{-1}}(s_2\bar{z})|$ is sufficiently small 
as long as $|s_2|$ is sufficiently small. 
Therefore, we can always find $r>0$ so that
$|z_j|-|z_i|>|f_i^{-1}(t_2z)|$ and 
$|z_j|-|z_i|>|\overline{f_i^{-1}}(s_2\bar{z})|$ 
for all $j\neq i$ and $r>|t_2|, |s_2|>0$.   
By the convergence property of intertwining operators
of $V^L, V^R$, it is clear that the following series,
for $u_1, \dots, u_{m+1}\in F$,  
\beq   \label{switch-series-equ-3}
\sum_{p_1, q_1} 
\nu_m(\tilde{P}') \, _{^i}*_{^0} P_{p_1, q_1} 
(\mathbb{Y}(T_i \, \cdot; 
f_i^{-1}(t_2z), \overline{f_i^{-1}}(s_2\bar{z}))T_{i+1} \, \cdot)
(u_1\otimes \dots \otimes u_{m+1}) \, \, t_1^{p_1}s_1^{q_1},
\eeq
is absolutely convergent, when $r> |t_2|,|s_2|>0$
and $1\geq |t_1|, |s_1|>0$, to an multivalued 
analytic functions of $t_1, s_1, t_1, s_2$, whose restriction
on $s_1=\bar{t}_1, s_2=\bar{t}_2$ is
\bea
&&\hspace{-1cm}T_0
m_{m+1}(T_1u_1, \dots, t_1^{L^L(0)}s_1^{L^R(0)}T_iu_i|_{s_1=\bar{t}_1},  
t_1^{L^L(0)}s_1^{L^R(0)}T_{i+1}u_{i+1}|_{s_1=\bar{t}_1}, 
\dots,  T_{m+1}u_{m+1}; \nn 
&&\hspace{-0.5cm} z_1, \bar{z}_1, \dots, 
z_i+ t_1f_i^{-1}(t_2z), \bar{z}_i + s_1
\overline{f_i^{-1}}(s_2\bar{z})|_{s_1=\bar{t}_1, s_2=\bar{t}_2}, 
z_i, \bar{z}_i, \dots, z_m, \bar{z}_m)
\eea
where $T_i, T_{i+1}$ are given by (\ref{T-ii+1}) and 
\bea
T_0 &=& e^{-L_-^L(A^{(0)})-L_-^R(\overline{A^{(0)}})},  \nn
T_j &=& e^{-L_+^L(A^{(j)})-L_+^R(\overline{A^{(j)}})}(a_0^{(j)})^{-L^L(0)}
(\overline{a_0^{(j)}})^{-L^R(0)},   \quad \quad \quad\quad j<i, \nn
T_j &=& e^{-L_+^L(A^{(j-1)})-L_+^R(\overline{A^{(j-1}})}
(a_0^{(j-1)})^{-L^L(0)}(\overline{a_0^{(j-1)}})^{-L^R(0)}, \quad j>i+1. 
\eea

By the 
general property of analytic functions, we obtain that
the multiple series (\ref{switch-series}) is absolutely 
convergent when $r> |t_2|,|s_2|>0$
and $1\geq |t_1|, |s_1|>0$. Hence the series (\ref{preswitch-series})
is also absolutely convergent when $r> |t_2|,|s_2|>0$
and $1\geq |t_1|, |s_1|>0$. As a special case ($t_1=s_1=1$), 
we obtain that the following series: 
\beq  \label{series-case-4-final}
\sum_{p_2, q_2} \nu_m(\tilde{P})   \, _{^i}*_{^0}
P_{p_2,q_2} \nu_2(\tilde{Q})  (u_1\otimes \dots u_{m+1}) 
 \, \,  t_2^{p_2}s_2^{q_2}
\eeq
is absolutely convergent, when $r> |t_2|,|s_2|>0$, to 
a multivalued analytic function, denoted as $G(t_2, s_2)$. 
In particular, when $s_2=\bar{t}_2$, (\ref{series-case-4-final})
converges to 
\bea  \label{m-m+1}
&&T_0m_{m+1}(T_1u_1, \dots, T_{m+1}u_{m+1}; z_1, \bar{z}_1, \dots, 
z_{i-1}, \bar{z}_{i-1}, \nn
&&\hspace{1cm} z_i+f_i^{-1}(t_2z), \bar{z}_i+ 
\overline{f_i^{-1}}(s_2\bar{z})|_{s_2=\bar{t}_2}, z_i, \bar{z}_i, 
\dots, z_m, \bar{z}_m).
\eea

Let $g: t_2 \mapsto z_i+f_i^{-1}(t_2z)$. 
By the definition of sewing operation, for fixed 
$z_1, \dots, z_m, z$, $P _{^i}\infty_{^0} Q$ exists 
is equivalent to the statement 
that $g$ is well-defined on the unit disk
$B(0;1)=\{ 1\geq |t_2|>0 \}$
and $z_j \notin g(B(0;1))$ (or equivalently
$\bar{z}_j \notin \bar{g}(B(0,1))$) ) for all $j\neq i$. 
Therefore, the analytic function $G(t_2, s_2)$ is free of 
singularities in $\{ (t_2, s_2)| 1\geq |t_2|, |s_2| >0 \}$. 
By the property of generalized power series, 
(\ref{series-case-4-final}) must be absolutely convergent
when $1\geq |t_2|,|s_2|>0$. In particular, when 
$t_2=s_2=1$, the series (\ref{series-case-4-final})
converges absolutely to 
\bea
&&T_0 m_{m+1}(T_1u_1, \dots, T_{m+1}u_{m+1}; z_1, \bar{z}_1, \dots,
z_{i-1}, \bar{z}_{i-1}, \nn
&&\hspace{2cm} z_i+f_i^{-1}(z), \bar{z}_i+ \overline{f_i^{-1}(z)}, 
z_i, \bar{z}_i, \dots, z_m, \bar{z}_m),
\eea
which is nothing but 
$\nu_{m+1}(\tilde{P} _{^i}\widetilde{\infty}_{^0} \tilde{Q})$. 
We have then finished the proof in this case. 

The proofs of general cases are essentially same as that of 
Proposition 5.4.1 in \cite{H4}. Only difference is that 
the analytic family with respect to the variables $t_1, t_2$ 
there is replaced by $t_1, s_1, t_2, s_2$ here, 
where $t_1, t_2$ are for chiral part and $s_1,s_2$ are for 
antichiral part just as we have done in the fourth case. 
\epf

\begin{thm}  \label{cffa-K-iso-thm}
The category of conformal full field algebra over $V^L\otimes V^R$
is isomorphic to the category of smooth 
$\tilde{K}^{c^L} \otimes \overline{\tilde{K}^{\overline{c^R}}}$-algebras
over $V^L\otimes V^R$. 
\end{thm}

\pf
The proof is again similar to that of the 
Theorem 5.4.5 in \cite{H4}. 
One can use the proof of Proposition \ref{K-hat-ffa-prop} 
and Theorem \ref{K-hat-ffa-iso-thm}
as a guidance. We only outline the proof here. 

Given a conformal full field algebras over $V^L\otimes V^R$,
denoted as $(F, m, \rho)$,
we have a smooth 
$\tilde{K}^{c^L} \otimes \overline{\tilde{K}^{\overline{c^R}}}$-algebra
over $V^L\otimes V^R$ given
by Proposition \ref{cffa-K-alg-prop}. Moreover, this
construction (\ref{nu-0-L-one}) and (\ref{nu-cffa-K-alg}) is 
functorial. Hence we obtain a functor from the category of 
conformal full field algebras over $V^L\otimes V^R$ 
to the category of smooth 
$\tilde{K}^{c^L}\otimes \overline{\tilde{K}^{\overline{c^R}}}$-algebras
over $V^L\otimes V^R$.

Conversely, the canonical section $\psi\otimes \bar{\psi}$ 
gives a natural embedding from $\widehat{K}$
to $\tilde{K}^{c^L}\otimes \overline{\tilde{K}^{\overline{c^R}}}$
as partial operads.  Hence any smooth 
$\tilde{K}^{c^L}\otimes \overline{\tilde{K}^{\overline{c^R}}}$-algebra 
over $V^L\otimes V^R$, 
denoted as $(F,\nu,\rho)$, automatically gives a 
smooth $\widehat{K}$-algebra $(F, \C \one, \nu)$, where
$\one$ is the vacuum state in $F$ defined by 
\beq
\one = \nu_0( \psi\otimes \bar{\psi} (\mathbf{0}) ). 
\eeq
By Proposition \ref{K-hat-ffa-prop}, 
there is a natural structure of grading-restricted
$\R\times \R$-graded full field algebra on $F$.

In this case, we also have 
two elements $\omega^L, \omega^R \in \overline{F}$ given by
\beq
\omega^L = - \frac{\partial}{\partial \epsilon} 
\nu( \psi\otimes \bar{\psi} (A(\epsilon; 2))) 
\lbar_{\epsilon=0}; \quad 
\omega^R = - \frac{\partial}{\partial \bar{\epsilon}} 
\nu( \psi\otimes \bar{\psi} (A(\epsilon; 2))) 
\lbar_{\epsilon=0}.
\eeq
We define, for $n>0$,  
\bea   \label{def-L-l-r-n}
L^L(-n) &:=& -\frac{\partial}{\partial A_n^{(0)}} 
\nu_1( \psi\otimes \bar{\psi} (A^{(0)}, (1, \mathbf{0}))) \nn
L^R(-n) &:=& -\frac{\partial}{\partial \overline{A_n^{(0)}}} 
\nu_1( \psi\otimes \bar{\psi} (A^{(0)}, (1, \mathbf{0}))) \nn
L^L(0) &:=& -\frac{\partial}{\partial a}  
\nu_1( \psi\otimes \bar{\psi} (\mathbf{0}, (a, \mathbf{0})))\nn
L^R(0) &:=& -\frac{\partial}{\partial \bar{a}} 
\nu_1( \psi\otimes \bar{\psi} (\mathbf{0}, (a, \mathbf{0}))) \nn
L^L(n) &:=& -\frac{\partial}{\partial A_n^{(1)}} 
\nu_1( \psi\otimes \bar{\psi} (\mathbf{0}, (1, A^{(1)}))) \nn
L^R(n) &:=& -\frac{\partial}{\partial \overline{A_n^{(1)}}} 
\nu_1( \psi\otimes \bar{\psi} (\mathbf{0}, (1, A^{(1)})))
\eea
Using the method in the proof of Proposition 5.4.4 in \cite{H4}
(or equivalently that in the proof of Proposition 2.9),
it is easily to show that $[L^L(m), L^R(n)]=0$ for $m,n\in \Z$
and the set $\{ L^L(n) \}_{n\in \Z}$ 
($\{ L^R(n) \}_{n\in \Z}$) generates
a Virasoro algebra of central charge $c^L$ ($c^R$). 
Moreover, following the proof of Proposition 5.4.4 in \cite{H4} 
and using (\ref{L-I-equ}), we can show that 
\bea
\mathbb{Y}(\omega^L; z, \bar{z}) &=& \sum_{n\in \Z} L^L(n) z^{-n-2}, \nn
\mathbb{Y}(\omega^R; z, \bar{z}) &=& \sum_{n\in \Z} L^R(n) 
\bar{z}^{-n-2},
\eea
where $\mathbb{Y}(\cdot; z, \bar{z})\cdot = 
\nu_2(\psi\otimes \bar{\psi} (P(z)) )$. 
The definition of $L^L(0)$ and $L^R(0)$ 
in (\ref{def-L-l-r-n}) exactly
coincide with $\mathbf{d}^L$ and $\mathbf{d}^R$ in (\ref{d-L-R-nu})
respectively.  Hence $L^L(0)$ and $L^R(0)$ are exactly 
the left and the right grading operators, respectively. 
Moreover, the definition of $L^L(-1), L^R(-1)$ coincides with 
$D^L, D^R$ in (\ref{D-L-R-const}) respectively. Therefore, 
$F$ has a structure of conformal full field algebra. 

Since these two Virasoro elements are completely determined by
a distinguished sphere with tube, they must coincide with
the Virasoro elements of $V^L$ and $V^R$ respectively.
Hence $F$ is a conformal full field algebra over $V^L\otimes V^R$.

Thus we obtain a functor from the category of 
smooth 
$\tilde{K}^{c^L}\otimes \overline{\tilde{K}^{\overline{c^R}}}$-algebras 
over $V^L\otimes V^R$
to the category of conformal full field algebras over $V^L\otimes V^R$.

Now we prove that these two functors give isomorphisms. 
Given a conformal full field algebra over $V^L\otimes V^R$, 
denoted as $(F, m, \rho)$. We obtain another such algebraic
structure $(F, \tilde{m}, \tilde{\rho})$ from 
$(F, \nu, \rho)$, which is a smooth 
$\tilde{K}^{c^L}\otimes \overline{\tilde{K}^
{\overline{c^R}}}$-algebras 
over $V^L\otimes V^R$ induced from $(F,m,\rho)$. 
$\tilde{m}=m$ is proved 
in (\ref{m-tilde-m}), and $\tilde{\rho}=\rho$ is automatic. 
Hence one way of composing these two functors gives the
identity functor. 

Given a smooth 
$\tilde{K}^{c^L}\otimes \overline{\tilde{K}^{\overline{c^R}}}$-algebras 
over $V^L\otimes V^R$, 
denoted as $(F,\nu, \rho)$. 
We have another such algebra $(F, \tilde{\nu}, \tilde{\rho})$ 
given by a conformal full field algebra over $V^L\otimes V^R$
$(F, m, \rho)$ which is further induced from 
$(F, \nu, \rho)$. 
Let $\tilde{Q} = \lambda \psi\otimes \bar{\psi}(Q) 
\in \tilde{K}^{c^L}\otimes \overline{\tilde{K}^{\overline{c^R}}}(n)$, 
where $Q$ is given in (\ref{ele-K-n}). We have
\bea
\tilde{\nu}_n
(\tilde{Q}): &=& \lambda e^{-L_-^L(A^{(0)}) - L_-^R(\overline{A^{(0)}})}  
m_{n}( e^{-L_+^L(A^{(1)}) - L_+^R( \overline{A^{(1)} } ) } 
(a_0^{(1)})^{-L^L(0)} \overline{a_0^{(1)}}^{-L^R(0)} u_1,  \dots, \nn
&&\hspace{0.2cm}  
e^{-L_+^L(A^{(n)}) - L_+^R( \overline{A^{(n)} })} 
(a_0^{(1)})^{-L^L(0)} \overline{a_0^{(1)}}^{-L^R(0)}u_n;  
z_1, \bar{z}_1, \dots, z_{n-1}, \bar{z}_{n-1}, 0,0),  \nn
&=& \nu_1(\tilde{Q}_0) _{^1}*_{^0} 
\big( ( \dots ( \nu_n(\tilde{P})_{^1}*_{^0}  
\nu_1(\tilde{Q}_1)) 
\dots  ) _{^n}*_{^0} \nu_1(\tilde{Q}_n) \big), 
\eea
where 
\bea
&\tilde{P} = \lambda \psi\otimes \bar{\psi}
\big( (z_1, \dots, z_{n-1}, 0; \mathbf{0}, (1, \mathbf{0}), 
\dots (1, \mathbf{0}))\big),&  \nn
&\tilde{Q}_0 = \psi\otimes \bar{\psi}(A^{(0)}, (1, \mathbf{0})), 
\quad 
\tilde{Q}_i = \psi\otimes \bar{\psi}(\mathbf{0}, (a_0^{(i)}, A^{(i)})), 
\quad i=1, \dots, n.&
\eea
Using the defining property of $\nu$, and the fact that
$$
\tilde{Q}= (\tilde{Q}_0) _{^1}\widetilde{\infty}_{^0} 
( \dots ( \tilde{P} _{^1}\widetilde{\infty}_{^0} Q_1 ) \dots \,
_{^n}\widetilde{\infty}_{^0} Q_n), 
$$
we obtain $\tilde{\nu}_n(\tilde{Q})=\nu_n(\tilde{Q})$ immediately. 
$\tilde{\rho}=\rho$ is obvious. 
Hence we have shown that 
the opposite way of composing these two functors 
also gives the identity functor. 

Therefore, the two categories are isomorphic. 
\epf

\renewcommand{\theequation}{\thesection.\arabic{equation}}
\renewcommand{\thethm}{\thesection.\arabic{thm}}
\setcounter{equation}{0}
\setcounter{thm}{0}

\section{Invariant bilinear forms}

In this section, we study the invariant bilinear form of 
conformal full field algebra, and give 
a geometric interpretation of a conformal full field
algebra over $V^L\otimes V^R$ 
equipped with a nondegenerate invariant bilinear form. 

An invariant bilinear form $(\cdot, \cdot)$ on 
a conformal full field algebra $F$ 
is a bilinear form on $F$ 
such that, for any $u, w_1, w_2\in F$, 
\bea \label{inv-form-ffa-1}
&&\hspace{-1cm}(w_2, \mathbb{Y}_f(u, x, \bar{x})w_1)  \nn
&&\hspace{-0.8cm}=(\mathbb{Y}_f(e^{-xL^L(1)}x^{-2L^L(0)}\otimes 
e^{-\bar{x}L^R(1)}\bar{x}^{-2L^R(0)} \, u, 
e^{\pi i} x^{-1}, e^{-\pi i} \bar{x}^{-1})w_2, w_1).
\eea 
or equivalently, 
\bea \label{inv-form-ffa-2}
&&\hspace{-1cm}(\mathbb{Y}_f(u, e^{\pi i}x, e^{-\pi i}\bar{x})w_2, w_1)  \nn
&&\hspace{-0.5cm}=(w_2, \mathbb{Y}_f(e^{xL^L(1)}x^{-2L^L(0)}\otimes 
e^{\bar{x}L^R(1)}\bar{x}^{-2L^R(0)} \, u, x^{-1}, \bar{x}^{-1})w_1). 
\eea

\begin{rema}  \label{rema-inv-biform}
{\rm The invariance property of bilinear form defined in 
(\ref{inv-form-ffa-1}) and (\ref{inv-form-ffa-2}) is
different from that in \cite{HK2}. Its geometric meaning is
shown in (\ref{geo-inv-form}). 
The main reason for 
the new definition is that the old definition can not be
formulate categorically. 
With this new definition, we can show in Section 4 that
a conformal full field algebra over $V^L\otimes V^R$,
where $V^L$ and $V^R$ satisfy certain natural conditions,  
has a very nice categorical formulation. 
}
\end{rema}

\begin{prop} 
For all $m,n\in \Z$ and $w_1, w_2\in F$, 
\beq \label{L-lr-form}
((L^L(m)\otimes L^R(n)) w_2, w_1) = 
(w_2, (-1)^{m+n} (L^L(-m)\otimes L^R(-n)) w_1). 
\eeq
\end{prop}
\pf
Replace $u$ in (\ref{inv-form-ffa-1}) by $\omega^L \otimes \one^R$ 
and use the fact that $L^L(1)\omega^L = 0$, we obtain
$$
(w_2, \mathbb{Y}_f(\omega^L\otimes \one^R, x, \bar{x}) w_1)
= (\mathbb{Y}_f(x^{-4} \omega^L\otimes \one^R, e^{\pi i} x^{-1}, 
e^{-\pi i} \bar{x}^{-1})w_2, w_1). 
$$
Expanding above equation by components, we obtain 
$$
\sum_{n\in \Z} (w_2, (L^L(n)\otimes 1) w_1) x^{-n-2} 
= \sum_{n\in \Z} ( (L^L(n)\otimes 1) w_2, w_1) x^{n-2} e^{\pi i (-n-2)}, 
$$
which implies that 
$$
( (L^L(n)\otimes 1)w_2, w_1) = (w_2, (-1)^{n} (L^L(-n)\otimes 1) w_1). 
$$
Similarly, we can show that 
$$
( (1\otimes L^R(n)) w_2, w_1) = (w_2, (-1)^{n} (1\otimes L^R(-n)) w_1).
$$
Then it is clear that (\ref{L-lr-form}) is true.  
\epf

We will use (\ref{L-lr-form})
in many places in this section without pointing it out explicitly.

\begin{prop} \label{sym-form-cl}
An invariant bilinear form $(\cdot, \cdot)$ on $F$ 
is automatically symmetric. Namely,  For $w_1, w_2\in F$, we have
\begin{equation} \label{symm-form}
(w_1, w_2) = (w_2, w_1)
\end{equation}
\end{prop}
\pf
The proof is similar to that of Proposition 5.3.6. in \cite{FHL}.
By (\ref{inv-form-ffa-1}), skew symmetry 
and (\ref{inv-form-ffa-2}), we have
\bea
&&\hspace{-0.5cm}(w_2, \mathbb{Y}_f(u, x, \bar{x})w_1)  \nn
&&\hspace{-0.4cm} = (\mathbb{Y}_f(e^{-xL^L(1)}x^{-2L^L(0)}
\otimes e^{-\bar{x}L^R(1)}
\bar{x}^{-2 L^R(0)}\, u, e^{\pi i}x^{-1}, e^{-\pi i}\bar{x}^{-1})w_2, w_1)  \nn
&&\hspace{-0.4cm}= (e^{-x^{-1}L^L(-1)}\otimes e^{-\bar{x}^{-1}L^R(-1)}
\mathbb{Y}_f(w_2, x^{-1}, \bar{x}^{-1})  \nn
&&\hspace{5cm} 
e^{-xL^L(1)}x^{-2L^L(0)}\otimes e^{-\bar{x}L^R(1)} \bar{x}^{-2 L^R(0)}\, u, w_1) \nn
&&\hspace{-0.4cm}= (e^{-xL^L(1)}x^{-2L^L(0)}\otimes e^{-\bar{x}L^R(1)}
\bar{x}^{-2 L^R(0)}\, u,     \nn
&&\hspace{2cm} \mathbb{Y}_f(e^{-x^{-1}L^L(1)} x^{2L^L(0)}\otimes 
e^{-\bar{x}^{-1}L^R(1)} \bar{x}^{2L^R(0)} w_2,  \nn
&&\hspace{5cm} e^{\pi i}x, e^{-\pi i}\bar{x})e^{x^{-1}L^L(1)}\otimes
e^{\bar{x}^{-1}L^R(1)}w_1)  \nn
&&\hspace{-0.4cm}= (x^{-2L^L(0)}\otimes \bar{x}^{-2 L^R(0)}\, u, 
\mathbb{Y}_f(e^{x^{-1}L^L(1)}\otimes e^{\bar{x}^{-1}L^R(1)}w_1, x,\bar{x})  \nn
&&\hspace{3cm} e^{-x^{-1}L^L(1)}x^{2L^L(0)}\otimes e^{-\bar{x}^{-1}L^R(1)}
\bar{x}^{2L^R(0)} w_2)   \nn
&&\hspace{-0.4cm}=(\mathbb{Y}_f(e^{-xL^L(1)}x^{-2L^L(0)}e^{x^{-1}L^L(1)}\otimes 
e^{-\bar{x}L^R(1)}\bar{x}^{-2L^R(0)}e^{\bar{x}^{-1}L^R(1)}w_1, 
e^{\pi i}x^{-1}, e^{-\pi i}\bar{x}^{-1})  \nn
&&\hspace{1.5cm} x^{-2L^L(0)}\otimes \bar{x}^{-2 L^R(0)}\, u,
e^{-x^{-1}L^L(1)}x^{2L^L(0)}\otimes e^{-\bar{x}^{-1}L^R(1)}
\bar{x}^{2L^R(0)}w_2)  \nn
&&\hspace{-0.4cm}=(\mathbb{Y}_f(x^{-2L^L(0)}\otimes \bar{x}^{-2
  L^R(0)}\, u, 
x^{-1}, \bar{x}^{-1}) x^{-2L^L(0)}\otimes \bar{x}^{-2L^R(0)}w_1,
x^{2L^L(0)}\otimes \bar{x}^{2L^R(0)} w_2)  \nn
&&\hspace{-0.4cm}=(\mathbb{Y}_f(u, x, \bar{x})w_1, w_2)
\eea
\epf 

Now we consider conformal full field algebra over $V^L\otimes V^R$. 
In this case, $F$ has a countable basis. We choose it to be
$\{ e_i \}_{i\in \N}$. 
If an invariant bilinear form $(\cdot, \cdot)$ on $F$
is also nondegenerate, we also have 
the dual basis $\{ e^i \}_{i\in \N}$. In this case, we can define a 
linear map $\Delta: \C \rightarrow \overline{F\otimes F}$ as follow:
\beq
\Delta:  1\mapsto \sum_{i\in \R} e_i\otimes e^i. 
\eeq

A conformal full field algebra over $V^L\otimes V^R$ 
with a nondegenerate invariant bilinear form 
has a very nice geometric interpretation. This is what we will
discuss in the remaining part of this section. 

Recall the notion of sphere with tubes of type $(n_-,n_+)$, 
where $n_-$ ($n_+$) is the number of 
negatively (positively) oriented tubes. 
In section 1, we have only studied the moduli space of  
spheres with tubes of type $(1,n)$, the structure on which
is captured in a notion 
called sphere partial operad. 
In this section, we would like to 
consider sphere with tubes of type $(n_-,n_+)$ for all 
$n_-,n_+\in \N$. 

In the case $n_-=1$, we have described the moduli space of conformal 
equivalent classes of sphere with tubes of type $(1, n)$ 
by choosing a canonical 
representative for each conformal equivalent class. 
In particular, we choose 
the only negatively oriented puncture to sit at 
$\infty \in \hat{\C}$. For general $n_-\in \N$,  there is no
canonical choice for the positions of the negatively oriented
punctures. We will take a different approach.  
We label the $i$-th negatively oriented puncture as $-i$-th puncture
and $j$-th positively oriented puncture as $j$-th puncture. 

Let us use 
\bea  \label{pre-moduli}
&&
Q=(\,  (z_{-n_-}; a_0^{(-n_-)}, A^{(-n_-)}), \dots,  
(z_{-1}; a_0^{(-1)}, A^{(-1)}) | 
\nn 
&&\hspace{5cm} (z_{1}; a_0^{(1)}, A^{(1)}), (z_{-n_-}; a_0^{(n_+)}, A^{(n_+)}) \, )
\eea
to denote a sphere $\hat{\C}$ 
with positively (negative) oriented punctures 
at $z_{i}\in \hat{\C}$ for  $i=1,\dots, n_+$ ($i=-1,\dots, -n_- $) ,
and with local coordinate map $f_{i}$ around each puncture $z_i$:  
\bea
f_{i}(w) &=& e^{\sum_j A_j^{(i)} x^{j+1}\frac{d}{dx}} (a_0^{(i)})^{x\frac{d}{dx}}x 
\lbar_{x=w-z_i} \quad \mbox{if $z_i\in \C$ }, \label{case-1} \\
&=& e^{\sum_j A_j^{(i)} x^{j+1}\frac{d}{dx}} (a_0^{(i)})^{x\frac{d}{dx}}x 
\lbar_{x=\frac{-1}{w}}  \quad\quad \mbox{if $z_i=\infty$}. \label{case-2}
\eea

\begin{rema}
{\rm 
Notice that (\ref{case-2}) is different from (\ref{f-0}), even when 
we set $a_0^{(i)}=1$, by a factor $-1$.  It simply means that 
the parametrizations of the corresponding tubes, 
given by (\ref{case-2}) and (\ref{f-0}), are different. 
So we expect that their algebraic realizations are also different, 
as one can see later by comparing (\ref{nu-cffa-K-alg}) with
(\ref{def-psi}), (\ref{def-psi-1}), (\ref{def-psi-2}) 
and keeping in mind of (\ref{L-lr-form})! Therefore,
nothing is really changed except one should 
distinguish the notation 
(\ref{pre-moduli}) with (\ref{ele-K-n}). 
}
\end{rema}

We denote the set of all such $Q$ as 
$\mathcal{T}(n_-, n_+)$, and their disjoint union as
$$
\mathcal{T} = \{ \mathcal{T}(n_-, n_+) \}_{n_-, n_+\in \N}.
$$
It is clear that there is an 
action of $SL(2, \C)$, the group of Mobius 
transformations, on $\mathcal{T}(n_-, n_+)$.
Each orbit of this action represents a single 
conformal equivalent class of sphere with tubes. 
We denote the set of orbits as 
$$
K(n_-, n_+) = \mathcal{T}(n_-, n_+)/ SL(2, \C).
$$
Let $K(0,0)$ be the one element set consisting of
the standard sphere $\hat{\C}$ 
with no additional structure.
Then we denote the disjoint union as: 
$$
\mathbb{K} = \{ K(n_-, n_+)  \}_{n_-,n_+\in \N}. 
$$

For a pair of elements $P\in K(m_-,m_+)$ 
and $Q\in K(n_-,n_+)$, a partial 
sewing operation $_{^i}\infty_{^{-j}}$ between the $i$-th positively 
oriented tube in $P$ 
and the $j$-th negatively oriented tube in $Q$ 
can be defined same as in Section 1. 
We denote the sphere with tube after sewing as 
$$
P _{^i}\infty_{^{-j}} Q
$$  
where $m_+ \geq i>0>-j\geq -n_-$. 
This sewing operation is associative. 
Notice that if we sew two spheres with
tubes along more than one pair of tubes with opposite orientations, we
obtain a surface of higher genus. We will not discuss this
type of multiple sewing operations between two spheres with tubes in this
work. By restricting to only single sewing operation between any pair of 
elements in $\mathbb{K}$, it is clear that $\mathbb{K}$ is
closed under these sewing operations. 
Moreover,  $\mathbb{K}$ is generated from the following elements:
\bea   \label{spe-ele}
& (\,(\infty; 1, \mathbf{0}) |\,) & \nn
&Q_2=(\, |(\infty; 1, \mathbf{0}), (0; 1, \mathbf{0})\, )&  \nn
&Q_{-2}= (\, (\infty; 1, \mathbf{0}), (0; 1,\mathbf{0})|\,) &  \nn
&P(z) = (\, (\infty; 1, \mathbf{0}) | (z; 1, \mathbf{0}), 
(0; 1, \mathbf{0})\, )&
\eea
and elements in $K(1,1)$ by sewing operations.

There is also a natural $S_{n_-} \times S_{n_+}$-action on $K(n_-,
n_+)$. Namely, $S_{n_+}$ ($S_{n_-}$) permute 
the positively (negatively) oriented punctures
among each other. We refer to the action of 
permutation group on $\mathbb{K}$ for all $n_-, n_+\in \N$
as $S_-\times S_+$-action.  

The structure on $\mathbb{K}$ induced from sewing operations and 
$S_- \times S_+$-actions are
richer than the structure of a partial operad. As we can see that it 
contains the sphere partial operad $K$ as a substructure.

The determinant line bundle on $\mathbb{K}$, denoted as 
$\text{Det}_{\mathbb{K}}$, is a trivial bundle on $\mathbb{K}$. 
The determinant line bundle on $K(0,0)$ is just a complex line. 
We fix an element of $\text{Det}_{Q_2}$ as $\mu_{2}(Q_{2})$ 
and an element of $\text{Det}_{Q_2}$ as $\mu_{-2}(Q_{-2})$. These choices
will determine a canonical section on the determinant line bundle 
$\text{Det}_{\mathbb{K}}$ over 
$\mathbb{K}$ (see the section 6.5 in \cite{H4} for detail).  
The sewing operations and $S_-\times S_+$ actions 
on $\mathbb{K}$ can be extended to those 
on the determinant line bundle. Moreover, these sewing operations
on $\text{Det}_{\mathbb{K}}$ is associative and compatible with 
$S_{n_-}\times S_{n_+}$-action on $\text{Det}_{\mathbb{K}}$
for all $n_+, n_-\in \N$.  We denote the sewing operation on 
$\text{Det}_{\mathbb{K}}$ as $\widetilde{\infty}$.

We denote the $c/2$ power of the determinant line bundle on 
$\mathbb{K}$ as $\tilde{\mathbb{K}}^c$. We are mainly interested in  
the line bundle 
$\tilde{\mathbb{K}}^{c^L} \otimes 
\overline{\tilde{\mathbb{K}}^{\overline{c^R}}}$.
We denote the canonical section in 
$\tilde{\mathbb{K}}^{c^L} \otimes \overline{\tilde{\mathbb{K}}^{\overline{c^R}}}$
as 
$$
\psi^L \otimes \psi^R: \mathbb{K} \rightarrow 
\tilde{\mathbb{K}}^{c^L} \otimes 
\overline{\tilde{\mathbb{K}}^{\overline{c^R}}}.
$$
and the sewing operation as $\widetilde{\infty}$.

We can extend the definition of smooth function on 
$\tilde{K}^{c^L} \otimes 
\overline{\tilde{K}^{\overline{c^R}}}$ to that of smooth function on 
$\tilde{\mathbb{K}}^{c^L} \otimes 
\overline{\tilde{\mathbb{K}}^{\overline{c^R}}}$ in the obvious way. 

When we want to emphasis the structures on 
$\tilde{\mathbb{K}}^{c^L} \otimes 
\overline{\tilde{\mathbb{K}}^{\overline{c^R}}}$,
which include partially defined sewing operations
$\widetilde{\infty}$ and $S_{n_-}\times S_{n_+}$-actions, 
we will denote them as a triple
$$
(\tilde{\mathbb{K}}^{c^L} \otimes 
\overline{\tilde{\mathbb{K}}^{\overline{c^R}}},
\widetilde{\infty}, S_-\times S_+). 
$$

We will be interested in some algebraic realization of 
$(\tilde{\mathbb{K}}^{c^L} 
\otimes \overline{\tilde{\mathbb{K}}^{\overline{c^R}}},
\widetilde{\infty}, S_-\times S_+)$.   Consider the set 
$$
\mathbb{F} = \{  
\hom (F^{\otimes m_+}, \overline{F^{\otimes m_-}} ) \}_{m_+, m_-\in \N}.
$$
where $F^{\otimes 0}=\C$. 
Let $f\in \hom (F^{\otimes m_+}, \overline{F^{\otimes m_-}})$ and
$g\in \hom (F^{\otimes n_+}, \overline{F^{\otimes n_-}})$.  
The map $g$ can be expanded as follow:
$$
g: v_1\otimes \dots \otimes v_{n_+} \mapsto \sum_{l_1,\dots l_{n_-}\in \N} 
g_{l_1\dots l_{n_-}}
(v_{i}, \dots, v_{n_+}) \, \, e_{l_1}\otimes \dots \otimes  e_{l_{n_-}}.
$$ 
We can define 
$f _{^i}*_{^{-j}} g$ for $1\leq i\leq m_+, 1\leq j\leq n_-$ as
\bea  \label{def-*}
&&\hspace{-1cm}
\sum_{k\in \R} \sum_{l_1,\dots l_{n_-}} 
g_{l_1\dots l_{n_-}}(u_{i}, \dots, u_{i+n_+}) e_{l_1}\otimes
\dots \otimes e_{l_{j-1}}   \nn
&&\hspace{-0.5cm} \otimes
f(u_1, \dots, u_{i-1}, e_{l_j},u_{i+n_++1}, \dots, u_{m_++n_+-1} ) \otimes 
e_{l_{j+1}} \otimes \dots \otimes e_{l_{n_-}} 
\eea
whenever the sum is absolutely convergent. 
$S_{m_-}\times S_{m_+}$-action on 
$\hom (F^{\otimes m_+}, \overline{F^{\otimes m_-}} )$ can 
be defined in the obvious way. 

Similar to algebras over partial operads, 
what we are looking for here is a smooth map 
$$\Psi: 
\tilde{\mathbb{K}}^{c^L} \otimes 
\overline{\tilde{\mathbb{K}}^{\overline{c^R}}}
\rightarrow \mathbb{F}
$$
such that $\Psi$ is equivariant with respective to 
$S_-\otimes S_+$-actions and satisfies the following conditions:
\beq  \label{morph-cond-bi-ffa} 
\Psi (\tilde{P} _{^i} \widetilde{\infty}_{^{-j}} \tilde{Q} ) = 
\Psi (\tilde{P}) _{^i}*_{^{-j}} \Psi(\tilde{Q})  
\eeq
for any $\tilde{P}, \tilde{Q}\in \tilde{\mathbb{K}}^{c^R} \otimes \overline{\tilde{\mathbb{K}}^{\overline{c^R}}}$ whenever 
$\tilde{P} _{^i} \widetilde{\infty}_{^{-j}} \tilde{Q}$ exists. 
We will show that a conformal full field algebra $F$ together with
a nondegenerate invariant bilinear form $(\cdot, \cdot)$ canonically
give a map $\Psi$ which is $S_-\times S_+$-equivariant and satisfies
(\ref{morph-cond-bi-ffa}).

Let $(F,m,\rho)$ 
be a conformal full field algebra over $V^L\otimes V^R$
equipped with a
nondegenerate bilinear form $(\cdot, \cdot)$. 
Let $\{ e^i \}_{i\in \N}$
be the dual basis of $\{ e_i \}_{i\in \N}$.

For $Q\in \mathcal{T}(n_-, n_+)$ of the form (\ref{pre-moduli}) 
and $\lambda\in \C^{\times}$, we define
\beq  \label{def-psi-0}
\Psi ( \lambda (\psi^L\otimes \psi^R) (Q) )
( u_{-n_-}\otimes \dots \otimes u_{n_+} )
\eeq 
in the following three cases: 
\bnu
\item If $z_k\neq \infty$ for all $k=-n_-, \dots, -1, 1, \dots, n_+$,  
\bea  \label{def-psi}
&&\hspace{-0.5cm}\lambda \sum_{i_1, \dots, i_{n_-} \in \N} \big( 
\one, \, m_{n_-+n_+}(
e^{-L_+^L(A^{(-n_-)})-L_+^R(\bar{A}^{(-n_-)})}
(a_0^{(-n_-)})^{-L^L(0)} \overline{a_0^{(-n_-)}}^{-L^R(0)}
e_{i_{n_-}}, \nn  
&&\hspace{3cm} \dots,  
e^{-L_+^L(A^{(-1)})-L_+^R(\bar{A}^{(-1)})} 
(a_0^{(-1)})^{-L^L(0)} \overline{a_0^{(-1)}}^{-L^R(0)}e_{i_1}, \nn  
&&\hspace{3.7cm} 
e^{-L_+^L(A^{(1)})-L_+^R(\bar{A}^{(1)})} 
(a_0^{(1)})^{-L^L(0)} \overline{a_0^{(1)}}^{-L^R(0)}u_1,  \nn
&&\hspace{3cm} \dots,
e^{-L_+^L(A^{(n_+)})-L_+^R(\bar{A}^{(n_+)})}
(a_0^{(n_+)})^{-L^L(0)}\overline{a_0^{(n_+)}}^{-L^R(0)}u_{n_+}; \nn
&&\hspace{1cm}  z_{-n_-}, \bar{z}_{-n_-}, \dots, 
z_{-1}, \bar{z}_{-1},
z_1, \bar{z}_1, \dots, z_{n_+}, \bar{z}_{n_+}) \big)\,  
e^{i_{1}}\otimes \dots \otimes e^{i_{n_-}};
\eea

\item If $\exists \, k \in \{ -n_-, \dots, -1\}$ such that 
$z_k=\infty$ (recall (\ref{case-2})), 
we define (\ref{def-psi-0}) to be the formula obtained 
from (\ref{def-psi}) by switching the first $\one$ with
\beq  \label{def-psi-1}
e^{-L_+^L(A^{(k)})-L_+^R(\bar{A}^{(k)})}(a_0^{(k)})^{-L^L(0)} 
\overline{a_0^{(k)}}^{-L^R(0)} e_{i_{-k}};
\eeq

\item If $\exists \, k \in \{1, \dots, n_+\}$ such that 
$z_k=\infty$ (recall (\ref{case-2})), 
we define (\ref{def-psi-0}) to be the formula 
obtained from (\ref{def-psi}) by switching the first $\one$ with the 
\beq  \label{def-psi-2}
e^{-L_+^L(A^{(k)})-L_+^R(\bar{A}^{(k)})}(a_0^{(k)})^{-L^L(0)} 
\overline{a_0^{(k)}}^{-L^R(0)} u_k.
\eeq

\enu

We have finished the definition of $\Psi$ in all cases. 
Some interesting cases are listed explicitly below: 
\bea
\Psi( \psi^L \otimes \psi^R ( \hat{\C} )) &=& 
(\one, \one) I_{\C},  \nn
\Psi( \psi^L \otimes \psi^R  (\, (\infty, 1, \mathbf{0})|\, ) ) 
&=& \one \nn
\Psi (\psi^L\otimes \psi^R( Q_2 ) )(u\otimes v) &=& (u, v) 
\label{geo-inv-form} \\
\Psi (\psi^L\otimes \psi^R( Q_{-2} )) &=& \Delta.  \nn
\Psi(\psi^L \otimes \psi^R (P(z)))(u\otimes v) &=& 
\mathbb{Y}(u, z, \bar{z})v  \nonumber
\eea
where $\hat{\C}$ is the single element in $K(0,0)$. 
We can always choose $(\cdot, \cdot)$ to be
so that $(\one, \one)=1$.

First question we want to ask is whether such defined $\Psi$ induces a 
well-defined map, still denoted as $\Psi$,  
on $\tilde{\mathbb{K}}^{c^L} \otimes 
\overline{\tilde{\mathbb{K}}^{\overline{c^R}}}$.

\begin{lemma} \label{well-def-psi}
$\Psi$ is well-defined on 
$\tilde{\mathbb{K}}^{c^L} \otimes 
\overline{\tilde{\mathbb{K}}^{\overline{c^R}}}$.
\end{lemma}
\pf
Firstly, we need to show that the coefficient of each
$e^{i_{1}}\otimes \dots \otimes e^{i_{n_-}}$ in (\ref{def-psi}) 
is invariant under the action of $SL(2, \C)$ on $Q$. 
It is equivalent to prove Lemma for  
$Q\in \mathcal{T}(0, n_-+n_+)$. 
Hence we only need to consider the case when $n_-=0$. 

Secondly, if $z_k\neq \infty,0$ for all $k=-n_-, \dots, n_+$,  then 
(\ref{def-psi}) equals to  
$$
\Psi ( \lambda (\psi^L\otimes \psi^R) (P) )
( u_{-n_-}\otimes \dots \otimes u_{n_+} \otimes \one \otimes \one),
$$
where $P$ is obtained from $Q$ by adding 
the $n_++1$-th and $n_++2$-th positively
oriented tubes at $\infty$ and $0$ respectively
with arbitrary local coordinates maps. 
Hence we conclude that it is sufficient to 
consider the case $z_k=\infty$ and $z_l=0$ 
for some $-n_-\leq k,l \leq n_+$. 
For simplicity, we can assume $k=1$ and $l=n_+$.  
The proof for general $k,l$ is exactly the same. 

The group $SL(2, \C)$ 
is generated by the following three Mobius transformations:
\bea
w &\mapsto& w-a  \nn
w &\mapsto& aw   \nn
w &\mapsto& -\frac{1}{w}
\eea 
where $a\in \C^{\times}$. Hence it is enough to show the 
invariance of $\Psi$ with respect to above
three Mobius transformations.

We first consider the transformation $w\mapsto w-a$. It maps a 
$Q\in \mathcal{T}(0, n_+)$ of form:
$$
(\, | (\infty, a_0^{(1)}, A^{(1)}), (z_2, a_0^{(2)}, A^{(2)}), \dots, 
(z_{n_+}, a_0^{(n_+)}, A^{(n_+)}) \, ) 
$$
($z_{n_+}=0$) to 
$$
Q'=(\, | (\infty, b_0, B), (z_2-a, a_0^{(2)}, A^{(2)}), \dots, 
(z_{n_+}-a, a_0^{(n_+)}, A^{(n_+)}) \, ),
$$
where $b_0, B$ is so that the local coordinate map at $\infty$ 
is given by $f_1(w+a)$, where 
$f_1(w)$ is given by (\ref{case-2}). We have
\bea
f_1(w+a) &=& e^{a\frac{d}{dw}} f_1(w) \nn
&=& e^{ax^2\frac{d}{dx}} 
e^{\sum_j A_j^{(1)} x^{j+1}\frac{d}{dx}} (a_0^{(1)})^{x\frac{d}{dx}}x 
\lbar_{x=\frac{-1}{w}}
\eea
By the definition of $\Psi$ given in (\ref{def-psi}), we see 
that
\bea
&&\hspace{-1cm} 
\Psi \circ \psi^L\otimes \psi^R (Q') = \big( e^{-aL^L(1)-\bar{a}L^R(1)}u_1', \nn
&&\hspace{2cm}m_{n_+-1}(u_2', \dots, u_{n_+}'; 
z_2-a, \overline{z_2-a}, \dots, z_{n_+}-a, \overline{z_{n_+}-a} )\, \big)
\label{welldef-proof-1}
\eea
where 
\beq
u_j'= e^{-L_+^L(A^{(k)})-L_+^R(\bar{A}^{(k)})}(a_0^{(k)})^{-L^L(0)} 
\overline{a_0^{(k)}}^{-L^R(0)}u_j
\eeq
for $j=1, \dots, n_+$. 
Apply (\ref{L-lr-form}) and (\ref{D-m-n}) to (\ref{welldef-proof-1}), 
we obtain that 
$$
\Psi \circ \psi^L\otimes \psi^R (Q') =
\Psi \circ \psi^L\otimes \psi^R (Q).
$$

The transformation $w\mapsto aw$ maps $Q$ in the form 
(\ref{pre-moduli}) to $Q''$, which has only positively oriented 
punctures at $\infty, az_2, \dots, az_{n_+}$ with
local coordinate maps given by $f_i(a^{-1}w)$. We have 
$$
f_1(a^{-1}w) = a^{-w\frac{d}{dw}}f_1(w) = a^{x\frac{d}{dx}} 
e^{\sum_j A_j^{(1)} x^{j+1}\frac{d}{dx}} (a_0^{(1)})^{x\frac{d}{dx}}x 
\lbar_{x=\frac{-1}{w}}
$$
and $f_i(a^{-1}w)=a^{-w\frac{d}{dw}}f_i(w)$ for $i=2, \dots, n_+$. 
By the definition (\ref{def-psi}), we obtain
\bea
&&\Psi \circ \psi^L\otimes \psi^R (Q'') \nn
&&\hspace{1cm}= \big( a^{-L^L(0)}\bar{a}^{-L^R(0)}u_1', 
m_{n_+-1}(a^{L^L(0)}\bar{a}^{L^R(0)}u_2',  \dots 
a^{L^L(0)}\bar{a}^{L^R(0)}u_{n_+}';  \nn
&&\hspace{5cm}
az_2, \overline{az_2}, \dots, az_{n_+}, \overline{az_{n_+}})\, \big).
\eea  
Using (\ref{d-l-d-r-conj}), it is clear that 
$$
\Psi \circ \psi^L\otimes \psi^R (Q'') =
\Psi \circ \psi^L\otimes \psi^R (Q).
$$

Now we consider the last Mobius transformation 
$w\rightarrow -\frac{1}{w}$. It transforms $Q$ to
$Q'''$, which has punctures at $0, \frac{-1}{z_2}, \dots, 
\frac{-1}{z_{n_+-1}}, \infty$ with local coordinate maps 
$f_i(\frac{-1}{w})$. In particular, for $i=1$, we have
$$
f_1\left( -\frac{1}{w}\right) = 
e^{\sum_{j=1}^{\infty} A_j^{(1)} w^{j+1}\frac{d}{dw}} w; 
$$
and for $i=2, \dots, n_+-1$, 
\bea
f_i \left( \frac{-1}{w} \right) &=& e^{\sum_j A_j^{(i)}x^{j+1}\frac{d}{dx}} 
(a_0^{(i)})^{x\frac{d}{dx}} x \lbar_{x=\frac{-1}{w}-z_i}  \nn
&=& e^{z_ix^2\frac{d}{dx}} (z_i)^{2x\frac{d}{dx}} e^{\sum_j A_j^{(i)}x^{j+1}\frac{d}{dx}} 
(a_0^{(i)})^{x\frac{d}{dx}} x \lbar_{x=w+1/z_i};  \nonumber
\eea
and for $i=n_+$ , 
$$
f_{n_+}\left( -\frac{1}{w} \right) = e^{\sum_j A_j^{(n_+)}x^{j+1}\frac{d}{dx}} 
(a_0^{(n_+)})^{x\frac{d}{dx}} x \lbar_{x=-\frac{1}{w}}. 
$$
By the definition of $\Psi$, we obtain
\bea
&&\Psi \circ \psi^L\otimes \psi^R (Q''') \nn
&&\hspace{0.5cm}= \big( u_{n_+}', 
m_{n_+-1}(e^{-z_2L^L(1)-\bar{z}_2L^R(1)}z_2^{-2L^L(0)}\bar{z}_2^{-2L^R(0)}u_2', \nn
&&\hspace{2.5cm} e^{-z_{n_+-1}L^L(1)-\bar{z}_{n_+-1}L^R(1)}
z_{n_+-1}^{-2L^L(0)}\bar{z}_{n_+-1}^{-2L^R(0)}u_{n_+-1}', u_1';  \nn
&&\hspace{2cm}
-1/z_2, \overline{-1/z_2}, \dots,
-1/z_{n_+-1}, \overline{-1/z_{n_+-1}},0,0) \big).
\eea

In the cases that $z_2, \dots, z_{n_+-1}$ have 
pairwise distinct absolute values, we can expand $m_{n_+-1}$
as a $n_+-1$-products of $\mathbb{Y}$s in a proper order. 
Then using the fact that the bilinear form
$(\cdot, \cdot)$ is symmetric and invariant, 
it is easy to see that 
\beq \label{case-3-equ}
\Psi \circ \psi^L\otimes \psi^R (Q''') =
\Psi \circ \psi^L\otimes \psi^R (Q)
\eeq
holds in these cases. 
The remaining cases are the complement of an open dense subset,
all the cases in which have been proved.  
By the smoothness of $m_{n_+-1}$ with respect 
to $z_2, \dots, z_{n_+-1}$, 
(\ref{case-3-equ}) holds in all remaining cases.
\epf

The following Lemma follows immediately from the definition of $\Psi$.
\begin{lemma}
Let $Q$ be an element of $K(n_-, n_+)$ of the form (\ref{pre-moduli}). 
$Q_{\pm 2}$ and $Q$ are always sewable with $Q$ along oppositely oriented
tubes. Denote $\psi^L\otimes \psi^R(Q_{\pm 2})$ 
and $\psi^L\otimes \psi^R(Q)$ as $\tilde{Q}_{\pm 2}$ 
and $\tilde{Q}$ respectively. 
Then for all $1\leq k\leq n_-$ and $1\leq l\leq n_+$, we have 
\bea \label{lemma-2-1}
\Psi(\tilde{Q}_{2}) _{^i}\widetilde{\infty}_{^{-k}} \tilde{Q}  
&=& \Psi( \tilde{Q}_{2} )_{^i}*_{^{-k}} \Psi( \tilde{Q} ), 
\quad \forall i=1,2,\nn
\Psi( \tilde{Q} _{^l}\widetilde{\infty}_{^{i}} \tilde{Q}_{-2} )
&=& \Psi( \tilde{Q} ) _{^l}*_{^{-i}} 
\Psi(\tilde{Q}_{-2} ),  \quad \forall i=1,2.
\eea
\end{lemma}

\begin{thm}  \label{V-cl-hom}
$\Psi$ is smooth $S_-\times S_+$-equivariant and 
satisfies (\ref{morph-cond-bi-ffa}).  
\end{thm}
\pf
The smoothness of $\Psi$ follows from that of $m_n$ and our 
construction (\ref{def-psi}). 
It is also clear that $\Psi$ is $S_-\times S_+$-equivariant 
because of the permutation property of $m_n$. 

It remains to show (\ref{morph-cond-bi-ffa}). 
Let $P\in K(m_-, m_+)$ and $Q\in K(n_-, n_+)$ have the following form: 
\bea
P &=& ( (\xi_{-m_-}, a_0^{(-m_-)}, A^{(-m_-)}), \dots ,
(\xi_{-1}, a_0^{(-1)}, A^{(-1)}) | \nn
&&\hspace{3cm} (\xi_1, a_0^{(1)}, A^{(1)}), 
\dots, (\xi_{m_+}, a_0^{(m_+)}, A^{(m_+)}) )  ,\nn
Q &=& ( (\eta_{-n_-}, a_0^{(-n_-)}, A^{(-n_-)}), \dots ,
(\eta_{-1}, a_0^{(-1)}, A^{(-1)}) | \nn
&&\hspace{3cm} (\eta_1, a_0^{(1)}, A^{(1)}), 
\dots ,(\eta_{n_+}, a_0^{(n_+)}, A^{(n_+)}) ).
\eea
Let $\tilde{P}=\psi^L\otimes \psi^R(P)$ and $\tilde{Q}=\psi^L\otimes \psi^R(Q)$. 
Assume that $\tilde{P} _{^i}\widetilde{\infty}_{^{-j}} \tilde{Q}$ exists 
for $1\leq i\leq m_+$ and $1\leq j\leq n_-$.  
Because $\Psi$ is $S_-\times S_+$-equivariant, we can further assume that 
$i=m_+$ and $j=n_-$.  By Lemma \ref{well-def-psi},
we can choose $\xi_{-m_-}=\infty, \xi_{m_+}=0$ and 
$\eta_{-n_-}=\infty, \eta_{n_+}=0$. 

Let  $P'\in K(1,m_-+m_+-1)$ and $Q'\in K(1, n_-+n_+-1)$ be given as
\bea
&&P'=( (\infty, a_0^{(-m_-)}, A^{(-m_-)})| 
(\xi_{-m_-+1}, a_0^{(-m_-+1)}, A^{(-m_-+1)}),
\dots, (0, a_0^{(m_+)}, A^{(m_+)}), \nn
&&Q'=(\, (\infty, a_{0}^{(-n_-)}, A^{(-n_-)}) | 
(\eta_{-n_-+1}, a_0^{(-n_-+1)}, A^{(-n_-+1)}), \dots,  
(0, a_{0}^{(n_+)}, A^{(n_+)})\, ). \nonumber
\eea
Then it is clear that
\bea
P&=& (\dots (P' _{^1}\infty_{^{-1}} Q_{-2}) _{^1}\infty_{^{-1}} 
\dots \,) _{^{1}}\infty_{^{-1}} Q_{-2}, \nn
Q &=& (\dots (Q' _{^1}\infty_{^{-1}} Q_{-2}) _{^1}\infty_{^{-1}} 
\dots \,) _{^{1}}\infty_{^{-1}} Q_{-2}. 
\eea
Let $\tilde{P}'$ and $\tilde{Q}'$ 
be elements in the fiber over $P$ and $Q$ 
respectively such that
\bea
\tilde{P}&=& (\dots (\tilde{P}' _{^1}\widetilde{\infty}_{^{-1}} 
\tilde{Q}_{-2}) _{^1}\widetilde{\infty}_{^{-1}} 
\dots \,) _{^{1}}\widetilde{\infty}_{^{-1}} \tilde{Q}_{-2}, \nn
\tilde{Q}&=& (\dots (\tilde{Q}' _{^1}\widetilde{\infty}_{^{-1}} 
\tilde{Q}_{-2}) _{^1}\widetilde{\infty}_{^{-1}} 
\dots \,) _{^{1}}\widetilde{\infty}_{^{-1}} \tilde{Q}_{-2}. 
\eea
By (\ref{lemma-2-1}), we obtain
\bea  \label{P-P'-Q-Q'}
\Psi(\tilde{P})&=& 
(\dots (\Psi(\tilde{P}') _{^1}*_{^{-1}} \Psi(\tilde{Q}_{-2}))
_{^1}*_{^{-1}}  \dots \,) _{^{1}}*_{^{-1}} \Psi(\tilde{Q}_{-2}),\nn
\Psi(\tilde{Q})&=& 
(\dots (\Psi(\tilde{Q}') _{^1}*_{^{-1}} \Psi(\tilde{Q}_{-2}))
_{^1}*_{^{-1}}  \dots \,) _{^{1}}*_{^{-1}} \Psi(\tilde{Q}_{-2}).
\eea
Then we  have 
\bea
\Psi(\tilde{P}) _{^i}*_{^{-j}} \Psi(\tilde{Q} ) &=& 
(\, (\dots (\Psi(\tilde{P}') _{^1}*_{^{-1}} \Psi(\tilde{Q}_{-2}))
_{^1}*_{^{-1}}  \dots \,) _{^{1}}*_{^{-1}} \Psi(\tilde{Q}_{-2})\, )
\nn
&&\hspace{0.3cm} _{^i}*_{^{-j}} 
(\, (\dots (\Psi(\tilde{Q}') _{^1}*_{^{-1}} \Psi(\tilde{Q}_{-2}))
_{^1}*_{^{-1}}  \dots) _{^{1}}*_{^{-1}} \tilde{Q}_{-2})\, ). \nonumber
\eea
One can check easily that the right hand side of 
above equation can be obtained equivalently by first 
doing the contraction 
$\Psi(\tilde{P'}) _{^{m_-+m_+-1}}*_{^{-1}} \Psi(\tilde{Q}')$ then 
doing the remaining contractions with $\Psi(\tilde{Q}_{-2})$. 
By Theorem \ref{cffa-K-iso-thm}, we have
\beq  \label{H-thm-equ}
\Psi(\tilde{P'}) _{^{m_-+m_+-1}}*_{^{-1}} \Psi(\tilde{Q}')
=  \Psi( \tilde{P'} _{^{m_-+m_+-1}}\widetilde{\infty}_{^{-1}} \tilde{Q}').
\eeq
By (\ref{lemma-2-1}) again, the remaining contractions between
both sides of (\ref{H-thm-equ}) and $\Psi(\tilde{Q}_{-2})$ 
give exactly $\Psi(\tilde{P}) _{^{i}}*_{^{-j}} \Psi(\tilde{Q})$
and $\Psi(\tilde{P} \, _{^{i}}\widetilde{\infty}_{^{-j}} \tilde{Q})$
respectively. 
Therefore, we obtain the identity (\ref{morph-cond-bi-ffa}). 
\epf

\renewcommand{\theequation}{\thesection.\arabic{equation}}
\renewcommand{\thethm}{\thesection.\arabic{thm}}
\setcounter{equation}{0}
\setcounter{thm}{0}

\section{Commutative associative algebras in $\mathcal{C}_{V^L\otimes V^R}$}

In this section, we study the categorical 
formulation of conformal full field algebra over $V^L\otimes V^R$.

Let $V$ be a vertex operator algebra,  which satisfies the 
conditions in Theorem \ref{ioa} by our assumption. 
Then $\mathcal{C}_V$, the category of $V$-module, has
a structure of vertex tensor category. 
We review some of the ingredients of vertex tensor 
category $\mathcal{C}_V$, 
and set our notations along the way.

There is a tensor product bifunctor 
$\boxtimes_{P(z)}: \mathcal{C}_V \times \mathcal{C}_V \rightarrow 
\mathcal{C}_V$ 
for each $P(z), z\in \C^{\times}$ in the sphere partial operad $K$. 
We also denote $\boxtimes_{P(1)}$ simply as $\boxtimes$.

Let $W_1$ and $W_2$ be $V$-modules. For a given path 
$\gamma \in \C^{\times}$ from a point $z_1$ to $z_2$, there is
a parallel isomorphism associated to this path 
$$
\mathcal{T}_{\gamma} : W_1 \boxtimes_{P(z_1)} W_2 \longrightarrow 
W_1 \boxtimes_{P(z_2)} W_2.  
$$
Let $\Y$ be the intertwining operator corresponding 
to the intertwining 
map $\boxtimes_{P(z_2)}$ and $l(z_1)$ the value of the logarithm of 
$z_1$ determined by $\log z_2$ and analytic continuation 
along the path $\gamma$. 
For $w_1\in W_1,w_2\in W_2$, the map $\overline{\mathcal{T}_{\gamma}}: 
\overline{W_1\boxtimes_{P(z_1)} W_2} \rightarrow 
\overline{W_1\boxtimes_{P(z_1)} W_2}$ as an natural extension of 
$\mathcal{T}_{\gamma}$, is defined by
$$
\overline{\mathcal{T}_{\gamma}} (w_1\boxtimes_{P(z_1)} w_2)
=\Y(w_1, e^{l(z_1)})w_2,
$$
which uniquely determines
$\mathcal{T}_{\gamma}$ since the homogeneous components of 
$w_1\boxtimes_{P(z_1)} w_2\in 
\overline{W_1\boxtimes_{P(z_1)} W_2}$ span the module
$W_1\boxtimes_{P(z_2)}W_2$. Moreover, the parallel isomorphism
depends only on the homotopy class of $\gamma$. 

A remark of our notation, we use ``overline'' for both
complex conjugation and natural extension of morphisms in 
tensor categories. There should be no confusion because they
act on different things.

For each $V$-module $W$, there is a left unit isomorphism 
$l_W: V\boxtimes W \rightarrow W$ defined by
\beq  \label{l-unit-cat}
\overline{l_W}(v\boxtimes w) 
= Y_{W}(v,1)w, \quad \quad \forall v\in V, w\in W,
\eeq
where $Y_W$ is the vertex operator which defines the module
structure on $W$, and a right unit isomorphism 
$r_W: W\boxtimes V \rightarrow W$ defined by
\beq \label{r-unit-cat}
\overline{r_W}(w\boxtimes v) = e^{L(-1)}Y_{W}(v, -1)w, 
\quad \quad \forall v\in V, w\in W.
\eeq

There is an associativity isomorphism, for each triple of 
$V$-modules $W_1, W_2, W_3$,  
$$
\mathcal{A}^{P(z_{1}-z_{2}), P(z_{2})}_{P(z_{1}), P(z_{2})}: 
W_1\boxtimes_{P(z_{1})} (W_2\boxtimes_{P(z_{2})} W_3)\to 
(W_1\boxtimes_{P(z_{1}-z_{2})} W_2)\boxtimes_{P(z_{2})} W_3,
$$
which is characterized by 
\begin{equation}\label{assoc-iso}
\overline{\mathcal{A}^{P(z_{1}-z_{2}), P(z_{2})}_{P(z_{1}), P(z_{2})}}
(w_{(1)}\boxtimes_{P(z_{1})} (w_{(2)}\boxtimes_{P(z_{2})} w_{(3)}))
=(w_{(1)}\boxtimes_{P(z_{1}-z_{2})} w_{(2)})\boxtimes_{P(z_{2})}w_{(3)}
\end{equation}
for $w_{(i)} \in W_i, i=1,2,3$. 
Let $z_1>z_2>z_1-z_2>0$. The associativity isomorphism $\mathcal{A}$ 
of the braided tensor category is
\begin{equation}
\mathcal{A}=\mathcal{T}_{\gamma_{3}}\circ (\mathcal{T}_{\gamma_{4}}
\boxtimes_{P(z_{2})} I)\circ 
\mathcal{A}^{P(z_{1}-z_{2}), P(z_{2})}_{P(z_{1}), P(z_{2})}\circ
(I \boxtimes_{P(z_{1})} 
\mathcal{T}_{\gamma_{2}})\circ \mathcal{T}_{\gamma_{1}}.\label{A-assoc}
\end{equation}
where $\gamma_1$ and $\gamma_2$ are paths in $\R_+$
from $1$ to $z_1$ and $z_2$, respectively; and $\gamma_3$ and
$\gamma_4$ are paths in $\R_+$ from 
$z_2$ and $z_1-z_2$ to $1$, respectively.

There is also a braiding isomorphism 
$\mathcal{R}^{+}_{W_1W_2}: W_1\boxtimes W_2 \rightarrow W_2\boxtimes W_1$ 
for each pair of $V$-modules $W_1,W_2$, defined as 
\begin{equation}  \label{R-+-cat}
\overline{\mathcal{R}_{W_1W_2}^{+}}(w_1\boxtimes w_2) = e^{L(-1)} 
\overline{\mathcal{T}}_{\gamma_+} (w_2\boxtimes_{P(-1)} w_1),
\end{equation}
where $\gamma_+$ is a path from $-1$ to $1$ inside upper half
plane as shown in the following diagram.
\beq    \label{gamma_+-fig}
\epsfxsize  0.4\textwidth
\epsfysize  0.2\textwidth
\epsfbox{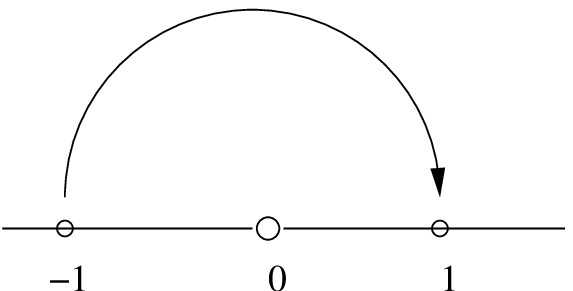}
\eeq
The inverse of $\mathcal{R}_{W_1W_2}^{+}$ is denoted as 
$\mathcal{R}_{W_1W_2}^-$, which is characterized by 
\begin{equation}   \label{R---cat}
\overline{\mathcal{R}_{W_1W_2}^{\, -}}(w_1\boxtimes w_2) = e^{L(-1)} 
\overline{\mathcal{T}}_{\gamma_-} (w_2\boxtimes_{P(-1)} w_1),
\end{equation}
where $\gamma_-$ is a path in the lower half plane as shown in the
following picture.
\beq        \label{gamma_--fig}
\epsfxsize  0.4\textwidth
\epsfysize  0.2\textwidth
\epsfbox{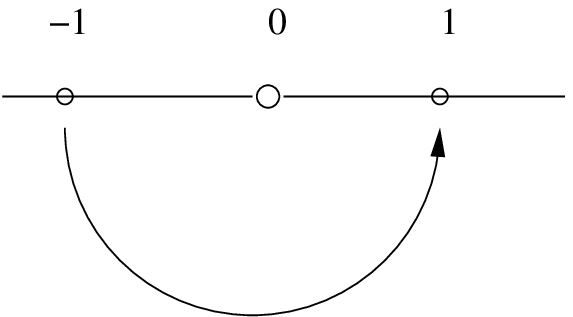}
\eeq 
For simplicity, we often write $\mathcal{R}_{W_1W_2}^{+}$ 
and $\mathcal{R}_{W_1W_2}^-$ as $\mathcal{R}_{+}$ 
and $\mathcal{R}_{-}$ respectively. Notice that
we change the superscripts to subscripts because 
we need save space for some superscripts needed later.

Let $W_i, i=1,2,3$ be $V$-modules and 
$\V_{W_1W_2}^{W_3}, \V_{W_2W_1}^{W_3}$ be the space of 
intertwining operators of type $\binom{W_3}{W_1W_2}$ 
and $\binom{W_3}{W_2W_1}$.
We have isomorphisms $\Omega_r: \V_{W_1W_2}^{W_3} \rightarrow
\V_{W_2W_1}^{W_3}, r\in \Z$ for $V$-modules $W_i, i=1,2,3$, 
given as follow:
\beq
\Omega_r(\Y)(w_2, z)w_1 = e^{zL(-1)} \Y(w_1, e^{(2r+1)\pi i}z)w_2,
\eeq
for $\Y\in \V_{W_1W_2}^{W_3}$. 

For any $\Y \in \V_{W_1W_2}^{W_3}$ and $z\in \C^{\times}$,  
by the universal property of tensor 
product\cite{HL1}-\cite{HL4}, 
there is a unique module map $W_1\boxtimes_{P(z)} W_2 \rightarrow W_3$
associated with $\Y$. We denote it as $m_{\Y}^{P(z)}$. We also
denote $m_{\Y}^{P(1)}$ simply as $m_{\Y}$. For $z_1, z_2>0$, let
$\gamma$ be a path in $\R_+$ from $z_1$ to $z_2$, then we have
\beq \label{m-Y-T}
m_{\Y}^{P(z_2)} \circ \mathcal{T}_{\gamma} = m_{\Y}^{P(z_1)}.
\eeq

\begin{prop} \label{Omega-Y-R}
$m_{\Omega_0 (\Y)} = m_{\Y} \circ \mathcal{R}_+$, 
$m_{\Omega_{-1} (\Y)} = m_{\Y} \circ \mathcal{R}_-$.
\end{prop} 
\pf
Let $u_1\in W_1, u_2\in W_2$ and $\Y_0$ be the 
intertwining operator associated with intertwining map
$\boxtimes: W_1\otimes W_2 \rightarrow W_1\boxtimes W_2$.
Let $\gamma_+$ be the path given in (\ref{gamma_+-fig}). 
Then we have  
\bea
\overline{m_{\Y}} \circ \overline{\mathcal{R}_+}(u_2 \boxtimes u_1)
&=& \overline{m_{\Y}}  
(e^{L(-1)} \overline{\mathcal{T}}_{\gamma_+} (u_1\boxtimes_{P(-1)} u_2))
\nn
&=& e^{L(-1)} \overline{m_{\Y}}(\Y_0(u_1, e^{\pi i})u_2)  \nn
&=& e^{L(-1)} \overline{m_{\Y}}(e^{\pi iL(0)} (e^{-\pi iL(0)}u_1
\boxtimes e^{-\pi iL(0)} u_2) \nn
&=& e^{L(-1)} e^{\pi iL(0)} \Y(e^{-\pi iL(0)}u_1, 1) e^{-\pi iL(0)} u_2 \nn
&=& e^{L(-1)} \Y(u_1, e^{\pi i})u_2  \nn
&=& \Omega_0(\Y) (u_2, 1) u_1   \nn
&=& m_{\Omega_0(\Y)} (u_2\boxtimes u_1).
\eea
Hence we obtain the first identity. 
The proof of the second identity is exactly same. 
\epf

The following Theorem is proved in \cite{HK2}.
\begin{thm}  \label{ffa-thm}
A conformal full field algebra over 
$V^L\otimes V^R$ is equivalent to a module $F$ 
for the vertex operator algebra $V^L\otimes V^R$ equipped with
an intertwining operator
$\Y$ of type $\binom{F}{FF}$ and an injective 
linear map $\rho: V^L\otimes V^R \to F$, 
satisfying the following conditions:
\begin{enumerate}

\item The {\it identity property}: 
$\Y(\rho(\one^L\otimes \one^R), x)=I_{F}$.

\item The {\it creation property}: For $u\in F$, 
$\lim_{x\rightarrow 0} \Y(u, x)\rho(\one^L\otimes \one^R)=u$.

\item The {\it associativity}: for $u,v,w\in F$ and $w'\in F'$,  
\bea
&&\langle w', \Y(u; z_1,\zeta_1)\Y(v; z_1, \zeta_1)w\rangle \nn
&&\hspace{1cm}=\langle w', \Y(\Y(u; z_1-z_2, \zeta_1-\zeta_2)v;
z_2, \zeta_2)w\rangle  \label{asso-ffa-ioa}
\eea
when $|z_{1}|>|z_{2}|>0$ and $|\zeta_{1}|>|\zeta_{2}|>0$.

\item The {\it single-valuedness property}: 
\begin{equation} \label{sing-val-1}
e^{2\pi i (L^L(0)-L^R(0))} = I_{F}.
\end{equation}   

\item The {\it skew symmetry}: 
\begin{equation} \label{skew}
\mathbb{Y}^{\Y}(u, 1, 1)v = e^{L^L(-1)+L^R(-1)}
\mathbb{Y}^{\Y}(v, e^{\pi i}, e^{-\pi i})u.
\end{equation}
\end{enumerate}
\end{thm}

\begin{cor} \label{ffa-ovoa-cor}
The category of conformal full field algebras over $V^L\otimes V^R$ 
is isomorphic to the category of 
open-string vertex operator algebras 
which contain the vertex operator algebra 
$V^L\otimes V^R$ in its meromorphic center 
and satisfy the single-valuedness condition (\ref{sing-val-1}) 
and the skew symmetry (\ref{skew}). 
\end{cor}
\pf
The proof is almost obvious, except that
the associativity axiom in Theorem 
\ref{ffa-thm} might looks stronger than that 
of open-string vertex operator algebra, 
which only requires (\ref{asso-ffa-ioa}) to hold
for $z_1=\zeta_1>0$ and $z_2=\zeta_2>0$. However, 
under the strong condition that both $V^L$ and 
$V^R$ satisfy the condition in Theorem \ref{ioa}
and the splitting property of $\Y$ in this case, 
one can show that these two associativities are
actually equivalent. 
\epf

The category of $V^L$- ($V^R$-) modules, 
denoted as $\mathcal{C}_{V^L}$ ($\mathcal{C}_{V^R}$), 
has a natural structure of vertex tensor category
and semisimple braided tensor category. 
The unit isomorphisms, the associativity isomorphisms 
and the braiding isomorphisms in $\mathcal{C}_{V^L}$ 
and $\mathcal{C}_{V^R}$ are given in the same way as
(\ref{l-unit-cat}), (\ref{r-unit-cat}), (\ref{assoc-iso}),
(\ref{R-+-cat}) and (\ref{R---cat}). We denote them as
$l^L, r^L, \mathcal{A}^L, \mathcal{R}_{\pm}^L$
and $l^R, r^R, \mathcal{A}^R, \mathcal{R}_{\pm}^R$ respectively.

Now we consider the category of 
$V^L\otimes V^R$-modules, denoted as 
$\mathcal{C}_{V^L\otimes V^R}$.  
It was proved in \cite{HK2} that
$V^L\otimes V^R$ also satisfies the condition in Theorem 
\ref{ioa}. Thus by Theorem \ref{ioa}, the category 
$\mathcal{C}_{V^L\otimes V^R}$ has a structure of 
semisimple braided tensor category.

We would like to take a closer look at the 
tensor product in $\mathcal{C}_{V^L\otimes V^R}$. 
Let $A=A^{L}\otimes A^{R}$, $B=B^{L}\otimes B^{R}$
and $C=C^{L}\otimes C^{R}$ be $V^L\otimes V^R$-modules.
Let $\Y^L \in \mathcal{V}_{A^LB^L}^{C^L}$ and 
$\Y^R \in \mathcal{V}_{A^RB^R}^{C^R}$. Let 
$\sigma : 
A\otimes B \rightarrow (A^{L}\otimes B^{L})\otimes 
( A^{R}\otimes B^{R} )$
be the map defined as 
$$
\sigma: (a^L\otimes a^R)\otimes (b^L\otimes b^R)
\mapsto (a^L\otimes b^L)\otimes (a^R\otimes b^R)
$$
for $a^L\in A^L, a^R\in A^R, b^L\in B^L, b^R\in B^R$. 
The following result is due to Dong, Mason and Zhu \cite{DMZ}. 

\begin{prop}  \label{prop-DMZ}
$(\Y^L\otimes \Y^R)\circ \sigma \in \mathcal{V}_{AB}^C$ and 
$\mathcal{V}_{A^LB^L}^{C^L}\otimes \mathcal{V}_{A^RB^R}^{C^R} \cong 
\mathcal{V}_{AB}^C$ canonically. 
\end{prop}

For $z\in \C^{\times}$, the $P(z)$-tensor product (\cite{HL1}-\cite{HL4}) 
of any two $V^L\otimes V^R$-modules
$A$ and $B$, is a $V^L\otimes V^R$-module $A\boxtimes_{P(z)} B$, 
together with a 
$P(z)$-intertwining map $\boxtimes_{P(z)}$ from 
$A\otimes B$ to $\overline{A\boxtimes_{P(z)} B}$, 
satisfying the universal property 
that given any $V^L\otimes V^R$-modules $C$ and 
a $P(z)$-intertwining map 
$G: A\otimes B \rightarrow \overline{C}$, 
there is a unique module map 
$\eta: A\boxtimes_{P(z)} B \rightarrow C$ 
such that $G= \overline{\eta} \circ \boxtimes_{P(z)}$. 

\begin{lemma} \label{lemma-1}
Let $A=A^L\otimes A^R$ and $B=B^L\otimes B^R$. Then 
$(A^L\boxtimes_{P(z)} B^L) \otimes (A^R\boxtimes_{P(z)} B^R)$ together with 
the $P(z)$-intertwining map 
$(\boxtimes_{P(z)} \otimes \boxtimes_{P(z)}) \circ \sigma$ from 
$A\otimes B$ to 
$\overline{(A^L\boxtimes_{P(z)} B^L) \otimes (A^R\boxtimes_{P(z)} B^R)}$
is a $P(z)$-tensor product of $A$ and $B$. 
\end{lemma}
\pf
We will only show the case $z=1$. 
The proof of general case is the same. 
Given any $V^L\otimes V^R$-module $C$ and a decomposition 
$C=\coprod_{i=1}^n C_i^L\otimes C_i^R$ or equivalently a family 
of morphism $\pi_i: C \rightarrow C_i^L\otimes C_i^R$ and 
$\iota_i: C_i^L\otimes C_i^R \rightarrow C$ satisfying 
the following conditions
\beq
\pi_i \circ \iota_i = I_{ C_i^L\otimes C_i^R}; 
\quad \pi_j\iota_i =0, \mbox{for $i\neq j$};
\quad \sum_i \iota_i\pi_i = I_C.
\eeq 
Assume $f$ is a $P(1)$-intertwining map from 
$(A^L\otimes A^R)\otimes (B^L\otimes B^R)$ 
to $\overline{C}$. It induces a 
$P(1)$-intertwining map $f_i: A\otimes B\rightarrow 
\overline{C_i^L\otimes C_i^R}$ for each $i$.  
By Proposition \ref{prop-DMZ}, 
$f_i = \sum_k (f_k^L \otimes f_k^R)\circ \sigma $ where 
$f_k^L, f_k^R$ are intertwining maps
$A^L\otimes B^L \rightarrow \overline{C_i^L}$, 
$A^R\otimes B^R \rightarrow \overline{C_i^R}$ respectively. 
By the universal property of 
$\boxtimes$, we know that there is a unique module map 
module map $\eta_{i}=\sum_{k}\eta_{i;k}^L\otimes \eta_{i;k}^R$, where
$\eta_{i;k}^L: A^L\boxtimes B^L \rightarrow C_i^L$ 
and $\eta_{i;k}^R: A^R\boxtimes B^R \rightarrow C_i^R$, 
such that the following diagram  
\beq  \label{lemma-diag}
\xymatrix{
A\otimes B 
\ar[rr]^{\hspace{-1.5cm}(\boxtimes \otimes \boxtimes) \circ \sigma} \ar[d]_{f}  
\ar[drr]^{f_i} &  &  
(A^L\boxtimes B^L)\otimes (A^R\boxtimes B^R) 
\ar@{.>}[d]^{\exists ! \, \eta_i} \\
C \ar[rr]^{\pi_i}   & & C_i^L\otimes C_i^R 
}
\eeq
is commutative. By the universal property of
 direct product, there is a unique
module map $\eta$ from $(A^L\boxtimes B^L)\otimes (A^R\boxtimes B^R)$
to $C$ such that $\pi_i \circ \eta = \eta_i$. 
Hence we have 
\beq  \label{lemma-int-equ-1}
\pi_i \circ f =
\pi_i \circ \eta \circ (\boxtimes \otimes \boxtimes) \circ \sigma
\eeq
Composing two sides of above equation 
with $\iota_i$ from left and sum up all $i$, 
we obtain the equality
\beq  \label{F-eta}
f= \eta \circ (\boxtimes \otimes \boxtimes) \circ \sigma.
\eeq

Moreover, the solution of $\eta$ satisfying (\ref{F-eta}) is 
unique because each $\pi_i\eta$ gives the unique
solution of (\ref{lemma-int-equ-1})
and $C$ as a direct product satisfies the universal property.

We have proved that $(A^L\boxtimes B^L)\otimes (A^R\boxtimes B^R)$
together with $(\boxtimes \otimes \boxtimes) \circ \sigma$ 
satisfies the universal properties of tensor product, 
thus gives a $P(1)$-tensor product of $A$ and $B$. 
\epf

There are only finite number of equivalent classes 
of simple objects in $\mathcal{C}_{V^L\otimes V^R}$. 
Simple objects in $\mathcal{C}_{V^L\otimes V^R}$ are 
objects of form $W^L\otimes W^R$, 
where $W^L$ and $W^R$ are simple $V^L$-module
and simple $V^R$-module respectively \cite{FHL}.  
Every object in $\mathcal{C}_{V^L\otimes V^R}$ is a 
direct sum of simple objects. If two objects $A$ and $B$ has 
the following decomposition: 
\beq  \label{ob-in-C-2}
A=\coprod_{i=1}^m A_{i}^{L}\otimes A_{i}^{R}, \quad
B=\coprod_{j=1}^n B_{j}^{L}\otimes B_{j}^{R},
\eeq 
then by Lemma \ref{lemma-1},  
$\coprod_{i,j}(A_i^L\boxtimes B_j^L) \otimes (A_i^L\boxtimes B_j^L)$
is a tensor product of $A$ and $B$. 
Therefore the universal property of tensor product provides
a canonical isomorphism: 
\beq \label{can-iso}
\coprod_{i,j}(A_i^L\boxtimes B_j^L) \otimes  
(A_i^L\boxtimes B_j^L) \cong A\boxtimes B.
\eeq

By Huang's construction given in Theorem \ref{ioa}, 
there is a braiding structure
in $\mathcal{C}_{V^L\otimes V^R}$. 
However, in this work, we are interested
in a different braiding structure in $\mathcal{C}_{V^L\otimes V^R}$. 
For each $A$ and $B$ as in (\ref{ob-in-C-2}), 
we define $\mathcal{R}_{+-}$ by the following commutative diagram:
\beq  \label{R-+--def}
\xymatrix{
\coprod_{i=1}^{M}\coprod_{j=1}^{N} 
(A_i^L\boxtimes B_j^L)\otimes (A_i^R\boxtimes B_j^R)
\ar@<0.4ex>[rr]^{\hspace{2cm} \cong} 
\ar[d]_{\mathcal{R}_+^L \otimes \mathcal{R}_-^R}   
&  & A\boxtimes B
\ar@{.>}[d]^{\exists !\, \, \mathcal{R}_{+-} }  \\
\coprod_{i=1}^{M}\coprod_{j=1}^{N} 
(B_j^L\boxtimes A_i^L) \otimes (B_j^R\boxtimes A_i^R) 
\ar@<0.4ex>[rr]^{\hspace{2cm}\cong} &  & B\boxtimes A  ,
}
\eeq
where the two horizontal maps are
the canonical isomorphisms as given in (\ref{can-iso}). 
$\mathcal{R}_+^L \otimes \mathcal{R}_-^R$ is an isomorphism. 
So is $\mathcal{R}_{+-}$.

\begin{rema}
{\rm 
If we replace the $\mathcal{R}_+^L\otimes \mathcal{R}_-^R$ 
in the diagram above 
by $\mathcal{R}_+^L\otimes \mathcal{R}_+^R$, or 
$\mathcal{R}_-^L\otimes \mathcal{R}_+^R$, or 
$\mathcal{R}_-^L\otimes \mathcal{R}_-^R$, 
we will obtain a new morphism in each case. 
We denoted them as $\mathcal{R}_{++}$,
$\mathcal{R}_{-+}$ and $\mathcal{R}_{--}$ respectively. They are
all isomorphisms. 
}
\end{rema}

\begin{rema}
{\rm 
The horizontal isomorphisms in the diagram (\ref{R-+--def})
are induced from universal property of tensor product. 
As a consequence,  
such defined $\mathcal{R}_{\pm\pm}$ and 
$\mathcal{R}_{\pm\mp}$ are all independent of 
decompositions. 
}
\end{rema}

\begin{prop}
Each of $\mathcal{R}_{++}$, $\mathcal{R}_{+-}$, $\mathcal{R}_{-+}$ and
$\mathcal{R}_{--}$ gives $\mathcal{C}_{V^L\otimes V^R}$ 
a structure of braided tensor category.  
\end{prop}
\pf
It is amount to show that all four isomorphisms 
satisfy both the functorial properties and the hexagon relations. 
Because the braiding isomorphism is naturally induced from 
$\mathcal{R}_+^L\otimes \mathcal{R}_-^R$, it is routine to check that
both the functorial properties and the hexagon relations follow from
those properties of $\mathcal{R}_+^L$ and $\mathcal{R}_-^R$. 
\epf

\begin{rema}
{\rm 
It is not hard to see that the braiding isomorphism constructed 
in Theorem \ref{ioa} is just $\mathcal{R}_{++}$.
}
\end{rema}

From now on, we will only consider the 
braiding tensor category structure of 
$\mathcal{C}_{V^L\otimes V^R}$ given by $\mathcal{R}_{+-}$. 
In order to emphasis
our choice of braiding, sometimes we will denote the 
category $\mathcal{C}_{V^L\otimes V^R}$ as 
$(\mathcal{C}_{V^L\otimes V^R}, \mathcal{R}_{+-})$.

For any object $A\in \mathcal{C}_{V^L\otimes V^R}$, we 
also define an isomorphism $\theta_A: A\rightarrow A$, called
twist, as follow
\beq  \label{twist-theta}
\theta_A = e^{2\pi iL^L(0)} \otimes e^{-2\pi iL^R(0)}.
\eeq

Let $A=\coprod_i A_i^L\otimes A_i^R$, 
$B=\coprod_j B_j^L\otimes B_j^R$, $C=\coprod_k C_k^L\otimes C_k^R$ 
be $V^L\otimes V^R$-modules. 
$\Omega_0\otimes \Omega_{-1}$ acts on the space 
$\coprod_{i,j,k} \mathcal{V}_{A_i^LB_j^L}^{C_k^L} \otimes 
\mathcal{V}_{A_i^RB_j^R}^{C_k^R}$ naturally as an automorphism. 
By the canonical isomorphism given in 
Proposition \ref{prop-DMZ}, $\coprod_{i,j,k} 
\mathcal{V}_{A_i^LB_j^L}^{C_k^L} 
\otimes \mathcal{V}_{A_i^RB_j^R}^{C_k^R}$ also canonically
isomorphic to $\mathcal{V}_{AB}^C$. 
Therefore, we obtain an action of
$\Omega_0\otimes \Omega_{-1}$ on $\mathcal{V}_{AB}^C$
as an automorphism. 

Let $\Y\in \mathcal{V}_{AB}^C$ and 
$m_{\Y}$ the unique module map $A\boxtimes B \rightarrow C$
associated with $\Y$. 
By Proposition \ref{Omega-Y-R}, it is easy to obtain the following
Lemma: 
\begin{lemma}\label{omega-pm-Y-R-pm}
$m_{\Omega_0\otimes \Omega_{-1}(\Y) }  = m_{\Y} \circ \mathcal{R}_{+-}$.
\end{lemma}

\begin{thm} \label{cfa-thm}
The following two notions are
equivalent in the sense that the categories 
given by these notions are isomorphic: 
\begin{enumerate}
\item A conformal full field algebra over $V^L\otimes V^R$, 
$(F, m, \rho)$. 
\item A commutative associative algebra $(F, m, \iota)$ in
the braided tensor category 
$(\mathcal{C}_{V^L\otimes V^R}, \mathcal{R}_{+-})$ 
satisfying $\theta_{F}=I_{F}$.
\end{enumerate}
\end{thm}
\pf
By the Corollary \ref{ffa-ovoa-cor} and the result in \cite{HK1},
the first three conditions in Theorem 
\ref{ffa-thm} exactly amounts to an associative algebra in 
$\mathcal{C}_{V^L\otimes V^R}$. In particular, data in the these
two structures are related as follow: $\rho=\iota$
and $m_{\mathbb{Y}}=m$, where
$\mathbb{Y}_f$ is the formal 
vertex operator associated with a conformal full field algebra
and is an intertwining operator of $V^L\otimes V^R$ in this case.

It is also obvious that single-valuedness property (\ref{sing-val-1})
is equivalent to $\theta_F= I_F$. 
Notice that the skew symmetry (\ref{skew}) is equivalent to
the following condition: 
$$
\Omega_{0}\otimes \Omega_{-1} (\mathbb{Y}_f) = \mathbb{Y}_f. 
$$
By Lemma \ref{omega-pm-Y-R-pm}, it is manifest that 
the skew symmetry (\ref{skew}) exactly 
amounts to the commutativity of the corresponding 
associative algebra $F$ in $\mathcal{C}_{V^L\otimes V^R}$. 
\epf

\renewcommand{\theequation}{\thesection.\arabic{equation}}
\renewcommand{\thethm}{\thesection.\arabic{thm}}
\setcounter{equation}{0}
\setcounter{thm}{0}

\section{Frobenius algebras in $\mathcal{C}_{V^L\otimes V^R}$}
 
In this section, we will give a categorical formulation of
a conformal full field algebra over $V^L\otimes V^R$
equipped with a nondegenerate invariant bilinear form.

\subsection{$\tilde{A}_r$ and $\hat{A}_r$}

In this subsection, we fix a vertex operator algebra $V$.  

Given an intertwining operator ${\cal Y}$ of type 
$\binom{W_{3}}{W_{1}W_{2}}$ of $V$ and an integer $r\in {\Z}$, 
the so called 
{\it $r$-contragredient operator of ${\cal Y}$} (\cite{HL2}) 
was defined to be the linear map
\begin{eqnarray}
W_1 \otimes  W'_3& \rightarrow&   W'_2\{z\}\nno\\
w_{(1)}\otimes w'_{(3)} &\mapsto  &A_{r}({\cal Y})(w_{(1)},x)w'_{(3)}
\end{eqnarray}
given by
\begin{eqnarray} \label{r-contra-1}
\lefteqn{\langle A_{r}({\cal Y})(w_{(1)},x)w'_{(3)},w_{(2)}\rangle=}
\nno\\
&&= \langle w'_{(3)}, {\cal Y}(e^{xL(1)}e^{(2r+1)\pi i L(0)}
x^{-2L(0)}w_{(1)},x^{-1})w_{(2)}
\rangle,
\end{eqnarray}
where $w_{(1)}\in W_{1}$, $w_{(2)}\in W_{2}$, $w'_{(3)}\in W'_{(3)}$. 
The following proposition was proved in \cite{HL4}. 
\begin{prop}  \label{A-prop}
The $r$-contragredient operator 
$A_{r}({\cal Y})$ of an intertwining operator
${\cal Y}$ of type $\binom{W_3}{W_1 W_2}$ is an
intertwining operator of type  $\binom{W'_2}{W_1 W'_3}$.
Moreover,
\begin{equation} \label{AA-id}
A_{-r-1}(A_{r}({\cal Y}))=A_{r}(A_{-r-1}({\cal Y}))={\cal Y}.
\end{equation}
In particular, the correspondence
${\cal Y} \mapsto  A_{r}({\cal Y})$ defines a linear isomorphism {from}
${\cal V}^{W_{3}}_{W_{1}W_{2}}$ to  ${\cal V}^{W'_{2}}_{W_{1}W'_{3}}$,
and we have
$N^{W_{3}}_{W_{1}W_{2}}= N^{W'_{2}}_{W_{1}W'_{3}}$.
\end{prop}

It turns out that $A_r$ is very hard to work with 
in the tensor category because they are simply not
categorical. To fix the problem, we will
introduce two slightly different operators.

Given an intertwining operator ${\cal Y}$ of type $\binom{W_{3}}
{W_{1}W_{2}}$ and an integer $r\in {\Z}$, we define two 
operators $\tilde{A}_r(\Y)$ and $\hat{A}_r(\Y)$ as
\bea \label{two-new-A-1}
&&\tilde{A}_r(\Y)( \cdot, x) = e^{-(2r+1)\pi i L(0)}A_r(\Y)(\cdot, x)e^{(2r+1)\pi i L(0)},  \nn
&&\hat{A}_r(\Y)(\cdot, x) = e^{-(2r+1)\pi i L(0)}A_{-r-1}(\Y)(\cdot, x)e^{(2r+1)\pi i L(0)}, 
\eea
or equivalently, 
\bea  \label{two-new-A-2}
&&\langle \tilde{A}_{r}(\Y)(w_{1}, e^{(2r+1)\pi i}x)w'_{3},w_{2}\rangle 
=\langle w'_{3}, \Y(e^{xL(1)}x^{-2L(0)}w_{1},  x^{-1})w_{2}\rangle,  
\nn
&&\langle \hat{A}_{r}(\Y)(w_{1}, x)w'_{3},w_{2}\rangle 
=\langle w'_{3}, \Y(e^{-xL(1)}x^{-2L(0)}w_{1},  
e^{(2r+1)\pi i}x^{-1})w_{2}\rangle,
\eea
for $w_{1}\in W_{1}$, $w_{2}\in W_{2}$, $w'_{3}\in W'_{3}$. 
In particular, in the case $W_{1}=V$ and $W_{2}=W_{3}=W$,
$\tilde{A}_{r}(Y_W)$ and $\hat{A}_r(Y_W)$ related to 
contragredient vertex operator $Y_W'$ 
for any $r\in {\Z}$ as follow
\begin{equation} \label{contra-1-2}
\tilde{A}_r(Y_W)(\cdot, x) = \hat{A}_r(Y_W)(\cdot, x)= 
(-1)^{-L(0)}Y_W'(\cdot, x)(-1)^{L(0)}.
\end{equation}

\begin{lemma}  \label{lemma-mod-tilde-A}
Let $(W, Y_W)$ be a $V$-module. Let $\tilde{Y}_W$ be defined as  
$$
\tilde{Y}_W(\cdot, x) = a^{-L(0)}Y_W(\cdot, x) a^{L(0)}.
$$
for $a\in \C^{\times}$. 
Then $(W, \tilde{Y}_W)$ gives $W$ 
another module structure which is isomorphic to 
$(W, Y_W)$. Moreover the isomorphism is 
given by $w \mapsto a^{-L(0)}w$. 
\end{lemma}
\pf
It is clear that $\tilde{Y}_W(\one, x)=Y_W(\one, a^{-1}x) = I_{W}$.
For any $u,v\in V$ and $N\in \N$ large enough, we have
\bea
&&(x_0+x_2)^N \tilde{Y}_W(u, x_0+x_2)\tilde{Y}_W(v, x_2)w  \nn
&& \hspace{1cm} = (x_0+x_2)^N a^{-L(0)}Y_W(u, x_0+x_2)Y_W(v,x_2) a^{L(0)}w \nn
&& \hspace{1cm} = (x_0+x_2)^N a^{-L(0)}Y_W(Y(u, x_0)v, x_2)a^{L(0)}w \nn
&& \hspace{1cm} = (x_0+x_2)^N \tilde{Y}_W(Y(u,x_0)v, x_2)w.
\eea
The Jacobi identity follows from the weak associativity (\cite{LL}). Hence
$(W, \tilde{Y}_W)$ is a $V$-module. Now we show that it is isomorphic
to $(W, Y)$. Let $f: W\rightarrow W$ be 
so that $f(w)=a^{-L(0)}w$ for $w\in W$. 
Then we have 
$$
f(Y_W(u,x)w) = a^{-L(0)}Y_W(u,x)w  = a^{-L(0)}Y_W(u,x) a^{L(0)} f(w) 
= \tilde{Y}(u, x)f(w).
$$
Hence $f$ gives an isomorphism between two module structures on $W$. 
\epf

By the above Proposition, we see that $(W', \tilde{A}_r(Y_W))$ 
(or equivalently $(W', \hat{A}_r(Y_W))$) 
gives on $W'$ another module structure , which
is isomorphic to $(W', Y')$. 
A difference of these two module structures on $W'$ 
worth of knowing and frequently used in the
later sections is the following relation: 
\beq  \label{L-n-dual}
\langle w', L(n)w\rangle = (-1)^n \langle L'(n)w', w\rangle, 
\quad n\in \Z,
\eeq
for $w\in W, w'\in W'$.

The notion of intertwining operator depends on the choices of 
module structures on the three modules involved. 
The following results is an analogue of
Proposition \ref{A-prop}. 
\begin{prop}   \label{prop-tilde-A}
If $\Y$ is an intertwining operator of type 
$\binom{(W_3,Y_3)}{(W_1,Y_1)(W_2,Y_2)}$, 
then $\tilde{A}_r(\Y)$, $\hat{A}_r(\Y)$ are intertwining operators of type 
$\binom{(W_2',\tilde{A}_r(Y_2))}{(W_1,Y_1)(W_3',\tilde{A}_r(Y_3))}$ for each $r\in \Z$. 
Moreover, we have
\begin{equation} \label{inverse-two-A}
\tilde{A}_r \circ \hat{A}_r (\Y) = 
\hat{A}_r \circ \tilde{A}_r (\Y) = \Y.
\end{equation}
\end{prop}
\pf
It is routine to show that $\tilde{A}_r(\Y)$, $\hat{A}_r(\Y)$ 
are intertwining operators of type 
$\binom{(W_2',\tilde{A}_r(Y_2))}
{(W_1,Y_1)(W_3',\tilde{A}_r(Y_3))}$ for each $r\in \Z$.
 
We only prove (\ref{inverse-two-A}) here. 
First, we have
\bea
\langle \hat{A}_r(\tilde{A}_r(\Y)(u,x)w', w\rangle 
&=&\langle w',
\tilde{A}_r(\Y)(e^{-xL(1)}x^{-2L(0)}u, e^{(2r+1)\pi i}x^{-1}) w\rangle
\nn
&=&\langle \Y(e^{x^{-1}L(1)}x^{2L(0)}e^{-xL(1)}x^{-2L(0)}u, x) w',  
w\rangle \nn
&=& \langle \Y(u,x)w', w\rangle.   \nonumber
\eea
Second, we have 
\bea
&&\langle \tilde{A}_r(\hat{A}_r(\Y))(u,e^{(2r+1)\pi i}x)w', w\rangle \nn
&&\hspace{2cm}=
\langle w', \hat{A}_r(\Y)(e^{xL(1)}x^{-2L(0)}u, x^{-1})w\rangle \nn
&&\hspace{2cm}= \langle \Y(e^{-x^{-1}L(1)}x^{2L(0)}e^{xL(1)}x^{-2L(0)}u, 
e^{(2r+1)\pi i}x) w',  w\rangle \nn 
&&\hspace{2cm}= \langle \Y(u, e^{(2r+1)\pi i}x)w', w\rangle. \nonumber
\eea
Therefore we obtain (\ref{inverse-two-A}). 
\epf

The correspondence
$$
\Y \mapsto  \tilde{A}_{r}(\Y), \quad\quad \Y \mapsto \hat{A}_r(\Y)
$$ 
defines two linear isomorphisms:
\beq
{\cal V}^{(W_{3},Y_{W_3})}_{(W_{1}, Y_{W_1})(W_{2},Y_{W_2})} 
\rightarrow 
{\cal V}^{(W'_{2},\tilde{A}_r(Y_{W_2}))}_{(W_{1},Y_{W_1})(W'_{3},\tilde{A}_r(Y_{W_3}))}.
\eeq
Obviously, we have
\begin{equation}
N^{(W_{3},Y_{W_3})}_{(W_{1}, Y_{W_1})(W_{2},Y_{W_2})}
=N^{(W'_{2},\tilde{A}_r(Y_{W_2}))}_{(W_{1},Y_{W_1})(W'_{3},\tilde{A}_r(Y_{W_3}))}.
\end{equation}

In this work, we only use $\tilde{A}_r$ or $\hat{A}_r$ 
instead of $A_r$. Therefore, to simplify notations, 
we simply denote $(W', \tilde{A}_r(Y_W))$ as $W'$.

In \cite{MS}\cite{H9}, the $S_3$-action on the 
space of intertwining operators $\V$ are used. 
In this work, we will use $\Omega_r$, $\tilde{A}_r$ and
$\hat{A}_r$ instead. So we won't see a $S_3$ action anymore. 
However, the cyclic subgroup $\Z_3$ of $S_3$ still appear here. 
Let $\Y \in \V_{a_1a_2}^{a_3}$. It is easy to see that 
\bea
&&\langle e^{-xL(-1)} \Omega_r(\tilde{A}_r(\Y))(w'_{a_3}, x)w_{a_1},
w_{a_2}\rangle \nn
&&\hspace{3cm}=
\langle w'_{a_3}, \Y(e^{xL(1)}x^{-2L(0)}w_{a_1}, x^{-1})w_{a_2}\rangle
\eea
for $w_{a_1}\in W_{a_1}, w_{a_2}\in W_{a_2}, w'_{a_3}\in W'_{a_3}$. 
It is clear that $\Omega_r \circ \tilde{A}_r$ is independent
of $r\in \Z$. We denote it as $\sigma_{123}$. 
\begin{prop}
$\sigma_{123}^3 = I_{\V}$.
\end{prop}
\pf
Keep in mind the relation (\ref{L-n-dual}). 
For $w'_3\in W'_3, w_1\in W_1, w_2\in W_2$ and $\Y\in \V_{a_1a_2}^{a_3}$, 
we have
\bea
&&\langle w'_3, \sigma_{123}^3(\Y)(w_1, x) w_2\rangle \nn
&&\hspace{1cm}= \langle \sigma_{123}^{2}(\Y)(e^{xL(1)}x^{-2L(0)}w_2, x^{-1})
e^{-xL(1)} w'_3, w_1\rangle \nn
&&\hspace{1cm}= \langle e^{xL(1)}x^{-2L(0)} w_2, 
\sigma_{123}(\Y)(e^{x^{-1}L(1)}x^{2L(0)}e^{-xL(1)}w'_3, x)
e^{-x^{-1}L(1)}w_1\rangle  \nn
&&\hspace{1cm}= \langle e^{xL(1)} x^{-2L(0)} w_2, 
\sigma_{123}(\Y )(x^{2L(0)}w'_3, x) e^{-x^{-1}L(1)}w_1\rangle \nn
&&\hspace{1cm}=\langle \Y(e^{xL(1)}x^{-2L(0)}e^{-x^{-1}L(1)}w_1, x^{-1})
e^{-xL(1)}e^{xL(1)}x^{-2L(0)}w_2, x^{2L(0)}w'_3\rangle \nn
&&\hspace{1cm}=\langle \Y(x^{-2L(0)}w_1, x^{-1}) x^{-2L(0)}w_2, 
x^{2L(0)}w'_3\rangle  \nn
&&\hspace{1cm}= \langle \Y(w_1, x)w_2, w'_3\rangle.
\eea
\epf

The inverse of $\sigma_{123}$ can be expressed
as $\hat{A}_r \circ \Omega_{-r-1}$. We denote it 
as $\sigma_{132}$. Both $\sigma_{123}$ and $\sigma_{132}$ 
play similar roles as those in \cite{H9}.

\subsection{Modular tensor categories}

In this subsection, we review Huang's construction of 
duality maps \cite{H10} in terms of our new conventions.

From now on, we assume $V$ satisfies the condition in 
Theorem \ref{MTC}. By assumption, $V'\cong V$, i.e. $e'=e$.   
From \cite{FHL}, there is a nondegenerate invariant bilinear form 
$(\cdot, \cdot)$ on $V$ such that $(\one, \one)=1$. In the rest of 
this work, we will simply identify $V'$ with $V$ without
making the isomorphism explicit.

Let $\I$ be the set of equivalent classes of irreducible $V$-modules
and $W^a$ the chosen representative for $a\in \I$. 
Let $\{ \Y_{ea}^a \}$ be a basis of $\V_{ea}^a$ for $a\in \I$ such
that it coincides with the vertex operator $Y_{W^a}$, which defines
the $V$-module structure on $W^a$, i.e. $\Y_{ea}^a=Y_{W^a}$.
We choose a basis $\{ \Y_{ae}^a \}$ of $\V_{ae}^a$ as 
\beq
\Y_{ae}^a = \Omega_{-1} (\Y_{ea}^a).
\eeq
We also choose a basis $\{ \Y_{aa'}^e \}$ of $\V_{aa'}^e$ as 
\beq
\Y_{aa'}^{e}= \Y_{aa'}^{e'} = \hat{A}_{0} (\Y_{ae}^a) = \sigma_{132}(\Y_{ea}^a). 
\eeq
Notice that these choices are made 
for all $a\in \I$. In particular, we have
$$
\Y_{a'e}^{a'} = \Omega_{-1} (\Y_{ea'}^{a'}), \quad \quad 
\Y_{a'a}^e = \Y_{a'a}^{e'} = \hat{A}_{0} (\Y_{a'e}^{a'}). 
$$

A remark of our notations, for an arbitrary basis of 
$\V_{a_1a_2}^{a_3}$, which is priori different 
from above specific choices, we will use notation
$\{ \Y_{a_1a_2;i}^{a_3;(p)} \}_{i=1}^{N_{a_1a_2}^{a_3}}$  
(with additional subscript and superscript!) 
for $p\in \N$. For example $\{ \Y_{ea;1}^{a;(1)} \}$.

We denote the matrix entries
of a fusing matrix with respect to some arbitrary basis as
$$
F(\Y_{a_1a_5;i}^{a_4;(1)} \otimes \Y_{a_2a_3;j}^{a_5;(2)}, 
\Y_{a_6a_3;k}^{a_4; (3)}\otimes \Y_{a_1a_2;l}^{a_6;(4)}). 
$$ 
We also use the following simple notation:
$$
F_a = F(\Y_{ae}^{a}\otimes \Y_{a'a}^e; \Y_{ea}^a \otimes \Y_{aa'}^e),
\quad a\in \I.
$$
It is proved in \cite{H9} that $F_a\neq 0$ for all $a\in \I$.

Now we are ready to give the construction of 
the duality maps \cite{H10}. 
Since $\mathcal{C}_V$ is semisimple, we only need to
discuss irreducible modules. We start with the right duals. 
For $a\in \I$, we need define maps
$e_{a}: (W^a)' \boxtimes W^a \rightarrow V$
for all $a\in \I$. Our choice is 
\beq
e_{a} = m_{\Y_{a'a}^e}.
\eeq 
We also need to define map
$i_{a}: V \rightarrow W^a \boxtimes (W^a)'$. 
Although there is a submodule of $W^a \boxtimes (W^a)'$
isomorphic to $V$, there is no canonical isomorphism. 
But we know that 
$\dim_{\C} \hom_{V}(V, W^a \boxtimes (W^a)')=1$. 
Hence we choose $i_{a}$ to be the unique morphism such that
\beq  \label{m-i-def-R}
m_{\Y_{aa'}^e} \circ i_{a} = \frac{1}{F_a}\, I_{V}. 
\eeq
It is obvious that $m_{\Y_{aa';i}^{b;(p)}}\circ i_a =0$ for all 
$b\neq e$. 

For left duals, we define the map 
$e'_{a}: W^a \boxtimes (W^a)' \rightarrow V$
to be 
\beq
e'_{a}= m_{\Y_{aa'}^e},
\eeq 
and the map $i'_{a}: V \rightarrow (W^a)' \boxtimes W^a$
to be the unique morphism such that 
\beq \label{m-i-def-L}
m_{\Y_{a'a}^e} \circ i'_{a} = \frac{1}{F_a}\, I_{V}. 
\eeq
It is clear that $m_{\Y_{a'a;i}^{b;(p)}}\circ i'_a =0$
for all $b\neq e$. 

Huang's main results in \cite{H14}\cite{H9}\cite{H10} is stated as follow: 
\begin{thm}   \label{MTC}
Let $V$ be a simple vertex operator
algebra satisfying the following conditions:
\bnu
\item $V_{n}=0$ for $n<0$, $V_{(0)}=\C \one$ and 
$V'$ is isomorphic to $V$ as $V$-module, 
\item Every $\N$-gradable weak $V$-module is 
completely reducible, 
\item $V$ is $C_2$-cofinite.
\enu
Then $\mathcal{C}_V$, together with its monoidal structure and 
above duality maps, is a rigid braided tensor 
category and $\dim a= \frac{1}{F_{a}}$ for $a\in \I$. 
Moreover, $\mathcal{C}_V$, together with the
twist $\theta_W= e^{2\pi iL(0)}$ for each object $W$, 
has a structure of modular tensor category \cite{T}\cite{BK}. 
\end{thm}

\begin{rema} {\rm
Our choice of duality maps is slightly different
from that of Huang because in this work the module structure on 
$W'$ for each $V$-module is different from that used in \cite{H10}. 
In the appendix, we will give a detailed proof of rigidity. 
}
\end{rema}

There are powerful tools, called graphic calculus, 
in a modular tensor category. In particular, 
the right duality maps $i_a$ and $e_a$ can be 
denoted by the following graphs:
$$
\begin{picture}(14,2)
\put(2,1){$i_{a} =$} \put(3.6, 1.8){$a$}\put(5.6,1.8){$a'$}
\put(4,0){\resizebox{1.5cm}{2cm}
{\includegraphics{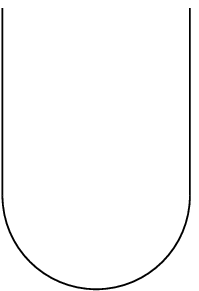}}}

\put(6.5,0){,}

\put(8,1){$e_{a} =$}\put(10,0)
{\resizebox{1.5cm}{2cm}{\includegraphics{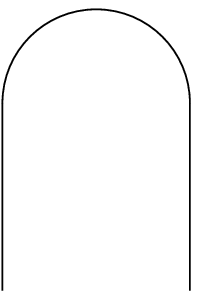}}}
\put(9.6, 0){$a'$}\put(11.6, 0){$a$} \put(12,0){,}
\end{picture}
$$
and the twist and its inverse, for any object $W$, 
are denoted by the following graphs
$$
\begin{picture}(14,2)
\put(2,1){$\theta_{W} =$} \put(3.4, 1.7){$W$}\put(3.4,0){$W$}
\put(4,0){\resizebox{0.3cm}{2cm}
{\includegraphics{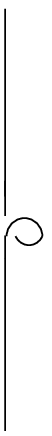}}}

\put(6.5,0){,}

\put(8,1){$\theta_{W}^{-1} =$}\put(9.4,1.7){$W$}\put(9.4,0){$W$}
\put(10,0){\resizebox{0.3cm}{2cm}
{\includegraphics{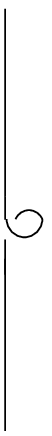}}}\put(11,0){.}

\end{picture}
$$

\begin{rema}
{\rm
In a ribbon category, given a right dual of an object $U^*$, 
there is automatically a left dual $^{*}U$ given as follows: 
\beq  \label{left-dual-rib}
\begin{picture}(14,2)
\put(0.5,1.8){$^{*}U$} \put(2.6, 1.8){$U$}
\put(1,0){\resizebox{1.5cm}{2cm}
{\includegraphics{left-dual-1.eps}}}
\put(3,1){$=$}
\put(3.3,1.8){$^{*}U$}\put(5.4, 1.8){$U$}
\put(3.8,0){\resizebox{1.5cm}{2cm}
{\includegraphics{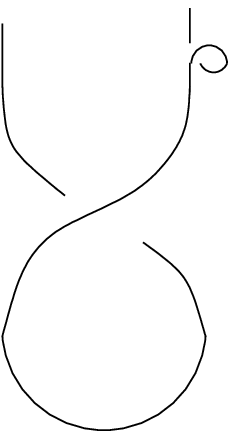}}}

\put(6,0){,}

\put(6.5, 0){$U$}\put(8.6,0){$^{*}U$}
\put(7,0){\resizebox{1.5cm}{2cm}{\includegraphics{left-dual-3.eps}}}
\put(9,1){$=$}
\put(9.4, 0){$U$}\put(11.4,0){$^{*}U$}
\put(9.8,0){\resizebox{1.5cm}{2cm}
{\includegraphics{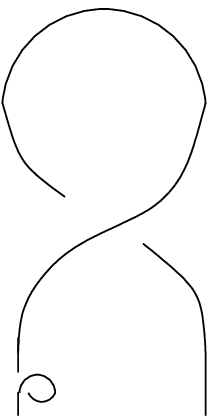}}}

\put(12,0){.}
\end{picture}
\eeq
In \cite{H10} and this work, the left duals and 
the right duals are constructed at the same time. 
It is certainly true that the left duals obtained from 
(\ref{left-dual-rib}) are compatible with our construction because
the left (or right) duals are unique up to isomorphisms. One can 
also see directly that the identity (\ref{3-equ-1}) proved
in the appendix is nothing but the following identity: 
\beq  \label{dual-twist}
\begin{picture}(14,2)
\put(2.6,0){$a'$}\put(4.6,0){$a$}
\put(3,0){\resizebox{1.5cm}{2cm}{\includegraphics{left-dual-3.eps}}}
\put(5,1){$=$}
\put(5.7,0){$a$}\put(7.6,0){$a'$}
\put(6,0){\resizebox{1.5cm}{2cm}{\includegraphics{left-dual-4.eps}}}
\put(8,1){$=$}
\put(8.4,0){$a$}\put(10.4,0){$a'$}
\put(8.8,0){\resizebox{1.5cm}{2cm}
{\includegraphics{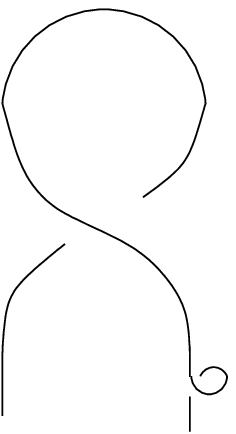}}}  \put(12,0){.}
\end{picture}
\eeq
The formula (\ref{dual-twist}) is implicitly 
used in many graphic calculations in this work. 
}
\end{rema}

\subsection{Categorical formulation of 
$\tilde{A}_0$ and $\tilde{A}_{-1}$}

Since $\mathcal{C}_V$ is a modular tensor category, we have a
powerful tool available. It is called graph calculus 
\cite{T}\cite{BK}. 
For this reason, we would like to 
express $\Omega_r$ and $\tilde{A}_r$ in terms of graphs. 

A basis $\{ \Y_{a_1a_2; i}^{a_3;(1)} \}_{i=1}^{N_{a_1a_2}^{a_3}}$ 
of $\V_{a_1a_2}^{a_3}$ for $a_1,a_2,a_3\in \I$
induces a basis in $\hom (a_1\boxtimes a_2, a_3)$, denoted
as $e_{a_1a_2;i}^{a_3}$ and as the following graph:
\beq  \label{basis-pic}
\begin{picture}(14, 1.8)
\put(6,0.2){\resizebox{1.5cm}{1.5cm}{\includegraphics{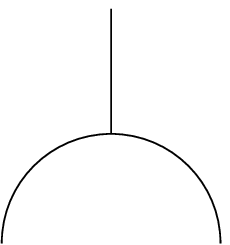}}}
\put(5.5, 0.2){$a_1$}\put(7.6, 0.2){$a_2$}\put(6.3, 1.5){$a_3$}
\put(6.6, 0.5){$i$}\put(8.5,0.2){.}
\end{picture}
\eeq
For simplicity, 
we always use $a$ to represent $W^a$ and $a'$ to represent
$(W^a)'$ in graphs. 

The map $\Omega_0: \V_{a_1a_2}^{a_3} \rightarrow \V_{a_2a_1}^{a_3}$ induces
a map $\hom_V(W^{a_1}\boxtimes W^{a_2}, W^{a_3}) \rightarrow
\hom_V(W^{a_2}\boxtimes W^{a_1}, W^{a_3})$, still denoted as $\Omega_0$.
By Proposition \ref{Omega-Y-R}, we have 
\beq
\begin{picture}(14, 1.8)
\put(2, 0.8){$ \Omega_0: $}
\put(4,0.2){\resizebox{1.5cm}{1.5cm}{\includegraphics{Y.eps}}}
\put(3.5, 0.2){$a_1$}\put(5.6, 0.2){$a_2$}\put(4.3, 1.5){$a_3$}
\put(4.6, 0.5){$i$}

\put(6.5, 0.8){$\mapsto$}

\put(9,0){\resizebox{1.5cm}{1.8cm}{\includegraphics{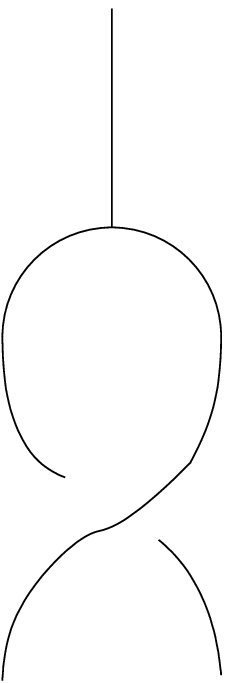}}}
\put(8.5, 0){$a_2$}\put(10.6, 0){$a_1$}\put(9.3, 1.7){$a_3$}
\put(9.7, 0.8){$i$}\put(11.5, 0){.}
\end{picture}
\eeq

\begin{lemma}
\beq  \label{L-z-1-2}
e^{(z_1-z_2)L(1)} z_1^{-2L(0)} e^{z_1L(-1)} e^{z_1^{-1}L(1)} =
z_2^{-2L(0)} e^{z_2L(-1)} e^{z_2^{-1}L(1)}
\eeq
\end{lemma}
\pf
It follows immediately from the following two identities:
\bea
e^{-(z_1-z_2)x^2\frac{d}{dx}} z_1^{2x\frac{d}{dx}} e^{-z_1\frac{d}{dx}} 
e^{-z_1^{-1}x^2\frac{d}{dx}} x &=& \frac{z_2x-1}{x} \nn
z_2^{2x\frac{d}{dx}} e^{-z_2\frac{d}{dx}} e^{-z_2^{-1}x^2\frac{d}{dx}}x
&=& \frac{z_2x-1}{x}.
\eea
\epf

\begin{prop}
By the universal property of tensor product, 
$\tilde{A}_0, \tilde{A}_{-1}: 
\V_{a_1a_2}^{a_3} \rightarrow \V_{a_1a'_3}^{a'_2}$ 
defined in (\ref{two-new-A-2}) also induce
two morphisms from $\hom_V(W^{a_1}\boxtimes W^{a_2}, W^{a_3})$ to 
$\hom_V(W^{a_1}\boxtimes (W^{a_3})', (W^{a_2})')$, still denoted as
$\tilde{A}_0, \tilde{A}_{-1}$. Then we have
\beq \label{tilde-A-0-graph}
\begin{picture}(14, 2)
\put(2, 0.8){$ \tilde{A}_0: $}
\put(4,0.2){\resizebox{1.5cm}{1.5cm}{\includegraphics{Y.eps}}}
\put(3.5, 0.2){$a_1$}\put(5.6, 0.2){$a_2$}\put(4.3, 1.5){$a_3$}
\put(4.6, 0.5){$i$}

\put(6.5, 0.8){$\mapsto$}
\put(9.5,0){\resizebox{2cm}{2.2cm}{\includegraphics{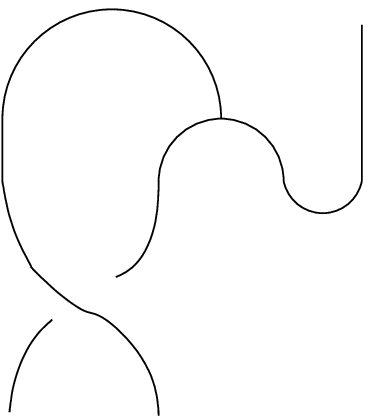}}}
\put(9, 0){$a_1$}\put(10.5, 0){$a'_3$}\put(11.6, 1.9){$a'_2$}
\put(10.6, 1.2){$i$}\put(12, 0){,}
\end{picture}
\eeq
\beq \label{tilde-A--1-graph}
\begin{picture}(14, 2)
\put(2, 0.8){$ \tilde{A}_{-1}: $}
\put(4,0.2){\resizebox{1.5cm}{1.5cm}{\includegraphics{Y.eps}}}
\put(3.5, 0.2){$a_1$}\put(5.6, 0.2){$a_2$}\put(4.3, 1.5){$a_3$}

\put(6.5, 0.8){$\mapsto$}
\put(9.5,0){\resizebox{2cm}{2.2cm}{\includegraphics{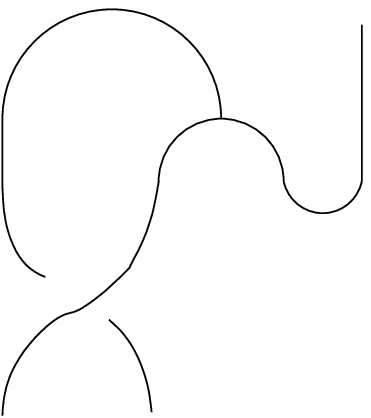}}}
\put(9, 0){$a_1$}\put(10.5, 0){$a'_3$}\put(11.6, 1.9){$a'_2$}
\put(10.6, 1.2){$i$}\put(12, 0){.}
\end{picture}
\eeq
\end{prop}
\pf
Let $\{ \Y_{a'_3a_1; i}^{a'_2;(2)} \}_{i=1}^{N_{a_1a_2}^{a_3}}$ 
be a basis of $\V_{(W^{a_3})' W^{a_1}}^{(W^{a_2})'}$ 
for $a_1,a_2,a_3\in \I$ such that 
the induced basis of $\hom ((W^{a_3})'\boxtimes W^{a_1}, (W^{a_2})')$ 
is given by the following graph:
\beq  
\begin{picture}(14, 2)
\put(6,0){\resizebox{2cm}{1.8cm}{\includegraphics{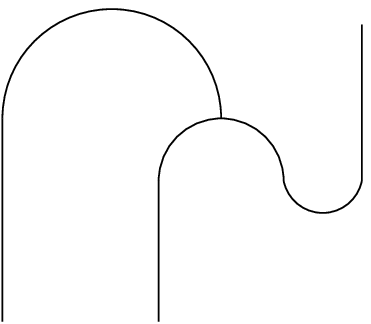}}}
\put(5.5, 0){$a'_3$}\put(7, 0){$a_1$}\put(8.1, 1.6){$a'_2$}
\put(7.2, 0.7){$i$}
\put(8.5, 0){.}
\end{picture}
\eeq
It is enough to show that 
\beq  \label{A-prop-equ-1}
\Y_{a'_3a_1; i}^{a'_2;(2)} = 
\Omega_0 (\tilde{A}_0 (\Y_{a_1a_2; i}^{a_3;(1)} )).
\eeq

The following identity is obvious.
\beq  
\begin{picture}(14, 2)
\put(4,0){\resizebox{2cm}{1.8cm}{\includegraphics{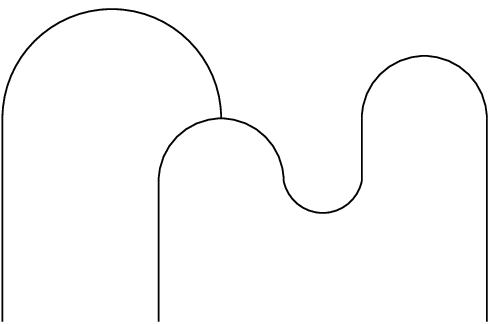}}}
\put(3.5, 0){$a'_3$}\put(4.8, 0){$a_1$}\put(6.1, 0){$a_2$}
\put(4.8, 0.7){$i$}

\put(7,1){$=$}

\put(8.5,0.2){\resizebox{2cm}{1.5cm}{\includegraphics{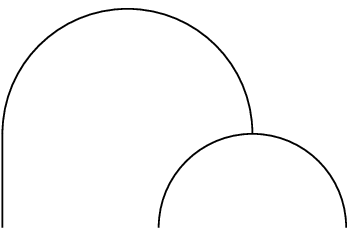}}}
\put(8.9, 0.2){$a_1$}\put(10.6, 0.2){$a_2$}\put(8, 0.2){$a'_3$}
\put(9.9, 0.5){$i$}\put(11,0.2){.}

\end{picture}
\eeq
It can be rewritten as 
\beq  \label{A-prop-equ-2-cat-1}
m_{\Y_{a'_2a_2}^e} \circ 
(m_{\Y_{a'_3a_1;i}^{a'_2; (2)}} \boxtimes I_{W^{a_2}} ) \circ
\mathcal{A} =  
m_{\Y_{a'_3a_3}^e}
\circ (I_{W^{a'_3}} \boxtimes m_{\Y_{a_1a_2;i}^{a_3;(1)}} ),
\eeq
where we have used the fact that $e_{a_k}=\dim a_k m_{\Y_{a'_ka_k}^e}, k=2,3$. 
Using (\ref{A-assoc}) and (\ref{m-Y-T}), it is easy to see that
\bea \label{A-prop-equ-2-cat-2}
&& m_{\Y_{a'_2a_2}^e}^{P(z_2)} \circ 
(m_{\Y_{a'_3a_1;i}^{a'_2; (2)}}^{P(z_1-z_2)} \boxtimes_{P(z_2)} I_{W^{a_2}} ) \circ
\mathcal{A}_{P(z_1),P(z_2)}^{P(z_1-z_2),P(z_2)}  \nn
&&\hspace{3cm} = m_{\Y_{a'_3a_3}^e}^{P(z_1)}
\circ (I_{W^{a'_3}} \boxtimes_{P(z_1)} m_{\Y_{a_1a_2;i}^{a_3;(1)}}^{P(z_2)} )
\eea
for $z_1>z_2>z_1-z_2>0$. Then (\ref{A-prop-equ-2-cat-2}) 
immediately implies the following identity: 
\bea  \label{A-prop-equ-2}
&&\Y_{a'_2a_2}^e ( 
\Y_{a'_3a_1;i}^{a'_2;(2)}(w_{a'_3}, z_1-z_2)w_{a_1}, z_2)w_{a_2} \nn 
&&\hspace{3cm} =
\Y_{a'_3a_3}^e(w_{a'_3}, z_1)\Y_{a_1a_2;i}^{a_3;(1)}(w_{a_1}, z_2) w_{a_2}.
\eea
for $w_{a_1}\in W^{a_1}, w_{a_2}\in W^{a_2}, w_{a'_3}\in (W^{a_3})'$.

On the one hand, we have 
\bea   \label{A-prop-equ-3}
&&\hspace{0cm}(\Y_{a'_2a_2}^e ( 
\Y_{a'_3a_1;i}^{a'_2;(2)}(w_{a'_3}, z_1-z_2)w_{a_1}, z_2)w_{a_2} , \one) \nn
&&\hspace{0.5cm}=
\langle \hat{A}_0(\Y_{a'_2e}^{a'_2}) ( 
\Y_{a'_3a_1;i}^{a'_2;(2)}(w_{a'_3}, z_1-z_2)w_{a_1}, z_2)w_{a_2} , \one \rangle 
\nn
&&\hspace{0.5cm}=
\langle w_{a_2}, \Y_{a'_2e}^{a'_2}(e^{-z_2L(1)}z_2^{-2L(0)}
\Y_{a'_3a_1;i}^{a'_2;(2)}(w_{a'_3}, z_1-z_2)w_{a_1}, e^{\pi i} z_2^{-1})\one
\rangle  \nn
&&\hspace{0.5cm}= \langle w_{a_2},e^{-z_2^{-1}L(-1)} e^{-z_2L(1)}z_2^{-2L(0)}
\Y_{a'_3a_1;i}^{a'_2;(2)}(w_{a'_3}, z_1-z_2)w_{a_1}\rangle.
\eea
On the other hand, we have
\bea  \label{A-prop-equ-4}
&&\hspace{-0.5cm}
(\Y_{a'_3a_3}^e(w_{a'_3}, z_1)\Y_{a_1a_2;i}^{a_3;(1)}(w_{a_1}, z_2) w_{a_2}, \one)
\nn
&&\hspace{0cm}=
\langle \hat{A}_0(\Y_{a'_3e}^{a'_3})(w_{a'_3}, z_1)
\Y_{a_1a_2;i}^{a_3;(1)}(w_{a_1}, z_2) w_{a_2}, \one\rangle \nn
&&\hspace{0cm}=
\langle \Y_{a_1a_2;i}^{a_3;(1)}(w_{a_1}, z_2) w_{a_2},
\Y_{a'_3e}^{a'_3}(e^{-z_1L(1)}z_1^{-2L(0)}w_{a'_3}, e^{\pi i}z_1^{-1})\one\rangle
\nn
&&\hspace{0cm}=
\langle \Y_{a_1a_2;i}^{a_3;(1)}(w_{a_1}, z_2) w_{a_2},
e^{-z_1^{-1}L(-1)}e^{-z_1L(1)}z_1^{-2L(0)}w_{a'_3}\rangle   \nn
&&\hspace{0cm}=
\langle z_1^{-2L(0)}e^{z_1L(-1)}e^{z_1^{-1}L(1)}
\Y_{a_1a_2;i}^{a_3;(1)}(w_{a_1}, z_2) w_{a_2}, w_{a'_3}\rangle \nn
&&\hspace{0cm}=
\langle \Y_{a_1a_2;i}^{a_3;(1)}( e^{(z_1-z_2)L(1)}(z_1-z_2)^{-2L(0)}w_{a_1}, 
(z_1-z_2)^{-1})  \nn
&&\hspace{4cm} \cdot z_1^{-2L(0)}e^{z_1L(-1)}e^{z_1^{-1}L(1)} 
w_{a_2}, w_{a'_3}\rangle \nn
&&\hspace{0cm}=
\langle z_1^{-2L(0)}e^{z_1L(-1)}e^{z_1^{-1}L(1)} w_{a_2},
\tilde{A}_0(\Y_{a_1a_2;i}^{a_3;(1)})(w_{a_1}, e^{\pi i}(z_1-z_2))
w_{a'_3}\rangle \nn
&&\hspace{0cm}=
\langle e^{(z_1-z_2)L(1)} z_1^{-2L(0)}e^{z_1L(-1)}e^{z_1^{-1}L(1)} w_{a_2},
\Omega_0(\tilde{A}_0(\Y_{a_1a_2;i}^{a_3;(1)}))(w_{a'_3}, z_1-z_2) w_{a_1}
\rangle \nn
&&\hspace{0cm}=
\langle z_2^{-2L(0)}e^{z_2L(-1)}e^{z_2^{-1}L(1)}w_{a_2},
\Omega_0(\tilde{A}_0(\Y_{a_1a_2;i}^{a_3;(1)}))(w_{a'_3}, z_1-z_2) w_{a_1}
\rangle,
\eea
where we have used (\ref{L-z-1-2}) in the last step.

Combining (\ref{A-prop-equ-2}), (\ref{A-prop-equ-3}) 
and (\ref{A-prop-equ-4}), we obtain (\ref{A-prop-equ-1}) 
immediately. 
\epf

Let us choose a basis $\{ f_{a_3;j}^{a_1a_2} \}_{j=1}^{N_{a_1a_2}^{a_3}}$ of 
$\hom_V(W^{a_3}, W^{a_1}\boxtimes W^{a_2})$, denoted as 
\beq  \label{dual-basis-pic}
\begin{picture}(14, 1.5)
\put(4, 0.7){$f_{a_3;j}^{a_1a_2} \, = $}
\put(6.8, 0){\resizebox{1.5cm}{1.5cm}{\includegraphics{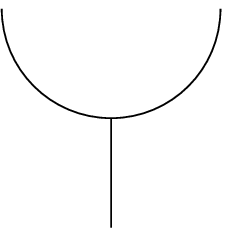}}}
\put(6.4, 1.2){$a_1$}\put(8.3, 1.2){$a_2$}\put(7, 0){$a_3$}
\put(7.4, 1){$j$}\put(8.5, 0){.}
\end{picture} 
\eeq
so that 
\beq \label{Y-dual-Y}
\begin{picture}(14,2)
\put(2,0.8){$\ds \frac{1}{\dim a_3}$}
\put(4, 0){\resizebox{2.5cm}{2cm}{\includegraphics{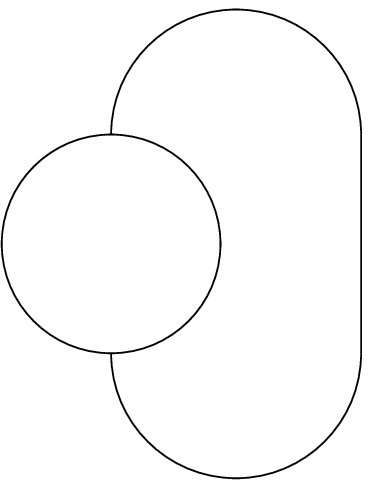}}}
\put(3.6, 0.8){$a_1$}\put(5.6, 0.8){$a_2$}\put(4.4, 1.8){$a_3$}
\put(4.4,0.1){$a_3$} \put(4.7, 1.13){$i$}\put(4.7, 0.7){$j$}
\put(7, 0.9){$=$} \put(8, 0.9){$\delta_{ij}.$} 
\end{picture}
\eeq

\begin{prop}
\beq \label{exp-id}
\begin{picture}(14,2)
\put(2.5, 0){$a_1$}\put(4.6,0){$a_2$}
\put(3, 0){\resizebox{1.5cm}{2cm}{\includegraphics{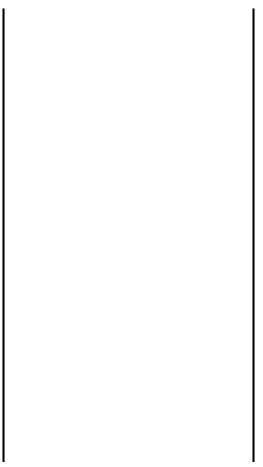}}}

\put(5.5, 1){$\ds =$}
\put(6.5, 1){$\ds \sum_{a_3\in \I} \sum_{i=1}^{N_{a_1a_2}^{a_3}}$}

\put(9, 0){\resizebox{1.5cm}{2cm}{\includegraphics{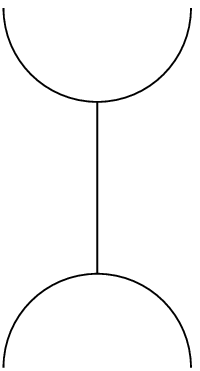}}}
\put(8.5,0){$a_1$}\put(10.6,0){$a_2$}
\put(8.5,1.8){$a_1$}\put(10.6,1.8){$a_2$}
\put(9.9,1){$a_3$}\put(9.7,0.2){$i$}\put(9.7,1.6){$i$}
\end{picture}
\eeq
\end{prop}
\pf
The basis $\{ e_{a_1a_2;i}^{a_3}\}_{i=1}^{N_{a_1a_2}^{a_3}}$ gives arise to
an isomorphism (\cite{HL1}-\cite{HL4}) 
$$
\phi_1:=\sum_{a_3}\sum_{i=1}^{N_{a_1a_2}^{a_3}} e_{a_1a_2;i}^{a_3} : \quad 
W^{a_1}\boxtimes W^{a_2}
\rightarrow \coprod_{a\in \I} \coprod_{i=1}^{N_{a_1a_2}^{a}} W^{a;(i)},
$$ 
where $W^{a;(i)}$ denotes the $i$-th copy of $W^a$ and 
$e_{a_1a_2;i}^{a}: W^{a_1}\otimes W^{a_2} \rightarrow W^{a;(i)}$. 

By the condition (\ref{Y-dual-Y}),
we have $e_{a_1a_2;i}^{a_3} \circ f^{a_1a_2}_{a_3;j} = \delta_{ij} I_{W^{a_3}}$.
Let 
$$
\phi_2:=\sum_{a_3}\sum_{j=1}^{N_{a_1a_2}^{a_3}} f^{a_1a_2}_{a_3;j}: \quad 
\coprod_{a\in \I} \coprod_{j=1}^{N_{a_1a_2}^{a}} 
W^{a;(j)} \rightarrow W^{a_1}\boxtimes W^{a_2}. 
$$ 
It is easy to see that we have 
$\phi_1\circ \phi_2 = I_{\coprod_{a\in \I} \coprod_{i=1}^{N_{a_1a_2}^{a}} W^{a;(i)}}$. 
Hence $\phi_2 = \phi_1^{-1}$. Then we also have 
$\phi_2 \circ \phi_1 = I_{W^{a_1}\boxtimes W^{a_2}}$, which implies
(\ref{exp-id}). 
\epf

Therefore $\{ f_{a_3;i}^{a_1a_2} \}$ can be viewed as the dual basis of  
$\{ \Y_{a_1a_2;i}^{a_3} \}_{i=1}^{N_{a_1a_2}^{a_3}}$. 
It is useful to figure out the 
action on $f_{a_3;j}^{a_1a_2;j}$ of $\hat{A}_0$ and $\hat{A}_{-1}$, 
which are the inverse of $\tilde{A}_0$ and $\tilde{A}_{-1}$ respectively.
The result is given in the following Proposition. 

\begin{prop}
\beq  \label{A-dual-1}
\begin{picture}(14,2)
\put(2, 0.8){$ \hat{A}_0^*: $}

\put(4, 0.2){\resizebox{1.5cm}{1.5cm}{\includegraphics{dual-Y.eps}}}
\put(3.6, 1.5){$a_1$}\put(5.6, 1.5){$a_2$}\put(4.3, 0.2){$a_3$}
\put(4.6, 1.2){$j$}

\put(6.5, 0.8){$\mapsto$}

\put(8,0.8){$\frac{\dim a_2}{\dim a_3}$}
\put(9.5,0){\resizebox{2cm}{2.2cm}{\includegraphics{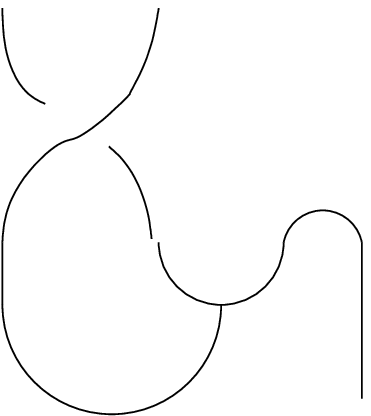}}}
\put(9.1, 2){$a_1$}\put(10.4, 2){$a'_3$}\put(11.6, 0){$a'_2$}
\put(10.6, 0.8){$j$}\put(12, 0){.}
\end{picture}
\eeq
\beq  \label{A-dual-2}
\begin{picture}(14,2)
\put(2, 0.8){$ \hat{A}_{-1}^*: $}

\put(4, 0.2){\resizebox{1.5cm}{1.5cm}{\includegraphics{dual-Y.eps}}}
\put(3.6, 1.5){$a_1$}\put(5.6, 1.5){$a_2$}\put(4.3, 0.2){$a_3$}
\put(4.6, 1.2){$j$}

\put(6.5, 0.8){$\mapsto$}

\put(8,0.8){$\frac{\dim a_2}{\dim a_3}$}
\put(9.5,0){\resizebox{2cm}{2.2cm}{\includegraphics{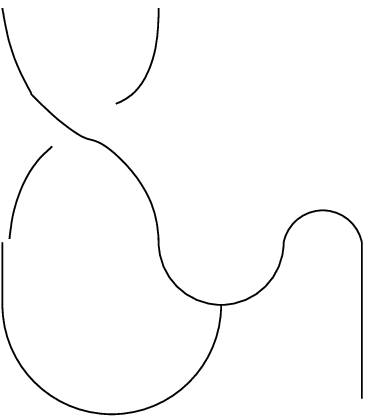}}}
\put(9.1, 2){$a_1$}\put(10.4, 2){$a'_3$}\put(11.6, 0){$a'_2$}
\put(10.6, 0.8){$j$}\put(12, 0){.}
\end{picture}
\eeq
\end{prop}
\pf
We will only prove (\ref{A-dual-1}). The proof of
(\ref{A-dual-2}) is completely an analogue. 

It is enough to show that the pairing between the 
image of (\ref{tilde-A-0-graph}) and (\ref{A-dual-1}) 
still gives $\delta_{ij}$. It is proved as follow:
\beq
\begin{picture}(14,2.5)
\put(1,1){$\frac{1}{\dim a'_2}\,\frac{\dim a_2}{\dim a_3}$}
\put(3.5, 0)
{\resizebox{2.5cm}{2.5cm}{\includegraphics{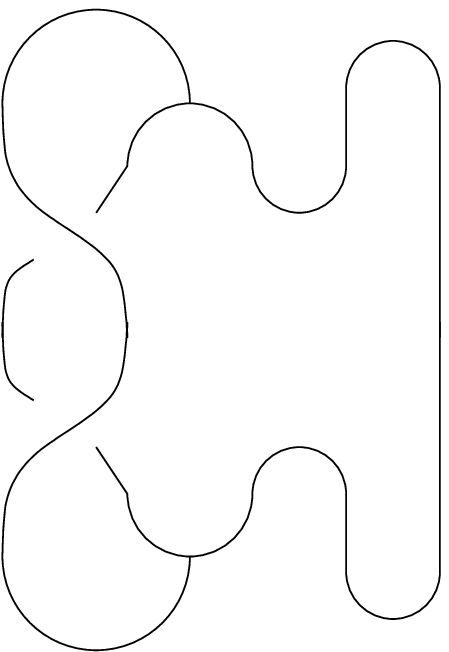}}}
\put(3.1, 1.2){$a_1$}\put(6.1, 1.5){$a_2$}\put(4.3, 1.2){$a'_3$}
\put(4.5, 0.6){$j$}\put(4.5, 1.7){$i$}
\put(1,1){$\frac{1}{\dim a'_2}\,\frac{\dim a_2}{\dim a_3}$}
\put(6.5, 1){$=$}

\put(7,1){$\frac{1}{\dim a_3}$}
\put(8.5,0)
{\resizebox{2cm}{2cm}{\includegraphics{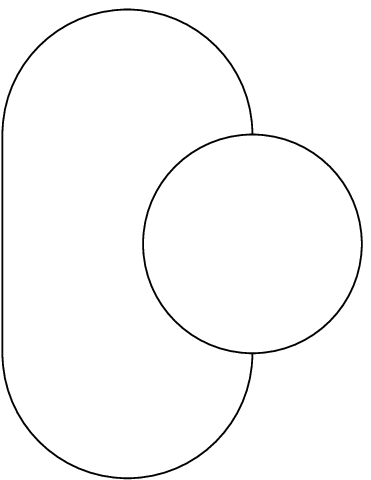}}}
\put(8.9, 1){$a_1$}\put(9.5, 2){$a_3$}\put(10.6, 1){$a_2$}
\put(9.8, 0.7){$j$}\put(9.8, 1.1){$i$}    
\put(11.5,1){$=\, \, \delta_{ij}.$}
\end{picture}
\eeq
\epf

\subsection{Frobenius algebras in $\mathcal{C}_{V^L\otimes V^R}$}

First let us recall the notion of coalgebra and 
Frobenius algebra (\cite{FS}) in a tensor category.  
\begin{defn}
{\rm 
A coalgebra $A$ in a tensor category 
$\mathcal{C}$ is an object with a
coproduct $\Delta\in \text{Mor}(A, A\otimes A)$ and a counit 
$\epsilon \in \text{Mor}(A, \one_{\mathcal{C}})$ such that 
\beq  \label{co-alg-def}
(\Delta \otimes I_{A}) \circ \Delta = 
(I_A \otimes \Delta) \circ \Delta, 
\hspace{1cm} (\epsilon \otimes I_{A})\circ \Delta 
= I_{A} = (I_A \otimes \epsilon)
\circ \Delta,
\eeq
which can also be expressed in term of the following
graphic equations:
\beq 
\begin{picture}(14,2)
\put(0.5,0){\resizebox{13cm}{2cm}{\includegraphics{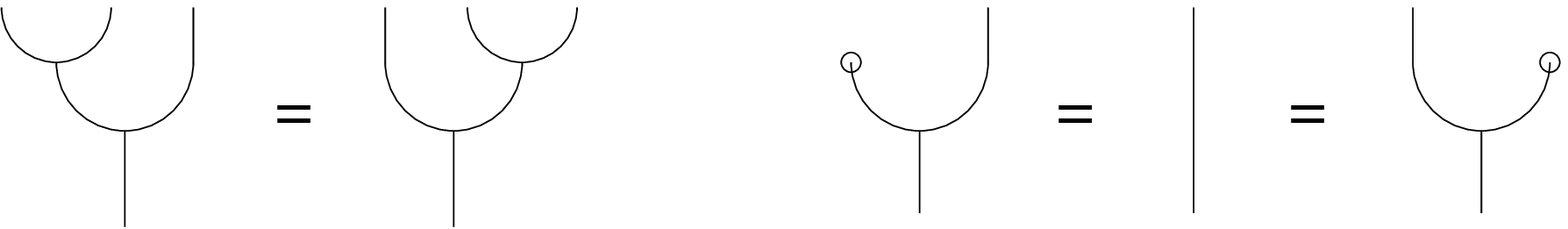}}}
\end{picture}  \nonumber
\eeq
}
\end{defn}

\begin{defn}
{\rm 
Frobenius algebra in $\mathcal{C}$ is 
an object that is both an algebra and a co-algebra
and for which the product and coproduct are related by 
\begin{equation}  \label{Frob-def-equ}
(I_A \otimes m) \circ (\Delta \otimes I_{A}) = \Delta \circ m = 
(m\otimes I_A) \circ (I_A \otimes \Delta),
\end{equation}
or as the following graphic equations, 
\beq \label{Frob-alg-def-fig}
\epsfxsize  0.8\textwidth
\epsfysize  0.2\textwidth
\epsfbox{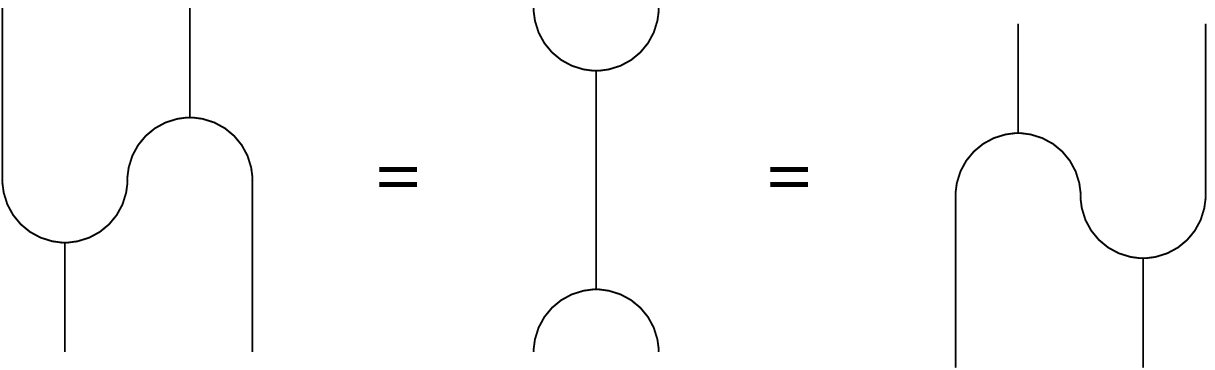}.
\eeq
}
\end{defn}

Let $(F,\mathbb{Y},\rho)$ be a conformal full field algebra  
over $V^L\otimes V^R$. Then 
$\mathbb{Y}_f(\cdot, x,x): F \rightarrow \edo F \, \{ x \}$ is 
an intertwining operator of type $\binom{F}{FF}$ of $V^L\otimes V^R$
viewed as vertex operator algebra. 
Let $\mathbb{Y}_f|_{V^L\otimes V^R}$ denote the restriction of  
$\mathbb{Y}_f(\cdot, x,x)$ on $\rho(V^L\otimes V^R)$. 
It defines the module structure on $F$. By Lemma
\ref{lemma-mod-tilde-A}, we know that 
$\tilde{A}_0\otimes \tilde{A}_{-1}(\mathbb{Y}_f|_{V^L\otimes V^R})$ 
gives a module structure on $F'$. 
By Proposition \ref{prop-tilde-A}, 
$\tilde{A}_0\otimes \tilde{A}_{-1}(\mathbb{Y}_f)$ is an 
intertwining operator of type $\binom{F'}{FF'}$. 

Let $W_i, i=1,2,3,4,5$ be $V^L\otimes V^R$-modules and 
$\Y$ be an intertwining operator of type $\binom{W_3}{W_1W_2}$. 
Let $f: W_2\rightarrow W_4$ and $g: W_3\rightarrow W_5$ are two 
$V^L\otimes V^R$-module maps. If we define
\beq
\big( g \circ \Y \circ (I_{W_1} \otimes f) \big)(u_1, x)u_2 := 
g(\Y(u_1, x)f(u_2))
\eeq
for $u_1\in W_1, u_2\in W_2$, then it is clear that 
$g \circ \Y \circ (I_{W_1} \otimes f)$ is an intertwining 
operator of type $\binom{W_5}{W_1W_4}$.

\begin{lemma} 
For $(F,\mathbb{Y}, \rho)$ to be
equipped with an invariant bilinear form $(\cdot, \cdot)$ 
is equivalent to give a homomorphism $\varphi: F\rightarrow F'$ 
as a $V^L\otimes V^R$-module map such
that the following condition: 
\beq  \label{A-r--r}
\mathbb{Y}_f = \varphi^{-1} \circ \tilde{A}_0 
\otimes \tilde{A}_{-1}
(\mathbb{Y}_f) \circ (I_{F} \otimes \varphi))
\eeq
is satisfied.  
\end{lemma}
\pf
The invariant bilinear form $(\cdot, \cdot)$
and $\varphi$ are related as follow:
\beq
\langle \varphi(v), u\rangle = (v, u) .
\eeq
One can use above equation to obtain one from the other. 
Moreover, the invariant property of bilinear form 
given in (\ref{inv-form-ffa-1}) exactly amounts to $\varphi$ 
being a $V^L\otimes V^R$-module map and satisfying 
the condition (\ref{A-r--r}). 
\epf

\begin{thm}  \label{F-alg-thm}
The following two notions are
equivalent in the sense that corresponding categories are 
isomorphic.  
\begin{enumerate}
\item A conformal full field algebra $F$ over $V^L\otimes V^R$ 
with a nondegenerate invariant bilinear form.
\item A commutative Frobenius algebra $F$ in 
$(\mathcal{C}_{V^L\otimes V^R}, \mathcal{R}_{+-})$ with $\theta_{F}=I_{F}$.  
\end{enumerate}
\end{thm}
\pf
By Theorem \ref{cfa-thm}, 
a full field algebra $F$ over $V^L\otimes V^R$ gives 
an commutative associative algebra with $\theta_F=I_F$. 
An invariant bilinear form on $F$ amounts to give a morphism 
from $\varphi: F \rightarrow F'$ so that (\ref{A-r--r}) holds.

$F'$ have a natural coalgebra structure with 
comultiplication $\Delta_{F'}$ and 
counit $\epsilon_{F'}$ defined as 
\beq  \label{co-F-dual}
\begin{picture}(14,2)
\put(2,1){$\Delta_{F'}$}
\put(3,1){$=$}
\put(4,0){\resizebox{3cm}{2cm}{\includegraphics{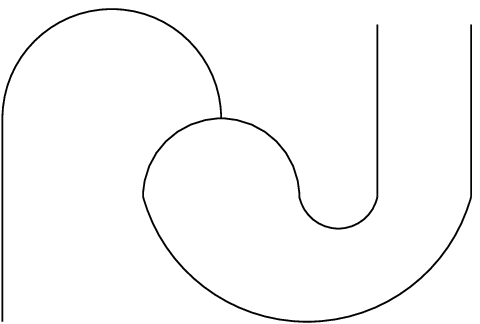}}}

\put(8,0){,}

\put(9,1){$\epsilon_{F'}$}
\put(10,1){$=$}
\put(11,0.5){\resizebox{1cm}{1cm}{\includegraphics{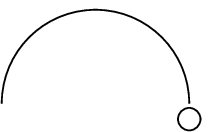}}}
\put(13,0){.}
\end{picture}
\eeq

When the bilinear form is nondegenerate, $\varphi$ is invertible. 
We will use the following graphic notation for $\varphi$ and 
$\varphi^{-1}$: 
\beq  \label{varphi}
\begin{picture}(14,2)
\put(2,1){$\varphi$} 
\put(3,1){$=$}
\put(4,0){\resizebox{0.8cm}{2cm}{\includegraphics{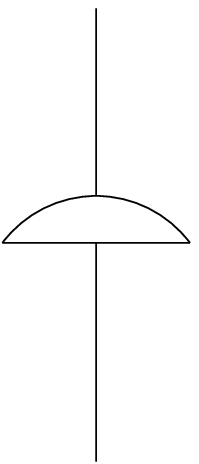}}}
\put(4.5,0){$F$}\put(4.5,1.8){$F'$}

\put(6,0){,}

\put (7,1){$\varphi^{-1}$} 
\put(8,1){$=$}
\put(9,0){\resizebox{0.8cm}{2cm}
{\includegraphics{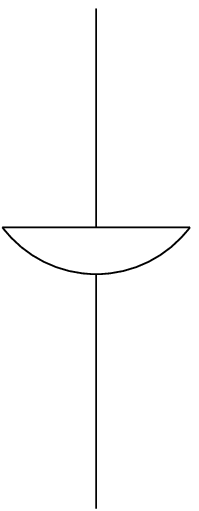}}}
\put(9.5,0){$F'$}\put(9.5,1.8){$F$}

\put(11,0){.}

\end{picture}
\eeq
Using the map $\varphi$ and its inverse, we can obtain
a natural coalgebra structure on $F$ defined as follow: 
\beq  \label{co-F}
\begin{picture}(14,2)
\put(2,0.7){$\Delta_{F}$}
\put(3,0.7){$=$}
\put(4,0){\resizebox{3cm}{2cm}{\includegraphics{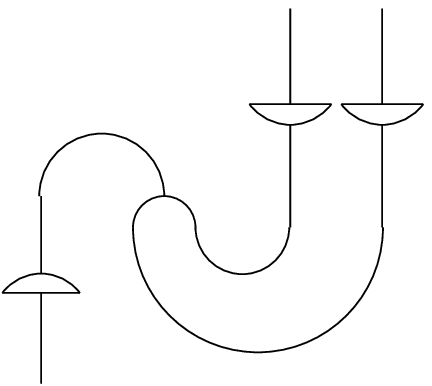}}}
\put(3.8,0){$F$}\put(5.7,1.7){$F$}\put(6.8,1.7){$F$}
\put(9,0.7){$\epsilon_{F}$}
\put(10,0.7){$=$}
\put(11,0){\resizebox{1cm}{1.5cm}{\includegraphics{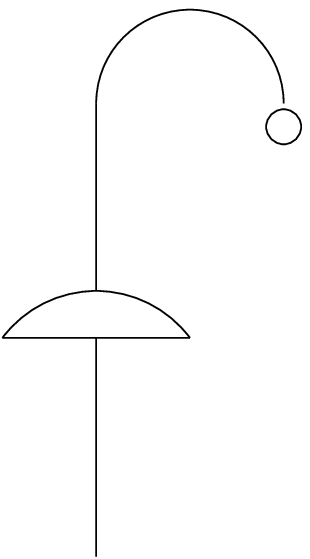}}}
\put(10.8,0){$F$}\put(13,0){.}
\end{picture}
\eeq

We claim that the multiple $(F, m, \iota, \Delta, \epsilon)$ gives 
a Frobenius algebra in $(\mathcal{C}, \mathcal{R}_{+-})$. 
We will prove the defining property (\ref{Frob-alg-def-fig}) of 
Frobenius algebra below. First, we define 
a left action of $F$ on $F'$ as follow: 
\beq  \label{F-l-act-F-dual}
\begin{picture}(14,2)
\put(3,0){\resizebox{1.5cm}{2cm}
{\includegraphics{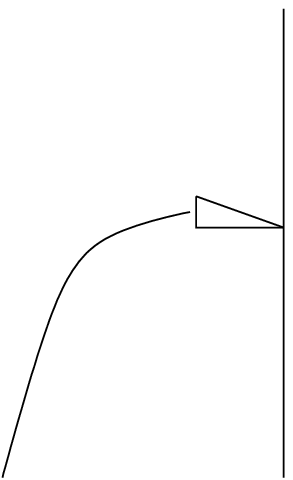}}}
\put(2.6,0){$F$}\put(4.6,0){$F'$}
\put(6,1){$:=$}
\put(7.5,0){\resizebox{2.5cm}{2cm}
{\includegraphics{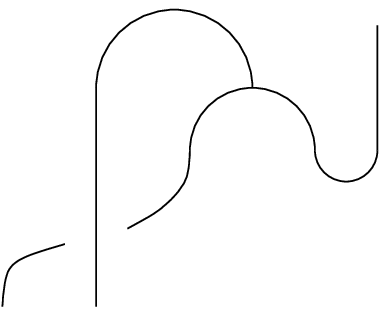}}}
\put(7.1,0){$F$}\put(8.2,0){$F'$}\put(10.1,1.8){$F'$}
\put(11,0){.}
\end{picture}
\eeq
By (\ref{tilde-A-0-graph}),
(\ref{tilde-A--1-graph}) and our choice of
braiding $\mathcal{R}_{+-}$, 
the right hands side of (\ref{F-l-act-F-dual}) is nothing but
$\tilde{A}_0\otimes \tilde{A}_{-1} (m)$. 
Then it not hard to see that (\ref{A-r--r}) exactly amounts to 
the following relation:
\beq  \label{F-act-F-varphi}
\begin{picture}(14,2)
\put(3,0){\resizebox{2cm}{2cm}{\includegraphics{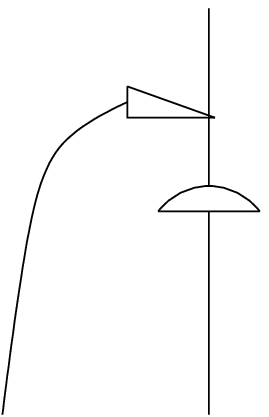}}}
\put(2.7,0){$F$}\put(4.8,0){$F$}\put(4.8,1.7){$F'$}
\put(6,1){$=$}
\put(7.5,0){\resizebox{1.5cm}{2cm}{\includegraphics{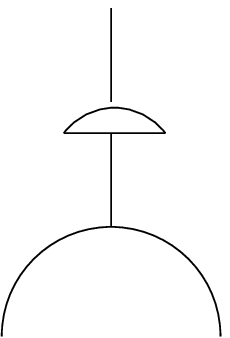}}}
\put(7.1,0){$F$}\put(9,0){$F$}\put(8.4,1.7){$F'$}
\put(9.5,0){.}
\end{picture}
\eeq
We also introduce a right action of $F$ on $F'$ as follow:
\beq  \label{F-r-act-F-dual}
\begin{picture}(14,2)
\put(3,0){\resizebox{1.5cm}{2cm}
{\includegraphics{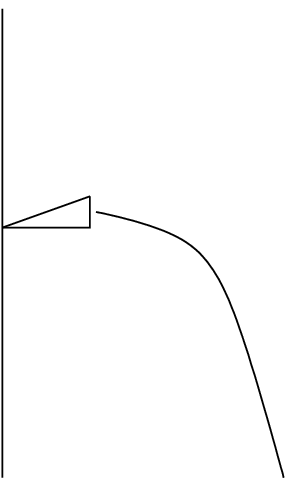}}}
\put(2.5,0){$F'$}\put(4.6,0){$F$}
\put(6,1){$:=$}
\put(7.5,0){\resizebox{2.5cm}{2cm}
{\includegraphics{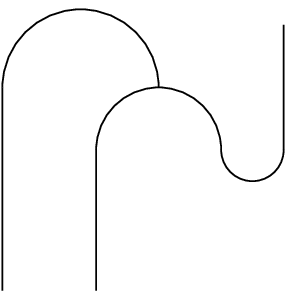}}}
\put(7,0){$F'$}\put(8.4,0){$F$}\put(10.1,1.8){$F'$}
\end{picture}
\eeq
Using (\ref{F-act-F-varphi}) and the commutativity of $F$, it 
is easy to obtain that 
\beq  \label{F-r-act-F-varphi}
\begin{picture}(14,2)
\put(3,0){\resizebox{2cm}{2cm}{\includegraphics{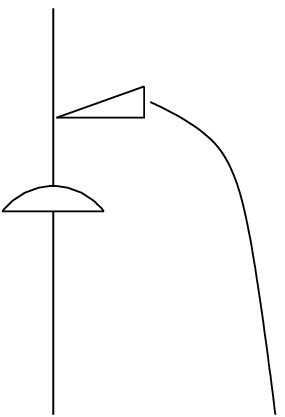}}}
\put(2.9,0){$F$}\put(5.1,0){$F$}\put(3.5,1.7){$F'$}
\put(6,1){$=$}
\put(7.5,0){\resizebox{1.5cm}{2cm}
{\includegraphics{varphi-2.eps}}}
\put(7.1,0){$F$}\put(9.1,0){$F$}\put(8.4,1.7){$F'$}\put(10,0){.}
\end{picture}
\eeq
Using (\ref{F-r-act-F-varphi}), we see that the 
$\Delta$ can be rewritten as follow: 
\beq \label{Delta-2-equ}
\begin{picture}(14,3)

\put(0.5, 0.5){\resizebox{1cm}{1.2cm}{\includegraphics{dual-Y.eps}}}
\put(2,1){$:=$}
\put(2.8,0){\resizebox{3cm}{3cm}{\includegraphics{Delta.eps}}}
\put(6.3,1){$=$}
\put(7.2,0){\resizebox{3cm}{3cm}{\includegraphics{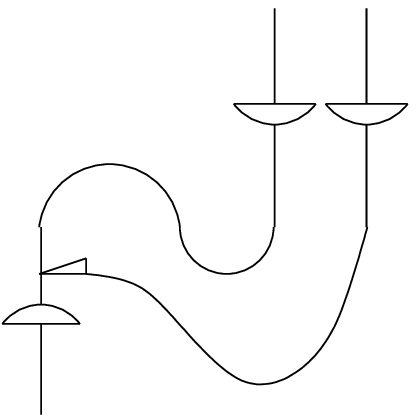}}}
\put(10.5,1){$=$}
\put(11.5,0){\resizebox{2cm}{2cm}{\includegraphics{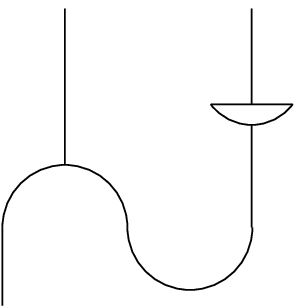}}}
\put(13.8,0){.}
\end{picture}
\eeq
By (\ref{Delta-2-equ}) and the associativity of $F$,  we have
$$
\begin{picture}(14,2)
\put(0.5,0){\resizebox{2cm}{2cm}{\includegraphics{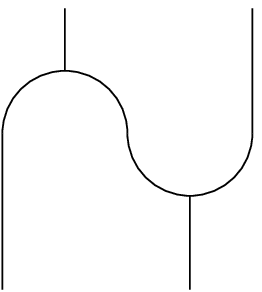}}}
\put(3.5,1){$=$}
\put(4.5,0){\resizebox{2.5cm}{2cm}{\includegraphics{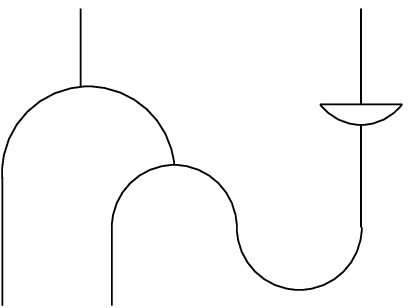}}}
\put(7.5,1){$=$}
\put(8.5,0){\resizebox{2.5cm}{2cm}{\includegraphics{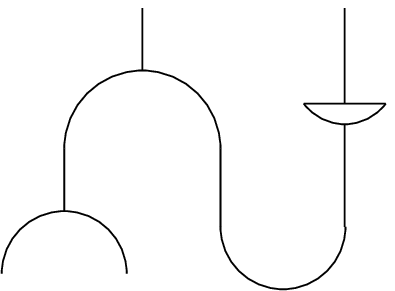}}}
\put(11.5,1){$=$}
\put(12.5,0){\resizebox{1cm}{2cm}{\includegraphics{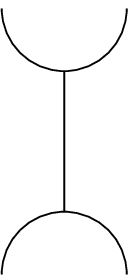}}}
\put(13.8,0){,}
\end{picture}
$$
and 
$$
\begin{picture}(14,2)
\put(0.5,0){\resizebox{2cm}{2cm}{\includegraphics{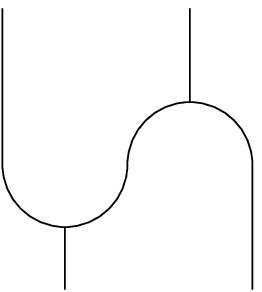}}}
\put(3.5,1){$=$}
\put(4.5,0){\resizebox{2cm}{2cm}{\includegraphics{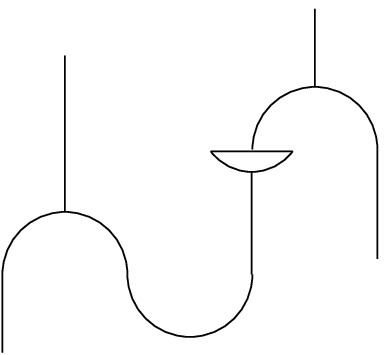}}}
\put(7.5,1){$=$}
\put(8.5,0){\resizebox{2cm}{2cm}{\includegraphics{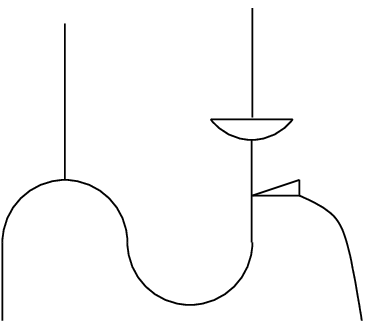}}}

\end{picture}
$$
$$
\begin{picture}(14,2)
\put(3.5,1){$=$}
\put(4.5,0){\resizebox{2cm}{2cm}{\includegraphics{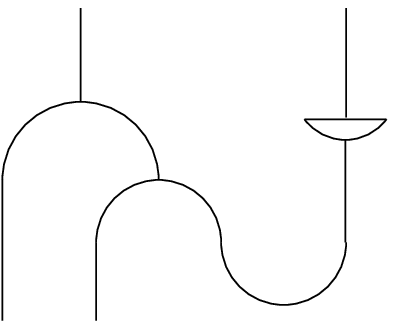}}}
\put(7.5,1){$=$}
\put(8.5,0){\resizebox{2.5cm}{2cm}{\includegraphics{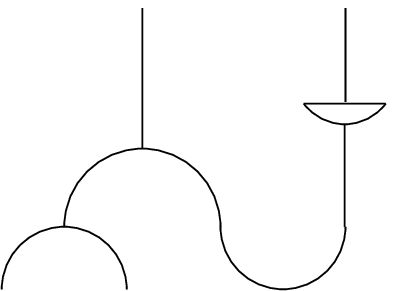}}}
\put(11.5,1){$=$}
\put(12.5,0){\resizebox{1cm}{2cm}{\includegraphics{Frob-M.eps}}}
\put(13.8,0){.}
\end{picture}
$$
Hence we have proved that $F$ is commutative Frobenius algebra
with $\theta_F=I_F$ in $(\mathcal{C}_{V^L\otimes V^R}, \mathcal{R}_{+-})$.

Conversely, given a commutative 
Frobenius algebra $F$ with $\theta_F=I_F$ in 
$\mathcal{C}_{V^L\otimes V^R}$, it is shown in \cite{FRS2}
that there is an isomorphism $\Phi$ from $F$ to $F'$. 
$\Phi$ and its inverse $\Phi^{-1}$ is given in the following diagrams:
\beq  \label{Phi-1}
\begin{picture}(14,2)
\put(2,0.7){$\Phi$} \put(3,0.7){$=$}
\put(4,0){\resizebox{2cm}{2cm}{\includegraphics{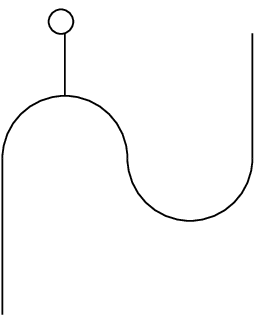}}}

\put(6.5,0){,}

\put(7,0.7){$\Phi^{-1}$}\put(8,0.7){$=$}
\put(9,0){\resizebox{2cm}{2cm}{\includegraphics{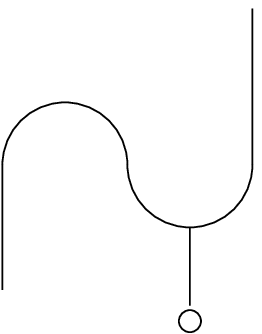}}}
\put(12,0){.}
\end{picture}
\eeq
$\Phi$ induces a nondegenerate bilinear from on $F$. 
We want to show that this bilinear form is also invariant. 
Let $\varphi:=\Phi$. Adopting the notation of $\varphi$ in 
(\ref{varphi}), we see that to prove 
the invariance property of $\Phi$ (or $\varphi$) 
amounts to prove that (\ref{F-act-F-varphi}) holds. 
We first prove (\ref{F-r-act-F-varphi}).  Indeed, we have
$$
\begin{picture}(14,2)
\put(1,0){\resizebox{2.5cm}{2cm}{\includegraphics{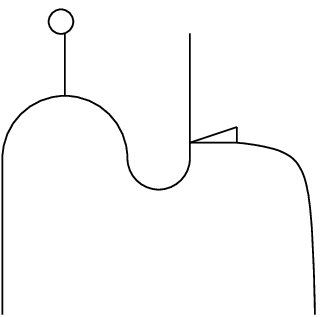}}} 
\put(4,1){$=$}
\put(5,0){\resizebox{2.5cm}{2cm}{\includegraphics{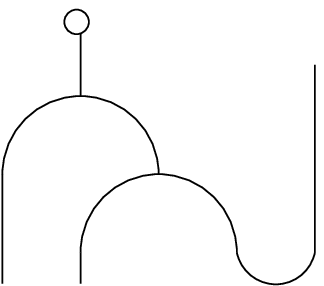}}}
\put(8.5,1){$=$}
\put(9.5,0){\resizebox{2.5cm}{2cm}{\includegraphics{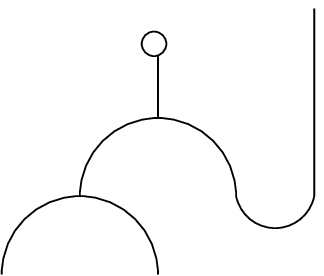}}}
\put(13,0){,}
\end{picture}
$$
which is nothing but (\ref{F-r-act-F-varphi}). Then 
(\ref{F-act-F-varphi}) follows from (\ref{F-r-act-F-varphi}) 
by the commutativity of $F$.

So far, we have constructed two functors between the two categories.
It remains to show that these two functors 
give isomorphisms between two categories.

One can check that if a Frobenius algebra 
$F$ is obtained from an associative algebra $F$ with an 
isomorphism $\varphi: F \rightarrow F'$ using (\ref{co-F}),
then $\Phi=\varphi$. Indeed, we have
\begin{center}
\epsfxsize 0.9\textwidth
\epsfysize  0.16\textwidth
\epsfbox{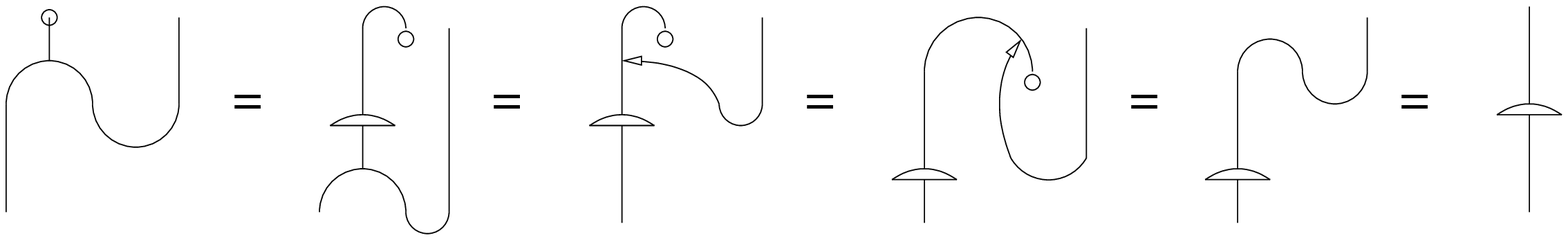}\hspace{0.3cm}.
\end{center}

Similarly, given a Frobenius algebra $F$. One have 
the isomorphism $\Phi: F \rightarrow F'$ given 
as (\ref{Phi-1}). 
We claim that the Frobenius algebra structure on $F$ induced
from $\Phi$ using the construction (\ref{co-F}) 
is exactly same as the original one. The proof is the
given in the following picture. 
$$
\begin{picture}(14,2)
\put(2.4,1){$\scriptstyle \Phi$}\put(4,1){$=$}\put(8.1,1){$=$}
\put(10.5,1){$=$}
\put(2,0){\resizebox{10cm}{2cm}{\includegraphics{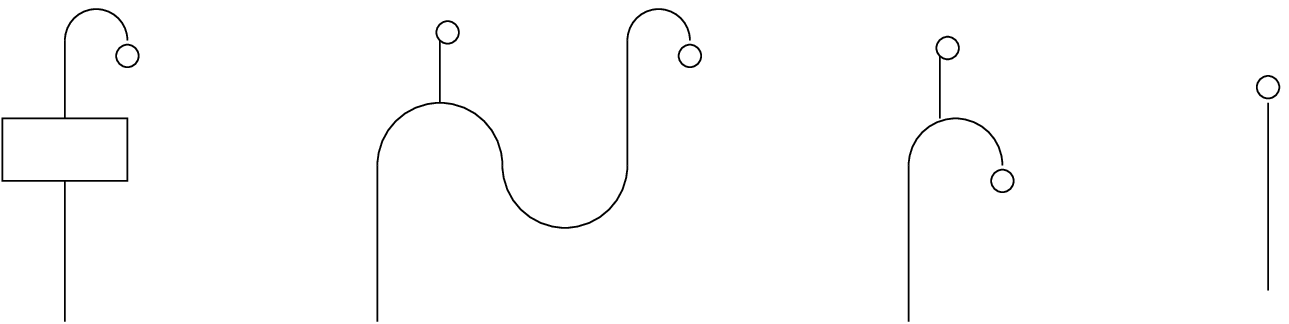}}} 
\end{picture}
$$
$$
\begin{picture}(14,2)
\put(2.4,1.1){$\scriptstyle \Phi^{-1}$}
\put(3.5,1){$=$}\put(7,1){$=$}\put(9.7,1){$=$}
\put(11.6,1){$=$}
\put(1,0){\resizebox{12cm}{2cm}{\includegraphics{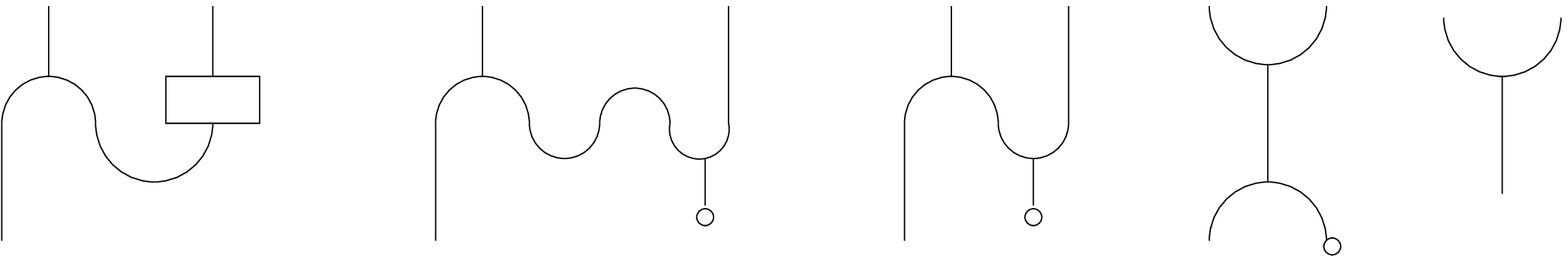}}} 
\put(13.5,0){.}
\end{picture}
$$
\epf

\begin{rema}
{\rm So far we have completely reformulated the genus-zero conformal
field theory in term of algebras in 
tensor category. At least, at the level 
of tensor category, the conformal field theories are the essentially
same as that of closed topological field theories. They
are both commutative Frobenius algebras, however, in 
very different categories. We will show in \cite{HK3} 
that genus-one conformal field theories are very different from
genus-one topological field theories 
even on the level of tensor category. }
\end{rema}

\section{Construction}

In this section, we consider the case $V^L=V^R=V$, where
$V$ is assumed to satisfy the conditions in Theorem \ref{MTC}. 
We will give a construction of 
commutative Frobenius algebra $F$ in 
the category $(\mathcal{C}_{V\otimes V}, \mathcal{R}_{+-})$
satisfying $\theta_F=I_F$. 
We assume all the notations used in Section 4.

Let $F$ is be the object in $\mathcal{C}_{V\otimes V}$ 
given as follow
\beq  \label{cvoa-const-1}
F = \coprod_{a\in \mathcal{I}} W^a\otimes (W^a)'.
\eeq
The decomposition of $F$ as a direct sum gives 
a natural embedding $V\otimes V \hookrightarrow F$. 
We denote this embedding as $\iota_{F}$. 
Now we need select a single morphism  
$m_{F} \in \hom_{V\otimes V}(F \boxtimes F, F)$.

By the universal property of tensor product, 
$\{ \Y_{a_1a_2; i}^{a_3;(1)} \}_{i=1}^{N_{a_1a_2}^{a_3}}$, 
a basis of $\V_{a_1a_2}^{a_3}$,
gives arise to $\{ e_{a_1a_2; i}^{a_3} \}_{i=1}^{N_{a_1a_2}^{a_3}}$,
a basis of
$\hom_V(W^{a_1}\boxtimes W^{a_2}, W^{a_3})$, as given in the graph 
(\ref{basis-pic}). And the dual basis 
$\{ f_{a_3; j}^{a_1a_2} \}_{j=1}^{N_{a_1a_2}^{a_3}}$ is given in 
(\ref{dual-basis-pic}), satisfying the duality 
condition (\ref{Y-dual-Y}).

In general, $m_{F}$ can always be written in the 
following form 
\beq  \label{cvoa-const-2}
m_{F} = \sum_{a_1,a_2,a_3\in \mathcal{A}} \sum_{i,j=1}^{N_{a_1a_2}^{a_3}} 
\langle f^{a_1a_2}_{a_3;i}, f^{a'_1a'_2}_{a'_3;j}\rangle 
e_{a_1a_2;i}^{a_3} \otimes e_{a'_1a'_2; j}^{a'_3},
\eeq
where $\langle \cdot, \cdot\rangle$ is a 
bilinear pairing between $\hom_V(a_1\boxtimes a_2, a_3)$ and
$\hom_V(a'_1\boxtimes a'_2, a'_3)$. Actually, using the duality
(\ref{Y-dual-Y}), $m_F$ can be
identified with the pairing $\langle \cdot, \cdot \rangle$ 
as an element in 
$$
(\hom_V(a_3, a_1\boxtimes a_2))^* \otimes 
(\hom_V(a'_3, a'_1\boxtimes a'_2))^*.
$$
In particular, it also means that $m_F$ defined in 
(\ref{cvoa-const-2}) is 
independent of the choice of basis. 
We define the pairing 
$\langle f^{a_1a_2}_{a_3;i}, f^{a'_1a'_2}_{a'_3;j}\rangle$ by
\beq  \label{pairing-equ}
\begin{picture}(14, 4)
\put(4,0){\resizebox{5cm}{3.5cm}{
\includegraphics{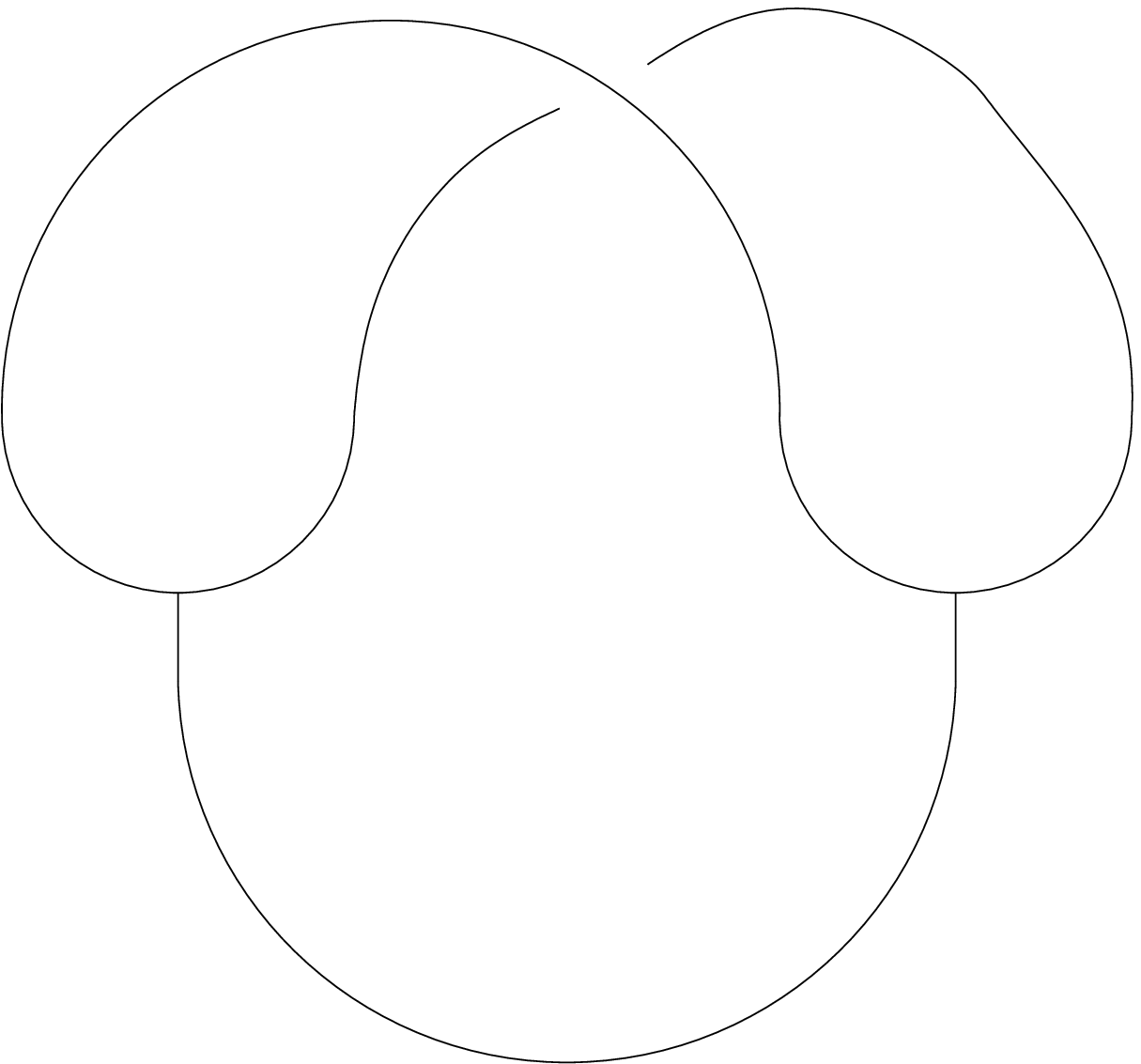}}}
\put(2, 1.7){$\ds \frac{1}{\dim a_3}$}
\put(3.6, 2.2){$a_1$}
\put(4.8, 1.7){$i$}\put(4.4, 1.2){$a_3$}
\put(5.8, 2.2){$a_2$}
\put(7, 2.2){$a'_1$}\put(8.1, 1.7){$j$}
\put(9.1, 2.2){$a'_2$}\put(8.3, 1.2){$a'_3$}
\put(10.5,0){.}
\end{picture}
\eeq

The following identity is manifest by itself, 
and will be useful later. 
\beq  \label{pair-dual-fig}
\begin{picture}(14,3)

\put(1.2,0){$a_3$}\put(5.1,2.8){$a_3$} 
\put(0.6,1.6){$a_1$}\put(2.4,1.6){$a_2$}
\put(1.6,1.3){$i$}\put(3.6,1.4){$j$}

\put(1,0){\resizebox{4cm}{3cm}{
\includegraphics{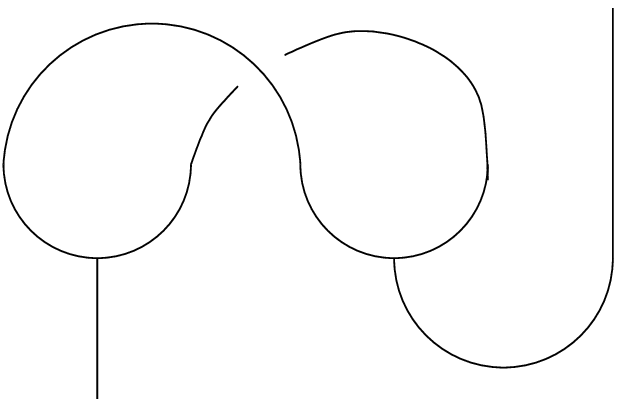}}}

\put(5.5,1.5){$\ds = \,  \frac{1}{\dim a_3} $}

\put(7.6, 1.6){$a_1$}\put(9.4,1.6){$a_2$}
\put(8.2, 0.8){$a_3$}\put(8.6, 1.5){$i$}\put(11.3,1.6){$j$}

\put(8, 0){\resizebox{4cm}{3cm}{
\includegraphics{pair-0.eps}}}

\put(12.9, 0){$a_3$}\put(12.9, 2.8){$a_3$}
\put(12.7, 0){\resizebox{0.1cm}{3cm}{
\includegraphics{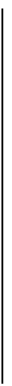}}}
\end{picture}
\eeq

The following Lemma is also manifest by itself. 
It has appeared in \cite{FFFS}.
\begin{lemma}
\beq  \label{F--1-fig}
\begin{picture}(14,3)

\put(0.5, 1.5)
{$F^{-1}(\Y_{a_6a_3;k}^{a_4;(1)} \otimes \Y_{a_1a_2;l}^{a_6;(1)},
\Y_{a_1a_5;i}^{a_4;(1)}\otimes \Y_{a_2a_3;j}^{a_5;(1)}) $}
\put(8.1, 0){\resizebox{0.1cm}{3cm}{
\includegraphics{id.eps}}}
\put(8.3,0){$a_4$}\put(8.3,2.8){$a_4$}

\put(9.3, 1.5){$=$}

\put(10.5,0){\resizebox{2.5cm}{3cm}{
\includegraphics{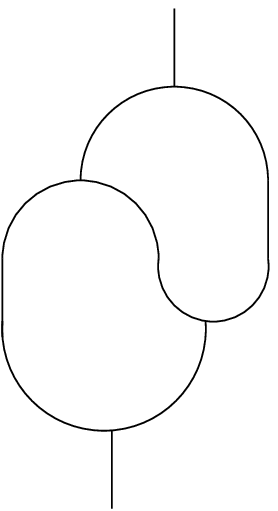}}}

\put(12.3,2.8){$a_4$}\put(11.8,0){$a_4$}
\put(11.4,0.7){$i$}\put(12.3,1.4){$j$}
\put(11.2,1.6){$l$}\put(12, 2.1){$k$}

\put(10.1,1.3){$a_1$}\put(11.5, 1.4){$a_2$}\put(13.1,1.6){$a_3$}
\put(12.3, 0.6){$a_5$}\put(11, 2.3){$a_6$}

\end{picture}
\eeq
\end{lemma}

Notice that $F'$ has the same contents as $F$. They
are isomorphic as $V\otimes V$-modules.  
There is, however, no canonical isomorphism. 
The space of isomorphisms between $F$ and $F'$ is a 
three dimensional vector space. 
Now we choose a particular isomorphism  
$\varphi_F: F\rightarrow F'$ given by 
\beq  \label{varphi-F}
\varphi_F=\oplus_{a\in \I} \,\,
\frac{1}{\dim a} \, e^{-2\pi ih_{a}} I_{W^{a}\otimes (W^a)'}.
\eeq
The isomorphism $\varphi_F$ induces a 
nondegenerate invariant bilinear form on $F$ viewed as
$V\otimes V$-module.

\begin{thm}
The triple $(F, m_{F}, \iota_{F})$, together with
the bilinear form induced from $\varphi_F$, is a 
commutative Frobenius algebra in 
$\mathcal{C}_{V\otimes V}$ satisfying $\theta_{F}= I_{F}$. 
\end{thm}
\pf
The left and right unit properties follow from the facts that 
$\langle l_{W^{a_1}}^{-1},  l_{(W^{a_1})'}^{-1} \rangle=1$ and 
$\langle r_{W^{a_1}}^{-1},  r_{(W^{a_1})'}^{-1} \rangle=1$. 
The condition $\theta_{F}=I_{F}$ is obvious by the definition
of $F$ in (\ref{cvoa-const-1}).

The commutativity amounts to 
the invariant property of the pairing (\ref{pairing-equ}) 
with respect to braiding isomorphisms. This invariant property 
of the pairing is rather easy to see, 
and was pointed out by Kirillov in (\cite{Kr}). 

Now we prove the associativity:
\beq \label{asso-F}
m_F \circ ( I_F \boxtimes m_F) = m_F \circ ( m_F \boxtimes I_F) \circ 
\mathcal{A}.
\eeq
Using (\ref{pair-dual-fig}), we have 
$$
\begin{picture}(14,2)
\put(0.5,0.8){$\ds m_F \circ ( I_F \boxtimes m_F) =
\sum_{a_1,a_2,a_3,a_4}\sum_{a_5}\sum_{i,j,k,l}$}

\put(7.1,0){$a_1$}\put(7.8,0){$a_2$}\put(9.6,0){$a_3$}
\put(7.8,1.8){$a_4$}\put(8.2, 1.1){$i$}\put(8.8,0.3){$j$}
\put(8.9, 1.1){$a_5$}

\put(7.5,0){\resizebox{2cm}{2cm}{
\includegraphics{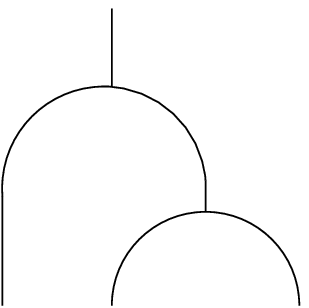}}}

\put(10, 0.8){$\otimes$}

\put(10.6,0){$a'_1$}\put(11.3,0){$a'_2$}\put(13.1,0){$a'_3$}
\put(11.3,1.8){$a'_4$}\put(11.7, 1.1){$k$}\put(12.3,0.3){$l$}
\put(12.4, 1.1){$a'_5$}

\put(11, 0){\resizebox{2cm}{2cm}{
\includegraphics{asso-left.eps}}}

\end{picture}
$$
\beq
\begin{picture}(14,3)

\put(4.6, 1.2){$\ds \cdot \, \, \frac{1}{\dim a_5}$}

\put(6.1,1.7){$a_2$}\put(7.5,1.7){$a_3$}\put(6.8, 1.4){$j$}
\put(6.6, 0.3){$a_5$} \put(8.5,1.4){$l$}

\put(6.5,0){\resizebox{2.5cm}{2.5cm}{
\includegraphics{pair-0.eps}}}

\put(9.5, 1.2){$\ds \frac{1}{\dim a_4}$}

\put(10.6,1.7){$a_1$}\put(12,1.7){$a_5$}\put(11.3, 1.4){$i$}
\put(11.1, 0.3){$a_4$} \put(13,1.4){$k$}

\put(11, 0){\resizebox{2.5cm}{2.5cm}{
\includegraphics{pair-0.eps}}}

\put(13.8,0){.}
\end{picture}
\eeq
$$
\begin{picture}(14,2)
\put(3.5,0.8){$\ds  =
\sum_{a_1,a_2,a_3,a_4}\sum_{a_5,a_6}\sum_{i,j,k,l}$}

\put(7.1,0){$a_1$}\put(7.8,0){$a_2$}\put(9.6,0){$a_3$}
\put(7.8,1.8){$a_4$}\put(8.2, 1.1){$i$}\put(8.8,0.3){$j$}
\put(8.9, 1.1){$a_5$}

\put(7.5,0){\resizebox{2cm}{2cm}{
\includegraphics{asso-left.eps}}}

\put(10, 0.8){$\otimes$}

\put(10.6,0){$a'_1$}\put(11.3,0){$a'_2$}\put(13.1,0){$a'_3$}
\put(11.3,1.8){$a'_4$}\put(11.7, 1.1){$k$}\put(12.3,0.3){$l$}
\put(12.4, 1.1){$a'_6$}

\put(11, 0){\resizebox{2cm}{2cm}{
\includegraphics{asso-left.eps}}}

\end{picture}
$$
\beq \label{const-fig-1-equ}
\begin{picture}(14,3.5)

\put(4.6, 1.7){$\ds \cdot \, \, \frac{1}{\dim a_4}$}

\put(6.6,1.2){$a_1$}\put(8,1.2){$a_5$}\put(7.4,1.3){$i$}
\put(7.2, 1.9){$a_2$}\put(9.4, 1.9){$a'_3$}\put(8.9, 1.8){$l$}
\put(7.9, 1.8){$j$}\put(10.6,1.6){$k$}\put(11.1, 1.8){$a'_6$}
\put(7.5, 0.3){$a_4$}

\put(7,0){\resizebox{4cm}{3cm}{
\includegraphics{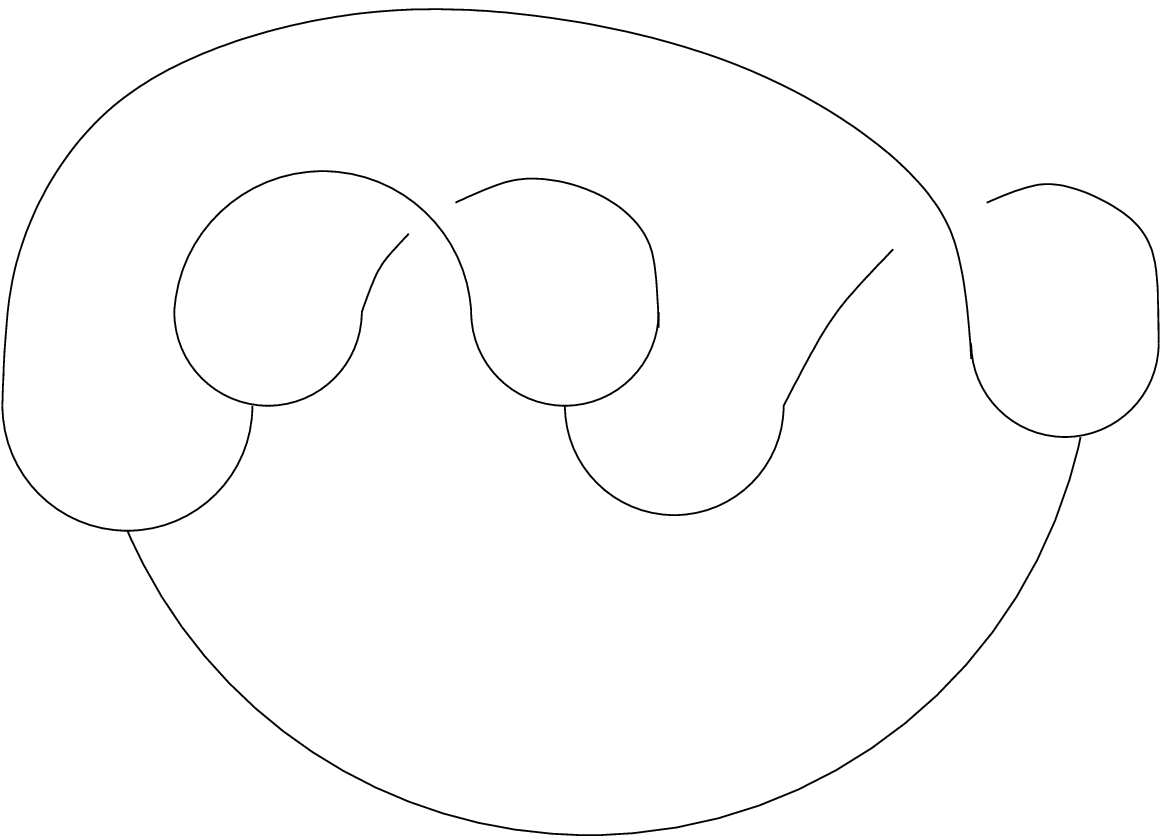}}}

\put(12,0){,}
\end{picture}
\eeq
where we have added a sum $\sum_{a_6}$ in the second step. 
It is harmless because the graph in the last 
line of (\ref{const-fig-1-equ}) is nonzero if and only if $a_6=a_5$. 

Now we can first deform the graph in 
the last line of (\ref{const-fig-1-equ}), 
then apply (\ref{exp-id}) to obtain the following identity: 
\beq  \label{const-asso-pair-equ}
\begin{picture}(14, 4.5)

\put(0.7, 1.9){$a_1$}\put(2.4, 1.9){$a_5$}\put(2.3, 2.5){$j$}
\put(1.7, 1.7){$i$}\put(1.4, 0.8){$a_4$}\put(1.5, 2.7){$a_2$}
\put(3, 2.7){$a_3$}
\put(3.4, 1.9){$a'_1$}\put(4.9, 2.5){$l$}
\put(4.3, 1.7){$k$}\put(5, 1.9){$a'_6$}

\put(1,0){\resizebox{4.5cm}{4.5cm}{
\includegraphics{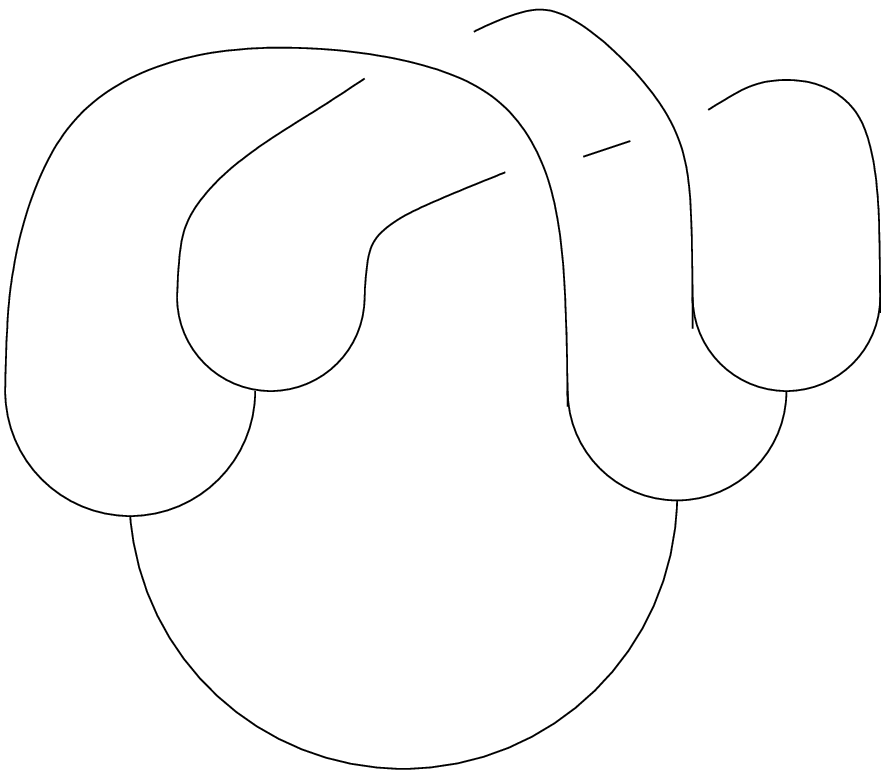}}}

\put(6, 2){$\ds =\sum_{a_7, a_8}\sum_{p,q,r,s}$}

\put(8.1, 1.2){$a_1$}\put(8.2, 3.7){$a_1$}
\put(9.6, 0.9){$a_5$}\put(9,1){$i$}\put(8.9, 1.4){$a_2$}
\put(8.8,1.7){$q$}\put(9.4, 2.2){$p$}\put(9.4,3.1){$p$}
\put(12.8, 4.2){$a'_2$}
\put(8.7, 2.4){$a_7$}\put(8.7, 2.9){$a_7$}\put(8.9, 3.5){$q$}
\put(10.1, 1.7){$a_3$}\put(10,3.2){$a_3$}
\put(9.5, 1.5){$j$}\put(9.3, 0.3){$a_4$}\put(9.5, 2.6){$a_4$}

\put(11.6, 1.2){$a'_1$}
\put(13, 1){$a'_6$}\put(12.5,1){$k$}\put(12.3, 1.5){$a'_2$}
\put(12.2,1.8){$s$}\put(12.8, 2.3){$r$}\put(12.8,3.2){$r$}
\put(12, 2.4){$a'_8$}\put(12.1, 2.9){$a'_8$}\put(12.3, 3.6){$s$}
\put(13.6, 1.7){$a'_3$}\put(13.4,3.2){$a'_3$}
\put(13, 1.6){$l$}\put(13,2.7){$a'_4$}

\put(8.5,0){\resizebox{5cm}{4.5cm}{
\includegraphics{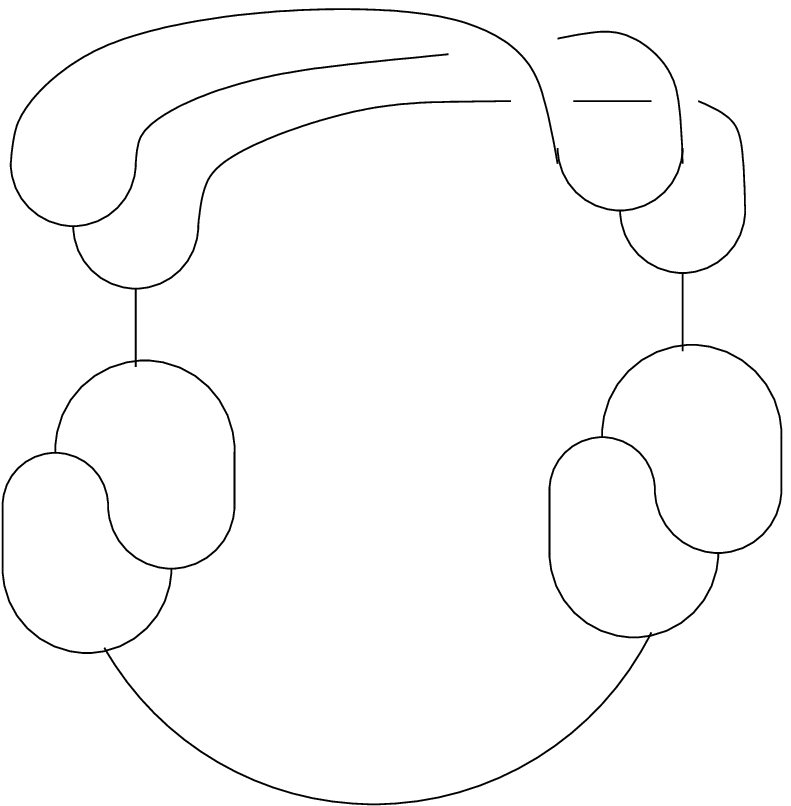}}}

\put(13.9,0){.}

\end{picture}
\eeq
Notice that when we apply (\ref{exp-id}) we need sum up the
intermediate index, such as the index $a_3$ in (\ref{exp-id}).
However, the second graph in (\ref{const-asso-pair-equ}) takes
nonzero value if and only if the intermediate indices are 
$a_4$ and $a'_4$ respectively.
  
Now using (\ref{F--1-fig}), we
can move the following two factors
\bea \label{two-F--1}
&F^{-1}(\Y_{a_7a_3;p}^{a_4;(1)} \otimes \Y_{a_1a_2;q}^{a_7;(1)}, 
\Y_{a_1a_5;i}^{a_4;(1)}\otimes \Y_{a_2a_3;j}^{a_5;(1)})&,  \nn
&F^{-1}(\Y_{a'_8a'_3;r}^{a'_4;(1)} 
\otimes \Y_{a'_1a'_2;s}^{a'_8;(1)}, 
\Y_{a'_1a'_6;k}^{a'_4;(1)} \otimes \Y_{a'_2a'_3;l}^{a'_6;(1)})&
\eea
out of the second graph in (\ref{const-asso-pair-equ}), 
and combine them with
$$
\begin{picture}(14,2)

\put(4.1,0){$a_1$}\put(4.8,0){$a_2$}\put(6.6,0){$a_3$}
\put(4.8,1.8){$a_4$}\put(5.2, 1.1){$i$}\put(5.8,0.3){$j$}
\put(5.9, 1.1){$a_5$}

\put(4.5,0){\resizebox{2cm}{2cm}{
\includegraphics{asso-left.eps}}}

\put(7, 0.8){$\otimes$}

\put(7.6,0){$a'_1$}\put(8.3,0){$a'_2$}\put(10.1,0){$a'_3$}
\put(8.3,1.8){$a'_4$}\put(8.7, 1.1){$k$}\put(9.3,0.3){$l$}
\put(9.4, 1.1){$a'_5$}

\put(8, 0){\resizebox{2cm}{2cm}{
\includegraphics{asso-left.eps}}}

\put(11,0){.}
\end{picture}
$$
Then we compose the resulting morphism from right by  
$\mathcal{A}^{-1}$, and sum up the indices $a_5, a_6, i,j,k,l$ 
(using (\ref{A-rel-F-2})).
As a result, we obtain the following identity:  
$$
\begin{picture}(14,2)
\put(0,0.8){$\ds m_F \circ ( I_F \boxtimes m_F) \circ \A^{-1} =
\sum_{a_1,a_2,a_3,a_4}\sum_{a_7,a_8}\sum_{p,q,r,s}$}

\put(7.6,0){$a_1$}\put(9.4,0){$a_2$}\put(10.1,0){$a_3$}
\put(9.4,1.8){$a_4$}\put(9.3, 1.1){$p$}\put(8.5,0.3){$q$}
\put(8.3, 1.1){$a_7$}

\put(8,0){\resizebox{2cm}{2cm}{
\includegraphics{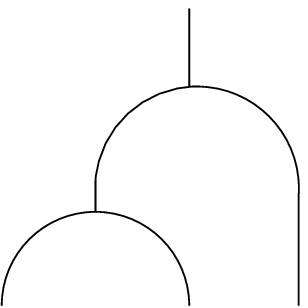}}}

\put(10.5, 0.8){$\otimes$}

\put(11.1,0){$a'_1$}\put(12.9,0){$a'_2$}\put(13.5,0){$a'_3$}
\put(12.9,1.8){$a'_4$}\put(12.8, 1.1){$r$}\put(12,0.3){$s$}
\put(11.8, 1.1){$a'_8$}

\put(11.5, 0){\resizebox{2cm}{2cm}{
\includegraphics{asso-right.eps}}}

\end{picture}
$$
\beq \label{asso-right-equ}
\begin{picture}(14,3.5)

\put(4.6, 1.7){$\ds \cdot \, \, \frac{1}{\dim a_4}$}

\put(7.1,1.3){$a_7$}\put(8.6,1.7){$a_3$}\put(7.9,1.4){$p$}
\put(6.6, 2.2){$a_1$}\put(7.9, 2.2){$a_2$}\put(9.6, 2){$s$}
\put(7.4, 2){$q$}\put(10,1.4){$r$}
\put(7.7, 0.3){$a_4$}\put(9.3, 1.4){$a'_8$}

\put(7,0){\resizebox{3.5cm}{3.5cm}{
\includegraphics{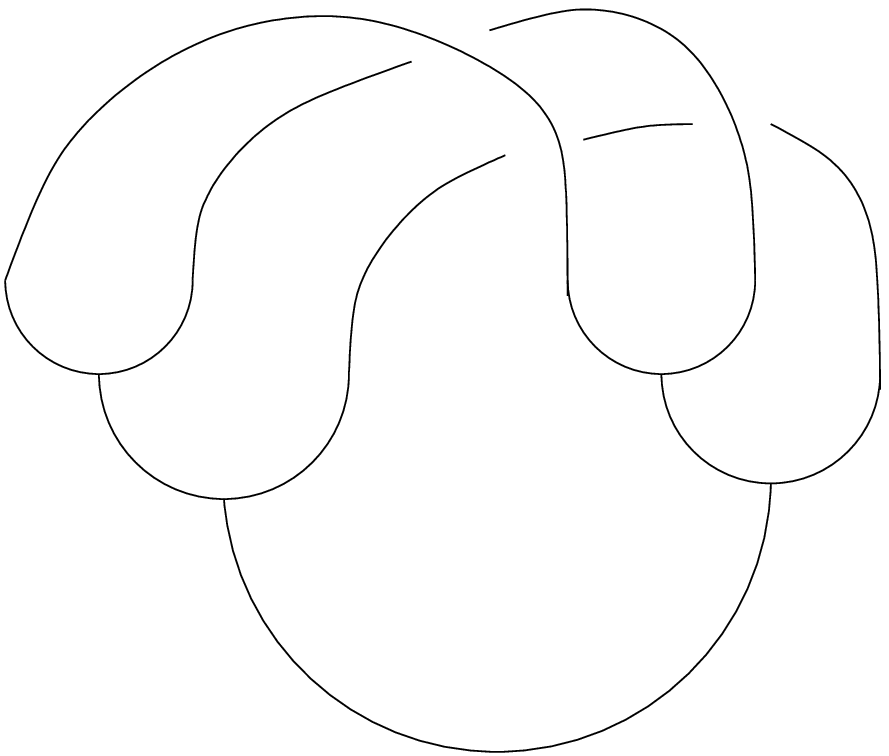}}}

\put(12,0){.}
\end{picture}
\eeq
Notice again that the graph in the second line of (\ref{asso-right-equ})
is nonzero if and only if $a_7=a_8$. 
It is easy to show that the right hand side of 
(\ref{asso-right-equ}) is nothing but
$m_F \circ ( m_F \boxtimes I_F)$.
Therefore, we obtain the associativity (\ref{asso-F}).

In summary, we have proved that  
$F$ is a commutative associative algebra with $\theta_F= I_F$. 
It remains to show that the bilinear form induced 
from $\varphi_F$ is invariant. This amounts to show 
that the following equation  
\beq  \label{m-F-inv-form}
m_F = \varphi_F^{-1} \circ \tilde{A}_0\otimes \tilde{A}_{-1} (m_F)
\circ (I_F \otimes \varphi_F)
\eeq
holds. Since the $m_F$ given in (\ref{cvoa-const-2}) 
is independent of choice of basis, we have 
\beq  \label{m-F-A-0-1}
m_F = \sum_{a_1, a_2, a_3\in \I} \sum_{i,j} 
\langle \hat{A}_0^* (f_{a_3;i}^{a_1a_2}), 
\hat{A}_{-1}^*(f_{a'_3;j}^{a'_1a'_2})\rangle 
\tilde{A}_0(e_{a_1a_2;i}^{a_3}) \otimes \tilde{A}_{-1}(e_{a'_1a'_2;j}^{a'_3}).
\eeq
Replacing the left hand side of (\ref{m-F-inv-form}) by
the right hand side of (\ref{m-F-A-0-1}) and the $m_F$ in the 
right hand side of (\ref{m-F-inv-form}) by the right hand side of 
(\ref{cvoa-const-2}), and using 
(\ref{A-dual-1}), (\ref{A-dual-2}) and (\ref{varphi-F}),  
we see that 
the validity of (\ref{m-F-inv-form}) follows from 
the following identity: 
$$
\begin{picture}(14,4.5)
\put(0.5,0)
{\resizebox{4.5cm}{4.5cm}{\includegraphics{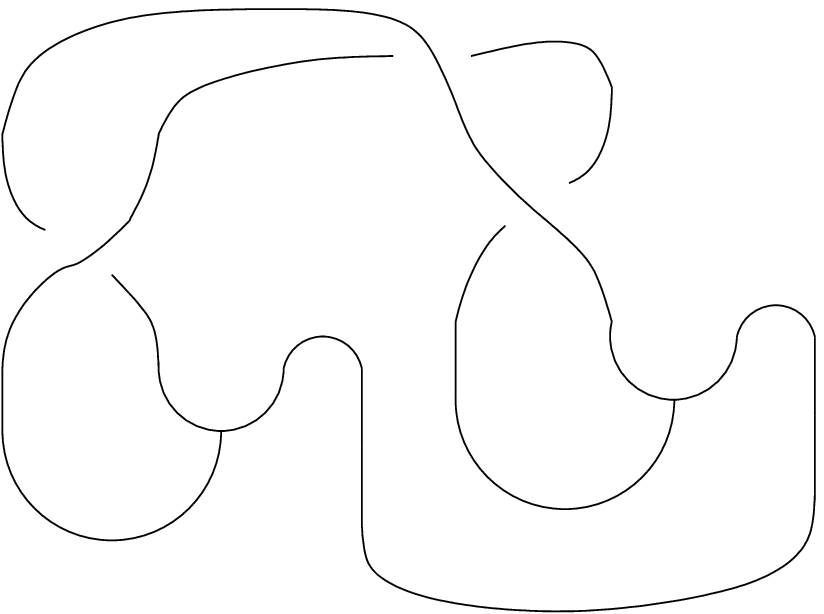}}}
\put(0, 3.2){$a_1$}
\put(1.6, 1.5){$i$}\put(1.7, 0.7){$a_3$}
\put(4.1, 1.8){$j$}
\put(5.1, 0.6){$a_2$}

\put(5.7,2){$= e^{2\pi i (h_{a_2}-h_{a_3})} $}

\put(9,0){\resizebox{4cm}{4cm}{\includegraphics{pair-0.eps}}}
\put(8.5, 2.2){$a_1$}
\put(9.5, 2){$i$}\put(9.2, 1.2){$a_3$}
\put(10.3, 2.2){$a_2$}
\put(11.3, 2.2){$a'_1$}\put(12.3, 2){$j$}
\put(13.2, 2.2){$a'_2$}\put(12.5, 1.2){$a'_3$}
\put(13.5,0){,}

\end{picture}
$$
which is manifestly true. 
Hence this bilinear form
induced from $\varphi_F$ is invariant. 
\epf

\appendix

\section{The Proof of rigidity}

In this section, we prove that
our new choice of duality maps, defined in Section 4.2, 
satisfies the rigidity axioms. 

\begin{lemma}  \label{3-equ-lemma}
For $a\in \I$, we have
\bea 
\Y_{a'a}^e &=& e^{2\pi i h_a} \Omega_{0} (\Y_{aa'}^e) =
 e^{-2\pi i h_a} \Omega_{-1} (\Y_{aa'}^e), \label{3-equ-1}  \\
\tilde{A}_0(\Y_{a'e}^{a'}) &=& e^{2\pi i h_a} \Y_{a'a}^{e'}, 
\label{3-equ-2}  \\
\tilde{A}_0(\Y_{ea'}^{a'}) &=& \Y_{ea}^a. \label{3-equ-3}
\eea
\end{lemma}
\pf
For $w_{a'}\in (W^{a})'$ and $w_a\in W^a$, we have
\bea  \label{3-equ-1-1}
(\one, \Y_{a'a}^e(w_{a'}, x) w_a) &=& \langle \one, 
\hat{A}_0(\Y_{a'e}^{a'})(w_{a'}, x) w_a\rangle  \nn
&=& \langle \Y_{a'e}^{a'}(e^{-xL(1)}x^{-2L(0)}w_{a'}, 
e^{\pi i}x^{-1})\one, w_a\rangle \nn
&=& \langle e^{-x^{-1}L(-1)}e^{-xL(1)}x^{-2L(0)} w_{a'}, w_a\rangle,
\eea
and
\bea  \label{3-equ-1-2}
(\one, \Omega_0(\Y_{aa'}^e)(w_{a'}, x) w_a) &=& 
\langle \one, e^{xL(-1)}\hat{A}_0(\Y_{ae}^a)(w_a, e^{\pi i}x)w_{a'}
\rangle \nn
&=& \langle  \Y_{ae}^a(e^{xL(1)}(e^{\pi i}x)^{-2L(0)} w_a, x^{-1})\one, 
w_{a'} \rangle   \nn
&=& \langle e^{x^{-1}L(-1)}e^{xL(1)}(e^{\pi i}x)^{-2L(0)} w_a, 
w_{a'} \rangle   \nn
&=& \langle x^{-2L(0)} e^{xL(-1)}e^{x^{-1}L(1)} e^{-2\pi iL(0)} w_a, w_{a'}
\rangle\nn
&=& \langle e^{-2\pi i L(0)}w_a, e^{-x^{-1}L(-1)}e^{-xL(1)}x^{-2L(0)} w_{a'}
\rangle.
\eea
Combining (\ref{3-equ-1-1}) and (\ref{3-equ-1-2}),
we obtain the first equality of (\ref{3-equ-1}). The second equality of
(\ref{3-equ-1}) can be proved similarly. 

Now we prove (\ref{3-equ-2}). On the one hand, we have
\bea  \label{3-equ-2-1}
\langle \tilde{A}_0(\Y_{a'e}^{a'})(w_{a'}, e^{\pi i} x) w_a, \one\rangle
&=& \langle w_a, \Y_{a'e}^{a'}(e^{xL(1)}x^{-2L(0)}w_{a'}, x^{-1})\one\rangle 
\nn 
&=& \langle w_a, e^{x^{-1}L(-1)}e^{xL(1)}x^{-2L(0)}w_{a'}\rangle. 
\eea
On the other hand, we have 
\bea \label{3-equ-2-2}
\langle \Y_{a'a}^{e'}(w_{a'}, e^{\pi i}x)w_a, \one\rangle &=&
\langle \hat{A}_0(\Y_{a'e}^{a'})(w_{a'}, e^{\pi i}x)w_a, \one\rangle \nn
&=& \langle w_a, \Y_{a'e}^{a'}(e^{zL(1)}(e^{\pi i}z)^{-2L(0)}w_{a'}, x^{-1})
\one\rangle \nn
&=& \langle w_a, e^{x^{-1}L(-1)}e^{xL(1)}x^{-2L(0)}w_{a'}\rangle  
e^{-2\pi i h_a}.
\eea
Combining (\ref{3-equ-2-1}) and (\ref{3-equ-2-2}),
we obtain (\ref{3-equ-2}). 

The last identity (\ref{3-equ-3}) can be proved as
follow: 
\bea \label{3-equ-3-1}
\langle \tilde{A}_0(\Y_{ea'}^{a'})(\one, e^{\pi i}x)w_a, w_{a'}\rangle 
&=& \langle w_a, \Y_{ea'}^{a'}(e^{xL(1)}x^{-2L(0)}\one, x^{-1})w_{a'}\rangle
\nn
&=& \langle w_a, w_{a'}\rangle \nn
&=& \langle \Y_{ea}^a(\one, e^{\pi i}x)w_a, w_{a'}\rangle . \nonumber
\eea
\epf

The following formula is proved in \cite{FHL}.
\begin{lemma}
Let $\Y$ be an intertwining operator of type 
$\binom{W_3}{W_1W_2}$ and $v\in W_1$. We have
\beq \label{L-1-Y}
e^{xL(1)} \Y(v, x_0)e^{-xL(1)} = 
\Y(e^{x(1-xx_0)L(1)} (1-xx_0)^{-2L(0)}
v, x_0/(1-xx_0)). 
\eeq
\end{lemma}

The fusing matrices and their symmetries under the $S_3$ actions 
are studied in detail in \cite{H9}. Here similar
results hold. They are stated in the following Lemma:
\begin{lemma}
\bea
&&F(\Y_{a_1a_5;i}^{a_4;(1)}\otimes \Y_{a_2a_3;j}^{a_5;(2)}, 
\Y_{a_6a_3;k}^{a_4;(3)} \otimes \Y_{a_1a_2;l}^{a_6;(4)}) \nn
&&\hspace{1cm}= F^{-1}(\Omega_0(\Y_{a_1a_5;i}^{a_4;(1)})
\otimes \Omega_0(\Y_{a_2a_3;j}^{a_5;(2)}), 
\Omega_0(\Y_{a_6a_3;k}^{a_4;(3)}) \otimes \Omega_0(\Y_{a_1a_2;l}^{a_6;(4)})) 
\label{F-omega-inv} \\
&&\hspace{1cm}= F(\tilde{A}_0(\Y_{a_2a_3;j}^{a_5;(2)}) 
\otimes \tilde{A}_0(\Y_{a_1a_5;i}^{a_4;(1)}), 
\tilde{A}_0(\Y_{a_6a_3;k}^{a_4;(3)}) \otimes 
\Omega_{-1}(\Y_{a_1a_2;l}^{a_6;(4)}))
\label{F-AAA-omega-inv}
\eea
\end{lemma}
\pf
(\ref{F-omega-inv}) is proved in \cite{H9}. 
We prove (\ref{F-AAA-omega-inv}) here. 
Let $w_i\in W^{a_i}, i=1,2,3$ and
$w'_4\in (W^{a_4})'$.  For $|z_1|>|z_2|>|z_1-z_2|>0$, we have

\bea
&&\langle w'_4, \Y_{a_1a_5;i}^{a_4;(1)}(w_1, z_1)
\Y_{a_2a_3;j}^{a_5;(2)}(w_2, z_2)w_3\rangle \nn
&&\hspace{0.2cm}=
\langle \tilde{A}_0(\Y_{a_2a_3;j}^{a_5;(2)})(e^{-z_2L(1)}z_2^{-2L(0)}w_2, 
e^{\pi i} z_2^{-1})  \nn
&&\hspace{2cm} \tilde{A}_0(\Y_{a_1a_5;i}^{a_4;(1)})(e^{-z_1L(1)}z_1^{-2L(0)}w_1,
e^{\pi i} z_1^{-1})w'_4, w_3\rangle \nn
&&\hspace{0.2cm}=
\sum_{a_6\in \I} \sum_{k,l} 
F\big( \tilde{A}_0 ( 
\Y_{a_2a_3;j}^{a_5;(2)})\otimes \tilde{A}_0(\Y_{a_1a_5;i}^{a_4;(1)}),
\tilde{A}_0(\Y_{a_6a_3;k}^{a_4;(3)})\otimes 
\Omega_{-1}(\Y_{a_1a_2;l}^{a_6;(4)}) \big)  \nn
&&\hspace{1.5cm} \langle \tilde{A}_0(\Y_{a_6a_3;k}^{a_4;(3)})
\big( \Omega_{-1}(\Y_{a_1a_2;l}^{a_6;(4)})(e^{-z_2L(1)}z_2^{-2L(0)}w_2, 
e^{\pi i}(z_2^{-1}-z_1^{-1})) \nn
&&\hspace{3cm} \cdot \, e^{-z_1L(1)}z_1^{-2L(0)}w_1, 
e^{\pi i}z_1^{-1}\big) w'_4, w_3\rangle  \nn
&&\hspace{0.2cm}=\sum_{a_6\in \I} \sum_{k,l} 
F\big( \tilde{A}_0 ( 
\Y_{a_2a_3;j}^{a_5;(2)})\otimes \tilde{A}_0(\Y_{a_1a_5;i}^{a_4;(1)}),
\tilde{A}_0(\Y_{a_6a_3;k}^{a_4;(3)})\otimes 
\Omega_{-1}(\Y_{a_1a_2;l}^{a_6;(4)}) \big)  \nn
&&\hspace{1.5cm} \langle \tilde{A}_0(\Y_{a_6a_3;k}^{a_4;(3)})\big( 
\Y_{a_1a_2;l}^{a_6;(4)}(e^{-z_1L(1)}z_1^{-2L(0)}w_1, z_2^{-1}-z_1^{-1}) \nn
&&\hspace{3cm} \cdot \, 
e^{-z_2L(1)}z_2^{-2L(0)}w_2, e^{\pi i}z_2^{-1}\big) w'_4, w_3\rangle. 
\label{F-AAA-omega-equ-1}
\eea
Applying (\ref{L-1-Y}), we have the follow result:
\bea
&&\Y_{a_1a_2;l}^{a_6;(4)}
(e^{-z_1L(1)}z_1^{2L(0)}w_1, z_2^{-1}-z_1^{-1})e^{-z_2L(1)}z_2^{-2L(0)} \nn
&&\hspace{0.5cm}= e^{-z_2L(1)}\Y_{a_1a_2;l}^{a_6;(4)}
(e^{\frac{z_2^2}{z_1}L(1)}\left(\frac{z_2}{z_1}\right)^{-2L(0)}
e^{-z_1L(1)}z_1^{-2L(0)}w_1, \frac{z_1-z_2}{z_2^2})z_2^{-2L(0)}  \nn
&&\hspace{0.5cm}=e^{-z_2L(1)}z_2^{-2L(0)}
\Y_{a_1a_2;l}^{a_6;(4)}(e^{z^{-1}L(1)} z_1^{2L(0)}e^{-z_1L(1)}z_1^{-2L(0)}w_1, 
z_1-z_2)\nn
&&\hspace{0.5cm}=e^{-z_2L(1)}z_2^{-2L(0)}
\Y_{a_1a_2;l}^{a_6;(4)}(w_1, z_1-z_2). 
\eea
Applying it to the right hand side of 
(\ref{F-AAA-omega-equ-1}), we obtain
\bea
&&\langle w'_4, \Y_{a_1a_5;i}^{a_4;(1)}(w_1, z_1)
\Y_{a_2a_3;j}^{a_5;(2)}(w_2, z_2)w_3\rangle \nn
&&\hspace{0.2cm}=\sum_{a_6\in \I} \sum_{k,l} 
F\big( \tilde{A}_0( 
\Y_{a_2a_3;j}^{a_5;(2)})\otimes \tilde{A}_0(\Y_{a_1a_5;i}^{a_4;(1)}),
\tilde{A}_0(\Y_{a_6a_3;k}^{a_4;(3)})\otimes 
\Omega_{-1}(\Y_{a_1a_2;l}^{a_6;(4)}) \big)  \nn
&&\hspace{1.5cm} \langle \tilde{A}_0(\Y_{a_6a_3;k}^{a_4;(3)})
\big( e^{-z_2L(1)}z_2^{-2L(0)} \Y_{a_1a_2;l}^{a_6;(4)}(w, z_1-z_2)w_2, e^{\pi i}z_2
\big) w'_4, w_3\rangle \nn
&&\hspace{0.2cm}=\sum_{a_6\in \I} \sum_{k,l} 
F\big( \tilde{A}_0( 
\Y_{a_2a_3;j}^{a_5;(2)})\otimes \tilde{A}_0(\Y_{a_1a_5;i}^{a_4;(1)}),
\tilde{A}_0(\Y_{a_6a_3;k}^{a_4;(3)})\otimes 
\Omega_{-1}(\Y_{a_1a_2;l}^{a_6;(4)}) \big)  \nn
&&\hspace{1.5cm} \langle w'_4, 
\Y_{a_6a_3;k}^{a_4;(3)}(\Y_{a_1a_2;l}^{a_6;(4)}(w_1, z_1-z_2)w_2, z_2)w_3\rangle.
\eea
By the definition of fusing matrices, we see that
(\ref{F-AAA-omega-inv}) must be true. 
\epf

\begin{lemma}
\beq  \label{F-a-a'}
F_{a}=F_{a'} = F^{-1}(\Y_{ea}^{a}\otimes \Y_{aa'}^e; 
\Y_{ae}^{a}\otimes \Y_{a'a}^e)
=F^{-1}(\Y_{ea'}^{a'}\otimes \Y_{a'a}^e; 
\Y_{a'e}^{a'}\otimes \Y_{aa'}^e) .
\eeq
\end{lemma}
\pf
Using (\ref{F-omega-inv}) and Lemma \ref{3-equ-lemma}, 
we have
\bea
F^{-1}(\Y_{ea}^{a}\otimes \Y_{aa'}^e; \Y_{ae}^{a}\otimes \Y_{a'a}^e)
&=& F(\Y_{ae}^{a} \otimes \Omega_0(\Y_{aa'}^e), 
\Y_{ea}^{a} \otimes \Omega_0(\Y_{a'a}^e)) \nn
&=& F(\Y_{ae}^{a} \otimes e^{-2\pi ih_a} \Y_{a'a}^e, 
\Y_{ea}^{a} \otimes e^{-2\pi ih_a}\Y_{aa'}^e)  \nn
&=& F_a  \nonumber
\eea
Similarly we have
$F^{-1}(\Y_{ea'}^{a'}\otimes \Y_{a'a}^e; 
\Y_{a'e}^{a'}\otimes \Y_{aa'}^e) = F_{a'}$.
It remains to show that $F_a=F_{a'}$. Using 
(\ref{F-AAA-omega-inv}) and Lemma \ref{3-equ-lemma}, we have
\bea
F_{a'} &=& F(\Y_{a'e}^{a'} \otimes \Y_{aa'}^e, \Y_{ea'}^{a'} \otimes
\Y_{a'a}^e) \nn
&=& F(\tilde{A}_0(\Y_{aa'}^e) \otimes \tilde{A}_0(\Y_{a'e}^{a'}),
\tilde{A}_0(\Y_{ea'}^{a'}) \otimes \Omega_{-1}(\Y_{a'a}^e) ) \nn
&=& F(\Y_{ae}^a \otimes e^{2\pi ih_a} \Y_{a'a}^e, 
\Y_{ea}^a \otimes e^{2\pi ih_a}\Y_{aa'}^e) \nn
&=& F_a.  \nonumber
\eea
\epf

To prove the rigidity \cite{T}\cite{BK}, it amounts to 
prove the following Proposition \cite{H11}.
\begin{prop}  \label{rigidity-prop}
\bea   \label{R-dual-axiom-1}
I_{W^a} &=& r_{W^a} \circ I_{W^a} \boxtimes e_{a}
\circ (\A)^{-1} \circ i_{a}\boxtimes I_{W^a} \circ l_{W^a}^{-1},  
\\
I_{(W^a)'} &=& l_{(W^a)'} 
\circ e_{a}  \boxtimes I_{(W^a)'}
\circ \A \circ I_{W^a} \boxtimes i_{a}  \circ r_{(W^a)'}^{-1}, 
\label{R-dual-axiom-2} \\
I_{W^a} &=& l_{a} 
\circ e'_{a}  \boxtimes I_{W^a}
\circ \A \circ I_{W^a} \boxtimes i'_{a}  \circ r_{W^a}^{-1},  
\label{L-dual-axiom-1}  \\
I_{(W^a)'} &=& r_{(W^a)'} \circ I_{(W^a)'} \boxtimes e'_{a} 
\circ (\A)^{-1} \circ i'_{a} \boxtimes I_{W^a} \circ l_{(W^a)'}^{-1} .
\label{L-dual-axiom-2}
\eea
\end{prop}
\pf 
Let $z_1>z_2>z_1-z_2>0$.
By the universal property of tensor product, we have the following
canonical isomorphisms 
\bea  \label{nat-iso-1}
&&\V_{a_1a_5}^{a_4}  \cong 
\hom(W^{a_1}\boxtimes_{P(z_1)} W^{a_5}, W^{a_4}), \quad 
\V_{a_2a_3}^{a_5} \cong 
\hom_V(W^{a_2}\boxtimes_{P(z_2)} W^{a_3}, W^{a_5}),  \nn
&&\V_{a_6a_3}^{a_4} \cong
\hom(W^{a_6}\boxtimes_{P(z_2)} W^{a_3}, W^{a_4}), \quad
\V_{a_1a_2}^{a_6} \cong \hom_V(W^{a_1}\boxtimes_{P(z_1-z_2)} W^{a_2}, W^{a_6}).\nn
\eea 
We also have the following canonical isomorphisms 
\bea \label{nat-iso-2}
&&\coprod_{a_5\in \I} \hom_V(W^{a_1}\boxtimes_{P(z_1)} W^{a_5}, W^{a_4}) \otimes 
\hom_V(W^{a_2}\boxtimes_{P(z_2)} W^{a_3}, W^{a_5}) \nn
&&\hspace{3.5cm}\cong 
\hom(W^{a_1}\boxtimes_{P(z_1)} (W^{a_2}\boxtimes_{P(z_2)} W^{a_3}), W^{a_4})  \nn
&&\coprod_{a_5\in \I} \hom(W^{a_6}\boxtimes_{P(z_2)} 
W^{a_3}, W^{a_4}) \otimes 
\hom_V(W^{a_1}\boxtimes_{P(z_1-z_2)} W^{a_2}, W^{a_6}) \nn
&& \hspace{3.5cm} \cong 
\hom((W^{a_1}\boxtimes_{P(z_1-z_2)} W^{a_2})\boxtimes_{P(z_2)} W^{a_3}, W^{a_4})
\eea
defined by the compositions of maps. Moreover, let $\gamma_1$, $\gamma_2$
and $\gamma_3$ be paths in $\R_+$ from $z_1, z_2, z_1-z_2$ 
to $1$ respectively. 
By the naturalness of the parallel isomorphisms 
$\mathcal{T}_{\gamma_i}, i=1,2,3$ \cite{HL3}
\cite{HKL}\cite{HK1}\cite{H10}, we also have the following 
commutative diagram:
\beq  \label{nat-iso-3}
\xymatrix{
\hom(W^{a_1}\boxtimes_{P(z_1)} (W^{a_2}\boxtimes_{P(z_2)} W^{a_3}), W^{a_4})
\ar[r]^{\hspace{0.5cm}\cong} 
\ar[d]^{((\A_{P(z_1),P(z_2)}^{P(z_1-z_2),P(z_2)})^{-1})^*}  &  
\hom_{V}(W^{a_1}\boxtimes (W^{a_2}\boxtimes W^{a_3}), W^{a_4})  
\ar[d]^{(\A^{-1})^*} \\
\hom((W^{a_1}\boxtimes_{P(z_1-z_2)} W^{a_2})\boxtimes_{P(z_2)} W^{a_3}, W^{a_4})
\ar[r]^{\hspace{0.5cm}\cong}    &  
\hom_{V}((W^{a_1}\boxtimes W^{a_2})\boxtimes W^{a_3}, W^{a_4}) 
} 
\eeq
where the top and bottom horizontal isomorphisms are given by 
$\mathcal{T}_{\gamma_1}\circ (I_{W^{a_1}} \boxtimes_{P(z_1)} 
\mathcal{T}_{\gamma_2})$ 
and $\mathcal{T}_{\gamma_2}\circ  (\mathcal{T}_{\gamma_3} 
\boxtimes_{P(z_2)} I_{W^{a_3}})$ respectively. 

Combining (\ref{nat-iso-1}), (\ref{nat-iso-2}) and 
(\ref{nat-iso-3}), we obtain the following commutative
diagram
\beq  \label{nat-F-A}
\xymatrix{
\coprod_{a_5\in \I} \V_{a_1a_5}^{a_4}\otimes \V_{a_2a_3}^{a_5}
\ar[r]^{\hspace{-0.5cm}\cong} \ar[d]^{\mathcal{F}}  &  
\hom_{V}(W^{a_1}\boxtimes (W^{a_2}\boxtimes W^{a_3}), W^{a_4})  
\ar[d]^{(\A^{-1})*} \\
\coprod_{a_6\in \I} \V_{a_6a_3}^{a_4}\otimes \V_{a_1a_2}^{a_6} 
\ar[r]^{\hspace{-0.5cm} \cong}    &  
\hom_{V}((W^{a_1}\boxtimes W^{a_2})\boxtimes W^{a_3}, W^{a_4}) 
} 
\eeq
with two horizontal maps being canonical isomorphisms.
In terms of any basis of intertwining operators, 
the fusing isomorphism $\mathcal{F}$ can be written as:  
\bea
&&\mathcal{F}(\Y_{a_1a_5;i}^{a_4;(1)}\otimes \Y_{a_2a_3;j}^{a_5;(2)}) \nn
&&\hspace{1cm}=\sum_{a_6\in \I}\sum_{k,l} 
F(\Y_{a_1a_5;i}^{a_4;(1)}\otimes \Y_{a_2a_3;j}^{a_5;(2)}, 
\Y_{a_6a_3;k}^{a_4;(3)} \otimes \Y_{a_1a_2;l}^{a_6;(4)}) 
\Y_{a_6a_3;k}^{a_4;(3)} \otimes \Y_{a_1a_2;l}^{a_6;(4)}. 
\eea
By (\ref{nat-F-A}), we must have 
\bea  \label{A-rel-F-1}
&&\hspace{-1.2cm}
m_{\Y_{a_1a_5;i}^{a_4;(1)}} \circ 
(I_{W^{a_1}} \boxtimes m_{\Y_{a_2a_3;j}^{a_5;(2)}})  \circ \mathcal{A}^{-1} \nn
&&\hspace{-1cm}=\sum_{a_6\in \I}\sum_{k,l} 
F(\Y_{a_1a_5;i}^{a_4;(1)}\otimes \Y_{a_2a_3;j}^{a_5;(2)}, 
\Y_{a_6a_3;k}^{a_4;(3)} \otimes \Y_{a_1a_2;l}^{a_6;(4)}) 
m_{\Y_{a_6a_3;k}^{a_4;(3)}} \circ 
(m_{\Y_{a_1a_2;l}^{a_6;(4)}} \boxtimes I_{W^{a_3}}),
\eea
or equivalently
\bea  \label{A-rel-F-2}
&&\hspace{-1.3cm}
m_{\Y_{a_6a_3;k}^{a_4;(3)}} \circ 
(m_{\Y_{a_1a_2;l}^{a_6;(4)}} \boxtimes I_{W^{a_3}})\circ \mathcal{A} \nn
&&\hspace{-1.3cm}=\sum_{a_5\in \I}\sum_{i,j} 
F^{-1}(\Y_{a_6a_3;k}^{a_4;(3)} \otimes \Y_{a_1a_2;l}^{a_6;(4)};
\Y_{a_1a_5;i}^{a_4;(1)}\otimes \Y_{a_2a_3;j}^{a_5;(2)}) 
m_{\Y_{a_1a_5;i}^{a_4;(1)}} \circ 
(I_{W^{a_1}} \boxtimes m_{\Y_{a_2a_3;j}^{a_5;(2)}}).
\eea

Notice that $r_{W^a}=m_{\Y_{ae}^a}$, $e_a = m_{\Y_{a'a}^e}$ and 
$l_{W^a}=m_{\Y_{ea}^a}$ by our construction. Hence we have
\bea
&&r_{a} \circ I_{W^a} \boxtimes e_{a}
\circ (\A)^{-1} \circ i_{a}\boxtimes I_{W^a} \circ l_{a}^{-1}  \nn
&&\hspace{1cm} = 
m_{\Y_{ae}^a} \circ (I_{W^a} \boxtimes m_{\Y_{a'a}^e})
\circ (\A)^{-1} \circ i_{a}\boxtimes I_{W^a} \circ 
m_{\Y_{ea}^a}^{-1}  \nn
&&\hspace{-1.5cm}\mbox{\small by (\ref{A-rel-F-1})}
\hspace{0.7cm} = \sum_{b\in \I}\sum_{i,j} 
F(\Y_{ae}^a \otimes \Y_{a'a}^e; \Y_{ba;i}^{a;(1)}\otimes \Y_{aa';j}^{b;(2)}) \nn
&&\hspace{4cm}
m_{\Y_{ba;i}^{a;(1)}} \circ (m_{\Y_{aa';j}^{b;(2)}} \boxtimes I_{W^{a}}) 
\circ (i_{a}\boxtimes I_{W^a}) \circ m_{\Y_{ea}^a}^{-1}  \nn
&&\hspace{-1.5cm}\mbox{\small by (\ref{m-i-def-R})}
\hspace{0.7cm} = \frac{1}{F_a}
F(\Y_{ae}^a \otimes \Y_{a'a}^e; \Y_{ea}^a\otimes \Y_{aa'}^e)\, \, 
m_{\Y_{ea}^a} \circ m_{\Y_{ea}^a}^{-1}  \nn
&&\hspace{1cm} = I_{W^a}.  \nonumber
\eea
which is just (\ref{R-dual-axiom-1}). 
We prove the remaining three identities below. 
\bea  \label{proof-R-2}
&&l_{(W^a)'} \circ e_{a}  \boxtimes I_{(W^a)'}
\circ \A \circ I_{W^a} \boxtimes i_{a}  \circ r_{(W^a)'}^{-1} \nn
&&\hspace{1cm}= m_{\Y_{ea'}^{a'}} \circ 
(m_{\Y_{a'a}^e} \boxtimes I_{(W^a)'} ) \circ \A 
\circ I_{W^a} \boxtimes i_{a} \circ m_{\Y_{a'e}^{a'}}^{-1}  \nn
&&\hspace{-1.5cm}\mbox{\small by (\ref{A-rel-F-2})}
\hspace{0.7cm}= \sum_{b\in \I} \sum_{ij} 
F^{-1}(\Y_{ea'}^{a'} \otimes \Y_{a'a}^e ;
\Y_{a'b;i}^{a',(3)}\otimes \Y_{aa';j}^{b;(4)}) \nn 
&&\hspace{4cm}
m_{\Y_{a'b;i}^{a',(3)}} \circ 
(I_{W^a} \boxtimes m_{\Y_{aa';j}^{b;(4)})}) \circ 
(I_{W^a} \boxtimes i_{a})  \circ m_{\Y_{a'e}^{a'}}^{-1}  \nn
&&\hspace{-1.5cm}\mbox{\small by (\ref{m-i-def-R})}
\hspace{0.7cm} = \frac{1}{F_a} 
F^{-1}(\Y_{ea'}^{a'} \otimes \Y_{a'a}^e ;
\Y_{a'e}^{a'}\otimes \Y_{aa'}^{e}) \, m_{\Y_{a'e}^{a'}} \circ
m_{\Y_{a'e}^{a'}}^{-1}  \nn
&&\hspace{-1.5cm}\mbox{\small by (\ref{F-a-a'})}
\hspace{0.7cm} = I_{(W^a)'} 
\eea
and 
\bea  \label{proof-L-1}
&&l_{W^a} \circ e'_{a}  \boxtimes I_{W^a}
\circ \A \circ I_{W^a} \boxtimes i'_{a}  \circ r_{W^a}^{-1} \nn
&&\hspace{1cm}=  m_{\Y_{ea}^{a}} \circ
(m_{\Y_{aa'}^e} \boxtimes I_{W^a} ) \circ \A 
\circ I_{W^a} \boxtimes i'_{a} \circ m_{\Y_{ae}^e}^{-1}  \nn
&&\hspace{-1.5cm}\mbox{\small by (\ref{A-rel-F-2})}
\hspace{0.7cm}= \sum_{b\in \I} \sum_{ij} 
F^{-1}(\Y_{ea}^{a}\otimes \Y_{aa'}^e ;
\Y_{ab;i}^{a,(5)}\otimes \Y_{a'a;j}^{b;(6)}) \nn
&&\hspace{4cm}
m_{\Y_{ab;i}^{a,(5)}} \circ (I_{W^a} \boxtimes m_{\Y_{a'a;j}^{b;(6)}})
\circ (I_{W^a} \boxtimes i'_{a}) \circ m_{\Y_{ae}^a}^{-1}  \nn
&&\hspace{-1.5cm}\mbox{\small by (\ref{m-i-def-L})}
\hspace{0.7cm} =  \frac{1}{F_a} 
F^{-1}(\Y_{ea}^{a}\otimes \Y_{aa'}^e 
\Y_{ae}^{a}\otimes \Y_{a'a}^{e}) \, m_{\Y_{ae}^{a}} 
\circ m_{\Y_{ae}^a}^{-1}  \nn
, &&\hspace{-1.5cm}\mbox{\small by (\ref{F-a-a'})}
\hspace{0.7cm} = I_{W^a} 
\eea
and
\bea  \label{proof-L-2}
&&r_{(W^a)'} \circ I_{(W^a)'} \boxtimes e'_{a} 
\circ \A^{-1} \circ i'_{a} \boxtimes I_{W^a} \circ l_{(W^a)'}^{-1} 
\nn
&&\hspace{1cm}= m_{\Y_{a'e}^{a'}} \circ
(I_{(W^a)'} \boxtimes m_{\Y_{aa'}^e}) \circ \A^{-1} 
\circ (i'_{a} \boxtimes I_{W^a}) \circ m_{\Y_{ea'}^{a'}}^{-1} \nn
&&\hspace{-1.5cm}\mbox{\small by (\ref{A-rel-F-1})}
\hspace{0.7cm}= \sum_{b\in \I} \sum_{ij} 
F(\Y_{a'e}^{a'}\otimes \Y_{aa'}^e, 
\Y_{ba';i}^{a';(7)}\otimes \Y_{a'a;j}^{b;(8)}) \nn
&&\hspace{4cm}
m_{\Y_{ba';i}^{a';(7)}} \circ (m_{\Y_{a'a;j}^{b;(8)}} \boxtimes I_{W^a}) 
\circ  (i'_{a} \boxtimes I_{W^a}) \circ m_{\Y_{ea'}^{a'}}^{-1}  \nn
&&\hspace{-1.5cm}\mbox{\small by (\ref{m-i-def-L})}
\hspace{0.7cm} = \frac{1}{F_a}
F(\Y_{a'e}^{a'}\otimes \Y_{aa'}^e, 
\Y_{ea'}^{a'} \otimes \Y_{a'a}^{e}) \, m_{\Y_{ea'}^{a'}} 
\circ m_{\Y_{ea'}^{a'}}^{-1}\nn
&&\hspace{-1.5cm}\mbox{\small by (\ref{F-a-a'})}
\hspace{0.7cm} = I_{(W^a)'} 
\eea
\epf

\noindent {\small \sc 
Max Planck Institute for Mathematics in the Sciences, 
Inselstrasse 22, D-04103 Leipzig, Germany}

\vspace{1em}

\noindent {\small \sc 
Institut Des Hautes \'{E}tudes Scientifiques, 
Le Bois-Marie, 35, Route De Chartres,
F-91440 Bures-sur-Yvette, France} (current address)

\noindent {\em E-mail address}: kong@ihes.fr

\end{document}